
\catcode`!=11 
 
  

\def\PiC{P\kern-.12em\lower.5ex\hbox{I}\kern-.075emC}
\def\PiCTeX{\PiC\kern-.11em\TeX}

\def\!ifnextchar#1#2#3{%
  \let\!testchar=#1%
  \def\!first{#2}%
  \def\!second{#3}%
  \futurelet\!nextchar\!testnext}
\def\!testnext{%
  \ifx \!nextchar \!spacetoken 
    \let\!next=\!skipspacetestagain
  \else
    \ifx \!nextchar \!testchar
      \let\!next=\!first
    \else 
      \let\!next=\!second 
    \fi 
  \fi
  \!next}
\def\\{\!skipspacetestagain} 
  \expandafter\def\\ {\futurelet\!nextchar\!testnext} 
\def\\{\let\!spacetoken= } \\  

\def\!tfor#1:=#2\do#3{%
  \edef\!fortemp{#2}%
  \ifx\!fortemp\!empty 
    \else
    \!tforloop#2\!nil\!nil\!!#1{#3}%
  \fi}
\def\!tforloop#1#2\!!#3#4{%
  \def#3{#1}%
  \ifx #3\!nnil
    \let\!nextwhile=\!fornoop
  \else
    #4\relax
    \let\!nextwhile=\!tforloop
  \fi 
  \!nextwhile#2\!!#3{#4}}

\def\!etfor#1:=#2\do#3{%
  \def\!!tfor{\!tfor#1:=}%
  \edef\!!!tfor{#2}%
  \expandafter\!!tfor\!!!tfor\do{#3}}

\def\!cfor#1:=#2\do#3{%
  \edef\!fortemp{#2}%
  \ifx\!fortemp\!empty 
  \else
    \!cforloop#2,\!nil,\!nil\!!#1{#3}%
  \fi}
\def\!cforloop#1,#2\!!#3#4{%
  \def#3{#1}%
  \ifx #3\!nnil
    \let\!nextwhile=\!fornoop 
  \else
    #4\relax
    \let\!nextwhile=\!cforloop
  \fi
  \!nextwhile#2\!!#3{#4}}

\def\!ecfor#1:=#2\do#3{%
  \def\!!cfor{\!cfor#1:=}%
  \edef\!!!cfor{#2}%
  \expandafter\!!cfor\!!!cfor\do{#3}}

\def\!empty{}
\def\!nnil{\!nil}
\def\!fornoop#1\!!#2#3{}

\def\!ifempty#1#2#3{%
  \edef\!emptyarg{#1}%
  \ifx\!emptyarg\!empty
    #2%
  \else
    #3%
  \fi}
 
\def\!getnext#1\from#2{%
  \expandafter\!gnext#2\!#1#2}%
\def\!gnext\\#1#2\!#3#4{%
  \def#3{#1}%
  \def#4{#2\\{#1}}%
  \ignorespaces}

%
\def\!getnextvalueof#1\from#2{%
  \expandafter\!gnextv#2\!#1#2}%
\def\!gnextv\\#1#2\!#3#4{%
  #3=#1%
  \def#4{#2\\{#1}}%
  \ignorespaces}

\def\!copylist#1\to#2{%
  \expandafter\!!copylist#1\!#2}
\def\!!copylist#1\!#2{%
  \def#2{#1}\ignorespaces}

\def\!wlet#1=#2{%
  \let#1=#2 
  \wlog{\string#1=\string#2}}
 
\def\!listaddon#1#2{%
  \expandafter\!!listaddon#2\!{#1}#2}
\def\!!listaddon#1\!#2#3{%
  \def#3{#1\\#2}}
 

\def\!rightappend#1\withCS#2\to#3{\expandafter\!!rightappend#3\!#2{#1}#3}
\def\!!rightappend#1\!#2#3#4{\def#4{#1#2{#3}}}

\def\!leftappend#1\withCS#2\to#3{\expandafter\!!leftappend#3\!#2{#1}#3}
\def\!!leftappend#1\!#2#3#4{\def#4{#2{#3}#1}}

\def\!lop#1\to#2{\expandafter\!!lop#1\!#1#2}
\def\!!lop\\#1#2\!#3#4{\def#4{#1}\def#3{#2}}



\def\!loop#1\repeat{\def\!body{#1}\!iterate}
\def\!iterate{\!body\let\!next=\!iterate\else\let\!next=\relax\fi\!next}
 
\def\!!loop#1\repeat{\def\!!body{#1}\!!iterate}
\def\!!iterate{\!!body\let\!!next=\!!iterate\else\let\!!next=\relax\fi\!!next}
 
\def\!removept#1#2{\edef#2{\expandafter\!!removePT\the#1}}
{\catcode`p=12 \catcode`t=12 \gdef\!!removePT#1pt{#1}}

\def\placevalueinpts of <#1> in #2 {%
  \!removept{#1}{#2}}
 
\def\!mlap#1{\hbox to 0pt{\hss#1\hss}}
\def\!vmlap#1{\vbox to 0pt{\vss#1\vss}}
 
\def\!not#1{%
  #1\relax
    \!switchfalse
  \else
    \!switchtrue
  \fi
  \if!switch
  \ignorespaces}


 

\let\!!!wlog=\wlog              
\def\wlog#1{}    

\newdimen\headingtoplotskip     
\newdimen\linethickness         
\newdimen\longticklength        
\newdimen\plotsymbolspacing     
\newdimen\shortticklength       
\newdimen\stackleading          
\newdimen\tickstovaluesleading  
\newdimen\totalarclength        
\newdimen\valuestolabelleading  

\newbox\!boxA                   
\newbox\!boxB                   
\newbox\!picbox                 
\newbox\!plotsymbol             
\newbox\!putobject              
\newbox\!shadesymbol            

\newcount\!countA               
\newcount\!countB               
\newcount\!countC               
\newcount\!countD               
\newcount\!countE               
\newcount\!countF               
\newcount\!countG               
\newcount\!fiftypt              
\newcount\!intervalno           
\newcount\!npoints              
\newcount\!nsegments            
\newcount\!ntemp                
\newcount\!parity               
\newcount\!scalefactor          
\newcount\!tfs                  
\newcount\!tickcase             

\newdimen\!Xleft                
\newdimen\!Xright               
\newdimen\!Xsave                
\newdimen\!Ybot                 
\newdimen\!Ysave                
\newdimen\!Ytop                 
\newdimen\!angle                
\newdimen\!arclength            
\newdimen\!areabloc             
\newdimen\!arealloc             
\newdimen\!arearloc             
\newdimen\!areatloc             
\newdimen\!bshrinkage           
\newdimen\!checkbot             
\newdimen\!checkleft            
\newdimen\!checkright           
\newdimen\!checktop             
\newdimen\!dimenA               
\newdimen\!dimenB               
\newdimen\!dimenC               
\newdimen\!dimenD               
\newdimen\!dimenE               
\newdimen\!dimenF               
\newdimen\!dimenG               
\newdimen\!dimenH               
\newdimen\!dimenI               
\newdimen\!distacross           
\newdimen\!downlength           
\newdimen\!dp                   
\newdimen\!dshade               
\newdimen\!dxpos                
\newdimen\!dxprime              
\newdimen\!dypos                
\newdimen\!dyprime              
\newdimen\!ht                   
\newdimen\!leaderlength         
\newdimen\!lshrinkage           
\newdimen\!midarclength         
\newdimen\!offset               
\newdimen\!plotheadingoffset    
\newdimen\!plotsymbolxshift     
\newdimen\!plotsymbolyshift     
\newdimen\!plotxorigin          
\newdimen\!plotyorigin          
\newdimen\!rootten              
\newdimen\!rshrinkage           
\newdimen\!shadesymbolxshift    
\newdimen\!shadesymbolyshift    
\newdimen\!tenAa                
\newdimen\!tenAc                
\newdimen\!tenAe                
\newdimen\!tshrinkage           
\newdimen\!uplength             
\newdimen\!wd                   
\newdimen\!wmax                 
\newdimen\!wmin                 
\newdimen\!xB                   
\newdimen\!xC                   
\newdimen\!xE                   
\newdimen\!xM                   
\newdimen\!xS                   
\newdimen\!xaxislength          
\newdimen\!xdiff                
\newdimen\!xleft                
\newdimen\!xloc                 
\newdimen\!xorigin              
\newdimen\!xpivot               
\newdimen\!xpos                 
\newdimen\!xprime               
\newdimen\!xright               
\newdimen\!xshade               
\newdimen\!xshift               
\newdimen\!xtemp                
\newdimen\!xunit                
\newdimen\!xxE                  
\newdimen\!xxM                  
\newdimen\!xxS                  
\newdimen\!xxloc                
\newdimen\!yB                   
\newdimen\!yC                   
\newdimen\!yE                   
\newdimen\!yM                   
\newdimen\!yS                   
\newdimen\!yaxislength          
\newdimen\!ybot                 
\newdimen\!ydiff                
\newdimen\!yloc                 
\newdimen\!yorigin              
\newdimen\!ypivot               
\newdimen\!ypos                 
\newdimen\!yprime               
\newdimen\!yshade               
\newdimen\!yshift               
\newdimen\!ytemp                
\newdimen\!ytop                 
\newdimen\!yunit                
\newdimen\!yyE                  
\newdimen\!yyM                  
\newdimen\!yyS                  
\newdimen\!yyloc                
\newdimen\!zpt                  

\newif\if!axisvisible           
\newif\if!gridlinestoo          
\newif\if!keepPO                
\newif\if!placeaxislabel        
\newif\if!switch                
\newif\if!xswitch               

\newtoks\!axisLaBeL             
\newtoks\!keywordtoks           

\newwrite\!replotfile           

\newhelp\!keywordhelp{The keyword mentioned in the error message in unknown. 
Replace NEW KEYWORD in the indicated response by the keyword that 
should have been specified.}    

\!wlet\!!origin=\!xM                   
\!wlet\!!unit=\!uplength               
\!wlet\!Lresiduallength=\!dimenG       
\!wlet\!Rresiduallength=\!dimenF       
\!wlet\!axisLength=\!distacross        
\!wlet\!axisend=\!ydiff                
\!wlet\!axisstart=\!xdiff              
\!wlet\!axisxlevel=\!arclength         
\!wlet\!axisylevel=\!downlength        
\!wlet\!beta=\!dimenE                  
\!wlet\!gamma=\!dimenF                 
\!wlet\!shadexorigin=\!plotxorigin     
\!wlet\!shadeyorigin=\!plotyorigin     
\!wlet\!ticklength=\!xS                
\!wlet\!ticklocation=\!xE              
\!wlet\!ticklocationincr=\!yE          
\!wlet\!tickwidth=\!yS                 
\!wlet\!totalleaderlength=\!dimenE     
\!wlet\!xone=\!xprime                  
\!wlet\!xtwo=\!dxprime                 
\!wlet\!ySsave=\!yM                    
\!wlet\!ybB=\!yB                       
\!wlet\!ybC=\!yC                       
\!wlet\!ybE=\!yE                       
\!wlet\!ybM=\!yM                       
\!wlet\!ybS=\!yS                       
\!wlet\!ybpos=\!yyloc                  
\!wlet\!yone=\!yprime                  
\!wlet\!ytB=\!xB                       
\!wlet\!ytC=\!xC                       
\!wlet\!ytE=\!downlength               
\!wlet\!ytM=\!arclength                
\!wlet\!ytS=\!distacross               
\!wlet\!ytpos=\!xxloc                  
\!wlet\!ytwo=\!dyprime                 

\!zpt=0pt                              
\!xunit=1pt
\!yunit=1pt
\!arearloc=\!xunit
\!areatloc=\!yunit
\!dshade=5pt
\!leaderlength=24in
\!tfs=256                              
\!wmax=5.3pt                           
\!wmin=2.7pt                           
\!xaxislength=\!xunit
\!xpivot=\!zpt
\!yaxislength=\!yunit 
\!ypivot=\!zpt
\plotsymbolspacing=.4pt
  \!dimenA=50pt \!fiftypt=\!dimenA     

\!rootten=3.162278pt                   
\!tenAa=8.690286pt                     
\!tenAc=2.773839pt                     
\!tenAe=2.543275pt                     

\def\!cosrotationangle{1}      
\def\!sinrotationangle{0}      
\def\!xpivotcoord{0}           
\def\!xref{0}                  
\def\!xshadesave{0}            
\def\!ypivotcoord{0}           
\def\!yref{0}                  
\def\!yshadesave{0}            
\def\!zero{0}                  

\let\wlog=\!!!wlog
%
  
\def\normalgraphs{%
  \longticklength=.4\baselineskip
  \shortticklength=.25\baselineskip
  \tickstovaluesleading=.25\baselineskip
  \valuestolabelleading=.8\baselineskip
  \linethickness=.4pt
  \stackleading=.17\baselineskip
  \headingtoplotskip=1.5\baselineskip
  \visibleaxes
  \ticksout
  \nogridlines
  \unloggedticks}
%
\def\setplotarea x from #1 to #2, y from #3 to #4 {%
  \!arealloc=\!M{#1}\!xunit \advance \!arealloc -\!xorigin
  \!areabloc=\!M{#3}\!yunit \advance \!areabloc -\!yorigin
  \!arearloc=\!M{#2}\!xunit \advance \!arearloc -\!xorigin
  \!areatloc=\!M{#4}\!yunit \advance \!areatloc -\!yorigin
  \!initinboundscheck
  \!xaxislength=\!arearloc  \advance\!xaxislength -\!arealloc
  \!yaxislength=\!areatloc  \advance\!yaxislength -\!areabloc
  \!plotheadingoffset=\!zpt
  \!dimenput {{\setbox0=\hbox{}\wd0=\!xaxislength\ht0=\!yaxislength\box0}}
     [bl] (\!arealloc,\!areabloc)}
%
\def\visibleaxes{%
  \def\!axisvisibility{\!axisvisibletrue}}

%

\def\!fixkeyword#1{%
  \errhelp=\!keywordhelp
  \errmessage{Unrecognized keyword `#1': \the\!keywordtoks{NEW KEYWORD}'}}

\!keywordtoks={enter `i\fixkeyword}

\def\fixkeyword#1{%
  \!nextkeyword#1 }


\def\axis {%
  \def\!nextkeyword##1 {%
    \expandafter\ifx\csname !axis##1\endcsname \relax
      \def\!next{\!fixkeyword{##1}}%
    \else
      \def\!next{\csname !axis##1\endcsname}%
    \fi
    \!next}%
  \!offset=\!zpt
  \!axisvisibility
  \!placeaxislabelfalse
  \!nextkeyword}

\def\!axisbottom{%
  \!axisylevel=\!areabloc
  \def\!tickxsign{0}%
  \def\!tickysign{-}%
  \def\!axissetup{\!axisxsetup}%
  \def\!axislabeltbrl{t}%
  \!nextkeyword}

\def\!axistop{%
  \!axisylevel=\!areatloc
  \def\!tickxsign{0}%
  \def\!tickysign{+}%
  \def\!axissetup{\!axisxsetup}%
  \def\!axislabeltbrl{b}%
  \!nextkeyword}

\def\!axisleft{%
  \!axisxlevel=\!arealloc
  \def\!tickxsign{-}%
  \def\!tickysign{0}%
  \def\!axissetup{\!axisysetup}%
  \def\!axislabeltbrl{r}%
  \!nextkeyword}

\def\!axisright{%
  \!axisxlevel=\!arearloc
  \def\!tickxsign{+}%
  \def\!tickysign{0}%
  \def\!axissetup{\!axisysetup}%
  \def\!axislabeltbrl{l}%
  \!nextkeyword}

\def\!axisshiftedto#1=#2 {%
  \if 0\!tickxsign
    \!axisylevel=\!M{#2}\!yunit
    \advance\!axisylevel -\!yorigin
  \else
    \!axisxlevel=\!M{#2}\!xunit
    \advance\!axisxlevel -\!xorigin
  \fi
  \!nextkeyword}

\def\!axisvisible{%
  \!axisvisibletrue  
  \!nextkeyword}

\def\!axisinvisible{%
  \!axisvisiblefalse
  \!nextkeyword}

\def\!axislabel#1 {%
  \!axisLaBeL={#1}%
  \!placeaxislabeltrue
  \!nextkeyword}

\expandafter\def\csname !axis/\endcsname{%
  \!axissetup 
  \if!placeaxislabel
    \!placeaxislabel
  \fi
  \if +\!tickysign 
    \!dimenA=\!axisylevel
    \advance\!dimenA \!offset 
    \advance\!dimenA -\!areatloc 
    \ifdim \!dimenA>\!plotheadingoffset
      \!plotheadingoffset=\!dimenA 
    \fi
  \fi}

\def\grid #1 #2 {%
  \!countA=#1\advance\!countA 1
  \axis bottom invisible ticks length <\!zpt> andacross quantity {\!countA} /
  \!countA=#2\advance\!countA 1
  \axis left   invisible ticks length <\!zpt> andacross quantity {\!countA} / }

\def\plotheading#1 {%
  \advance\!plotheadingoffset \headingtoplotskip
  \!dimenput {#1} [B] <.5\!xaxislength,\!plotheadingoffset>
    (\!arealloc,\!areatloc)}

\def\!axisxsetup{%
  \!axisxlevel=\!arealloc
  \!axisstart=\!arealloc
  \!axisend=\!arearloc
  \!axisLength=\!xaxislength
  \!!origin=\!xorigin
  \!!unit=\!xunit
  \!xswitchtrue
  \if!axisvisible 
    \!makeaxis
  \fi}

\def\!axisysetup{%
  \!axisylevel=\!areabloc
  \!axisstart=\!areabloc
  \!axisend=\!areatloc
  \!axisLength=\!yaxislength
  \!!origin=\!yorigin
  \!!unit=\!yunit
  \!xswitchfalse
  \if!axisvisible
    \!makeaxis
  \fi}

\def\!makeaxis{%
  \setbox\!boxA=\hbox{
    \beginpicture
      \!setdimenmode
      \setcoordinatesystem point at {\!zpt} {\!zpt}   
      \putrule from {\!zpt} {\!zpt} to
        {\!tickysign\!tickysign\!axisLength} 
        {\!tickxsign\!tickxsign\!axisLength}
    \endpicturesave <\!Xsave,\!Ysave>}%
    \wd\!boxA=\!zpt
    \!placetick\!axisstart}

\def\!placeaxislabel{%
  \advance\!offset \valuestolabelleading
  \if!xswitch
    \!dimenput {\the\!axisLaBeL} [\!axislabeltbrl]
      <.5\!axisLength,\!tickysign\!offset> (\!axisxlevel,\!axisylevel)
    \advance\!offset \!dp  
    \advance\!offset \!ht  
  \else
    \!dimenput {\the\!axisLaBeL} [\!axislabeltbrl]
      <\!tickxsign\!offset,.5\!axisLength> (\!axisxlevel,\!axisylevel)
  \fi
  \!axisLaBeL={}}

%


\def\arrow <#1> [#2,#3]{%
  \!ifnextchar<{\!arrow{#1}{#2}{#3}}{\!arrow{#1}{#2}{#3}<\!zpt,\!zpt> }}

\def\!arrow#1#2#3<#4,#5> from #6 #7 to #8 #9 {%
%
  \!xloc=\!M{#8}\!xunit   
  \!yloc=\!M{#9}\!yunit
  \!dxpos=\!xloc  \!dimenA=\!M{#6}\!xunit  \advance \!dxpos -\!dimenA
  \!dypos=\!yloc  \!dimenA=\!M{#7}\!yunit  \advance \!dypos -\!dimenA
  \let\!MAH=\!M
  \!setdimenmode
  \!xshift=#4\relax  \!yshift=#5\relax
  \!reverserotateonly\!xshift\!yshift
  \advance\!xshift\!xloc  \advance\!yshift\!yloc
%
  \!xS=-\!dxpos  \advance\!xS\!xshift
  \!yS=-\!dypos  \advance\!yS\!yshift
  \!start (\!xS,\!yS)
  \!ljoin (\!xshift,\!yshift)
%
  \!Pythag\!dxpos\!dypos\!arclength
  \!divide\!dxpos\!arclength\!dxpos  
  \!dxpos=32\!dxpos  \!removept\!dxpos\!!cos
  \!divide\!dypos\!arclength\!dypos  
  \!dypos=32\!dypos  \!removept\!dypos\!!sin
%
  \!halfhead{#1}{#2}{#3}
  \!halfhead{#1}{-#2}{-#3}
  \let\!M=\!MAH
  \ignorespaces}
%
  \def\!halfhead#1#2#3{%
    \!dimenC=-#1%
    \divide \!dimenC 2 
    \!dimenD=#2\!dimenC
    \!rotate(\!dimenC,\!dimenD)by(\!!cos,\!!sin)to(\!xM,\!yM)
    \!dimenC=-#1
    \!dimenD=#3\!dimenC
    \!dimenD=.5\!dimenD
    \!rotate(\!dimenC,\!dimenD)by(\!!cos,\!!sin)to(\!xE,\!yE)
    \!start (\!xshift,\!yshift)
    \advance\!xM\!xshift  \advance\!yM\!yshift
    \advance\!xE\!xshift  \advance\!yE\!yshift
    \!qjoin (\!xM,\!yM) (\!xE,\!yE) 
    \ignorespaces}

\def\betweenarrows #1#2 from #3 #4 to #5 #6 {%
  \!xloc=\!M{#3}\!xunit  \!xxloc=\!M{#5}\!xunit%
  \!yloc=\!M{#4}\!yunit  \!yyloc=\!M{#6}\!yunit%
  \!dxpos=\!xxloc  \advance\!dxpos by -\!xloc
  \!dypos=\!yyloc  \advance\!dypos by -\!yloc
  \advance\!xloc .5\!dxpos
  \advance\!yloc .5\!dypos
  \let\!MBA=\!M
  \!setdimenmode
  \ifdim\!dypos=\!zpt
    \ifdim\!dxpos<\!zpt \!dxpos=-\!dxpos \fi
    \put {\!lrarrows{\!dxpos}{#1}}#2{} at {\!xloc} {\!yloc}
  \else
    \ifdim\!dxpos=\!zpt
      \ifdim\!dypos<\!zpt \!dypos=-\!dypos \fi
      \put {\!udarrows{\!dypos}{#1}}#2{} at {\!xloc} {\!yloc}
    \fi
  \fi
  \let\!M=\!MBA
  \ignorespaces}

\def\!lrarrows#1#2{
  {\setbox\!boxA=\hbox{$\mkern-2mu\mathord-\mkern-2mu$}%
   \setbox\!boxB=\hbox{$\leftarrow$}\!dimenE=\ht\!boxB
   \setbox\!boxB=\hbox{}\ht\!boxB=2\!dimenE
   \hbox to #1{$\mathord\leftarrow\mkern-6mu
     \cleaders\copy\!boxA\hfil
     \mkern-6mu\mathord-$%
     \kern.4em $\vcenter{\box\!boxB}$$\vcenter{\hbox{#2}}$\kern.4em
     $\mathord-\mkern-6mu
     \cleaders\copy\!boxA\hfil
     \mkern-6mu\mathord\rightarrow$}}}

\def\!udarrows#1#2{
  {\setbox\!boxB=\hbox{#2}%
   \setbox\!boxA=\hbox to \wd\!boxB{\hss$\vert$\hss}%
   \!dimenE=\ht\!boxA \advance\!dimenE \dp\!boxA \divide\!dimenE 2
   \vbox to #1{\offinterlineskip
      \vskip .05556\!dimenE
      \hbox to \wd\!boxB{\hss$\mkern.4mu\uparrow$\hss}\vskip-\!dimenE
      \cleaders\copy\!boxA\vfil
      \vskip-\!dimenE\copy\!boxA
      \vskip\!dimenE\copy\!boxB\vskip.4em
      \copy\!boxA\vskip-\!dimenE
      \cleaders\copy\!boxA\vfil
      \vskip-\!dimenE \hbox to \wd\!boxB{\hss$\mkern.4mu\downarrow$\hss}
      \vskip .05556\!dimenE}}}

%

\def\putbar#1breadth <#2> from #3 #4 to #5 #6 {%
  \!xloc=\!M{#3}\!xunit  \!xxloc=\!M{#5}\!xunit%
  \!yloc=\!M{#4}\!yunit  \!yyloc=\!M{#6}\!yunit%
  \!dypos=\!yyloc  \advance\!dypos by -\!yloc
  \!dimenI=#2  
  \ifdim \!dimenI=\!zpt 
    \putrule#1from {#3} {#4} to {#5} {#6} 
  \else 
    \let\!MBar=\!M
    \!setdimenmode 
    \divide\!dimenI 2
    \ifdim \!dypos=\!zpt             
      \advance \!yloc -\!dimenI 
      \advance \!yyloc \!dimenI
    \else
      \advance \!xloc -\!dimenI 
      \advance \!xxloc \!dimenI
    \fi
    \putrectangle#1corners at {\!xloc} {\!yloc} and {\!xxloc} {\!yyloc}
    \let\!M=\!MBar 
  \fi
  \ignorespaces}

\def\setbars#1breadth <#2> baseline at #3 = #4 {%
  \edef\!barshift{#1}%
  \edef\!barbreadth{#2}%
  \edef\!barorientation{#3}%
  \edef\!barbaseline{#4}%
  \def\!bardobaselabel{\!bardoendlabel}%
  \def\!bardoendlabel{\!barfinish}%
  \let\!drawcurve=\!barcurve
  \!setbars}
\def\!setbars{%
  \futurelet\!nextchar\!!setbars}
\def\!!setbars{%
  \if b\!nextchar
    \def\!!!setbars{\!setbarsbget}%
  \else 
    \if e\!nextchar
      \def\!!!setbars{\!setbarseget}%
    \else
      \def\!!!setbars{\relax}%
    \fi
  \fi
  \!!!setbars}
\def\!setbarsbget baselabels (#1) {%
  \def\!barbaselabelorientation{#1}%
  \def\!bardobaselabel{\!!bardobaselabel}%
  \!setbars}
\def\!setbarseget endlabels (#1) {%
  \edef\!barendlabelorientation{#1}%
  \def\!bardoendlabel{\!!bardoendlabel}%
  \!setbars}

\def\!barcurve #1 #2 {%
  \if y\!barorientation
    \def\!basexarg{#1}%
    \def\!baseyarg{\!barbaseline}%
  \else
    \def\!basexarg{\!barbaseline}%
    \def\!baseyarg{#2}%
  \fi
  \expandafter\putbar\!barshift breadth <\!barbreadth> from {\!basexarg}
    {\!baseyarg} to {#1} {#2}
  \def\!endxarg{#1}%
  \def\!endyarg{#2}%
  \!bardobaselabel}

\def\!!bardobaselabel "#1" {%
  \put {#1}\!barbaselabelorientation{} at {\!basexarg} {\!baseyarg}
  \!bardoendlabel}
 
\def\!!bardoendlabel "#1" {%
  \put {#1}\!barendlabelorientation{} at {\!endxarg} {\!endyarg}
  \!barfinish}

\def\!barfinish{%
  \!ifnextchar/{\!finish}{\!barcurve}}

%
%
%
\def\putrectangle{%
  \!ifnextchar<{\!putrectangle}{\!putrectangle<\!zpt,\!zpt> }}
\def\!putrectangle<#1,#2> corners at #3 #4 and #5 #6 {%
%
  \!xone=\!M{#3}\!xunit  \!xtwo=\!M{#5}\!xunit%
  \!yone=\!M{#4}\!yunit  \!ytwo=\!M{#6}\!yunit%
  \ifdim \!xtwo<\!xone
    \!dimenI=\!xone  \!xone=\!xtwo  \!xtwo=\!dimenI
  \fi
  \ifdim \!ytwo<\!yone
    \!dimenI=\!yone  \!yone=\!ytwo  \!ytwo=\!dimenI
  \fi
  \!dimenI=#1\relax  \advance\!xone\!dimenI  \advance\!xtwo\!dimenI
  \!dimenI=#2\relax  \advance\!yone\!dimenI  \advance\!ytwo\!dimenI
  \let\!MRect=\!M
  \!setdimenmode
%
  \!shaderectangle
%
  \!dimenI=.5\linethickness
  \advance \!xone  -\!dimenI
  \advance \!xtwo   \!dimenI
  \putrule from {\!xone} {\!yone} to {\!xtwo} {\!yone} 
  \putrule from {\!xone} {\!ytwo} to {\!xtwo} {\!ytwo} 
%
  \advance \!xone   \!dimenI
  \advance \!xtwo  -\!dimenI%
  \advance \!yone  -\!dimenI
  \advance \!ytwo   \!dimenI
  \putrule from {\!xone} {\!yone} to {\!xone} {\!ytwo} 
  \putrule from {\!xtwo} {\!yone} to {\!xtwo} {\!ytwo} 
  \let\!M=\!MRect
  \ignorespaces}
 

\def\shaderectanglesoff{%
  \def\!shaderectangle{}%
  \ignorespaces}

\shaderectanglesoff
 
\def\!!shaderectangle{%
  \!dimenA=\!xtwo  \advance \!dimenA -\!xone
  \!dimenB=\!ytwo  \advance \!dimenB -\!yone
  \ifdim \!dimenA<\!dimenB
    \!startvshade (\!xone,\!yone,\!ytwo)
    \!lshade      (\!xtwo,\!yone,\!ytwo)
  \else
    \!starthshade (\!yone,\!xone,\!xtwo)
    \!lshade      (\!ytwo,\!xone,\!xtwo)
  \fi
  \ignorespaces}
  
\def\frame{%
  \!ifnextchar<{\!frame}{\!frame<\!zpt> }}
\long\def\!frame<#1> #2{%
  \beginpicture
    \setcoordinatesystem units <1pt,1pt> point at 0 0 
    \put {#2} [Bl] at 0 0 
    \!dimenA=#1\relax
    \!dimenB=\!wd \advance \!dimenB \!dimenA
    \!dimenC=\!ht \advance \!dimenC \!dimenA
    \!dimenD=\!dp \advance \!dimenD \!dimenA
    \let\!MFr=\!M
    \!setdimenmode
    \putrectangle corners at {-\!dimenA} {-\!dimenD} and {\!dimenB} {\!dimenC}
    \!setcoordmode
    \let\!M=\!MFr
  \endpicture
  \ignorespaces}
 
\def\rectangle <#1> <#2> {%
  \setbox0=\hbox{}\wd0=#1\ht0=#2\frame {\box0}}

%


\def\!plotfromfile"#1"{%
  \expandafter\!drawcurve \input #1 /}

\def\setquadratic{%
  \let\!drawcurve=\!qcurve
  \let\!!Shade=\!!qShade
  \let\!!!Shade=\!!!qShade}

\def\setlinear{%
  \let\!drawcurve=\!lcurve
  \let\!!Shade=\!!lShade
  \let\!!!Shade=\!!!lShade}

\def\sethistograms{%
  \let\!drawcurve=\!hcurve}

\def\!qcurve #1 #2 {%
  \!start (#1,#2)
  \!Qjoin}
\def\!Qjoin#1 #2 #3 #4 {%
  \!qjoin (#1,#2) (#3,#4)             
  \!ifnextchar/{\!finish}{\!Qjoin}}

\def\!lcurve #1 #2 {%
  \!start (#1,#2)
  \!Ljoin}
\def\!Ljoin#1 #2 {%
  \!ljoin (#1,#2)                    
  \!ifnextchar/{\!finish}{\!Ljoin}}

\def\!finish/{\ignorespaces}

\def\!hcurve #1 #2 {%
  \edef\!hxS{#1}%
  \edef\!hyS{#2}%
  \!hjoin}
\def\!hjoin#1 #2 {%
  \putrectangle corners at {\!hxS} {\!hyS} and {#1} {#2}
  \edef\!hxS{#1}%
  \!ifnextchar/{\!finish}{\!hjoin}}

\def\vshade #1 #2 #3 {%
  \!startvshade (#1,#2,#3)
  \!Shadewhat}

\def\hshade #1 #2 #3 {%
  \!starthshade (#1,#2,#3)
  \!Shadewhat}

\def\!Shadewhat{%
  \futurelet\!nextchar\!Shade}
\def\!Shade{%
  \if <\!nextchar
    \def\!nextShade{\!!Shade}%
  \else
    \if /\!nextchar
      \def\!nextShade{\!finish}%
    \else
      \def\!nextShade{\!!!Shade}%
    \fi
  \fi
  \!nextShade}
\def\!!lShade<#1> #2 #3 #4 {%
  \!lshade <#1> (#2,#3,#4)                 
  \!Shadewhat}
\def\!!!lShade#1 #2 #3 {%
  \!lshade (#1,#2,#3)
  \!Shadewhat} 
\def\!!qShade<#1> #2 #3 #4 #5 #6 #7 {%
  \!qshade <#1> (#2,#3,#4) (#5,#6,#7)      
  \!Shadewhat}
\def\!!!qShade#1 #2 #3 #4 #5 #6 {%
  \!qshade (#1,#2,#3) (#4,#5,#6)
  \!Shadewhat} 

\setlinear

\def\setdashpattern <#1>{%
  \def\!Flist{}\def\!Blist{}\def\!UDlist{}%
  \!countA=0
  \!ecfor\!item:=#1\do{%
    \!dimenA=\!item\relax
    \expandafter\!rightappend\the\!dimenA\withCS{\\}\to\!UDlist%
    \advance\!countA  1
    \ifodd\!countA
      \expandafter\!rightappend\the\!dimenA\withCS{\!Rule}\to\!Flist%
      \expandafter\!leftappend\the\!dimenA\withCS{\!Rule}\to\!Blist%
    \else 
      \expandafter\!rightappend\the\!dimenA\withCS{\!Skip}\to\!Flist%
      \expandafter\!leftappend\the\!dimenA\withCS{\!Skip}\to\!Blist%
    \fi}%
  \!leaderlength=\!zpt
  \def\!Rule##1{\advance\!leaderlength  ##1}%
  \def\!Skip##1{\advance\!leaderlength  ##1}%
  \!Flist%
  \ifdim\!leaderlength>\!zpt 
  \else
    \def\!Flist{\!Skip{24in}}\def\!Blist{\!Skip{24in}}\ignorespaces
    \def\!UDlist{\\{\!zpt}\\{24in}}\ignorespaces
    \!leaderlength=24in
  \fi
  \!dashingon}

\def\!dashingon{%
  \def\!advancedashing{\!!advancedashing}%
  \def\!drawlinearsegment{\!lineardashed}%
  \def\!puthline{\!putdashedhline}%
  \def\!putvline{\!putdashedvline}%
  \ignorespaces}%
\def\!dashingoff{%
  \def\!advancedashing{\relax}%
  \def\!drawlinearsegment{\!linearsolid}%
  \def\!puthline{\!putsolidhline}%
  \def\!putvline{\!putsolidvline}%
  \ignorespaces}

\def\setdots{%
  \!ifnextchar<{\!setdots}{\!setdots<5pt>}}
\def\!setdots<#1>{%
  \!dimenB=#1\advance\!dimenB -\plotsymbolspacing
  \ifdim\!dimenB<\!zpt
    \!dimenB=\!zpt
  \fi
\setdashpattern <\plotsymbolspacing,\!dimenB>}
 
\def\setdotsnear <#1> for <#2>{%
  \!dimenB=#2\relax  \advance\!dimenB -.05pt  
  \!dimenC=#1\relax  \!countA=\!dimenC 
  \!dimenD=\!dimenB  \advance\!dimenD .5\!dimenC  \!countB=\!dimenD
  \divide \!countB  \!countA
  \ifnum 1>\!countB 
    \!countB=1
  \fi
  \divide\!dimenB  \!countB
  \setdots <\!dimenB>}
 
\def\setdashes{%
  \!ifnextchar<{\!setdashes}{\!setdashes<5pt>}}
\def\!setdashes<#1>{\setdashpattern <#1,#1>}
 
\def\setdashesnear <#1> for <#2>{%
  \!dimenB=#2\relax  
  \!dimenC=#1\relax  \!countA=\!dimenC 
  \!dimenD=\!dimenB  \advance\!dimenD .5\!dimenC  \!countB=\!dimenD
  \divide \!countB  \!countA
  \ifodd \!countB 
  \else 
    \advance \!countB  1
  \fi
  \divide\!dimenB  \!countB
  \setdashes <\!dimenB>}
 
\def\setsolid{%
  \def\!Flist{\!Rule{24in}}\def\!Blist{\!Rule{24in}}%
  \def\!UDlist{\\{24in}\\{\!zpt}}%
  \!dashingoff}  
\setsolid


 
  
 
\def\!divide#1#2#3{%
  \!dimenB=#1
  \!dimenC=#2
  \!dimenD=\!dimenB
  \divide \!dimenD \!dimenC
  \!dimenA=\!dimenD
  \multiply\!dimenD \!dimenC
  \advance\!dimenB -\!dimenD
  \!dimenD=\!dimenC
    \ifdim\!dimenD<\!zpt \!dimenD=-\!dimenD 
  \fi
  \ifdim\!dimenD<64pt
    \!divstep[\!tfs]\!divstep[\!tfs]%
  \else 
    \!!divide
  \fi
  #3=\!dimenA\ignorespaces}

\def\!!divide{%
  \ifdim\!dimenD<256pt
    \!divstep[64]\!divstep[32]\!divstep[32]%
  \else 
    \!divstep[8]\!divstep[8]\!divstep[8]\!divstep[8]\!divstep[8]%
    \!dimenA=2\!dimenA
  \fi}

\def\!divstep[#1]{
  \!dimenB=#1\!dimenB
  \!dimenD=\!dimenB
    \divide \!dimenD by \!dimenC
  \!dimenA=#1\!dimenA
    \advance\!dimenA by \!dimenD%
  \multiply\!dimenD by \!dimenC
    \advance\!dimenB by -\!dimenD}
 
\def\Divide <#1> by <#2> forming <#3> {%
  \!divide{#1}{#2}{#3}}

 
 

 

\def\ellipticalarc axes ratio #1:#2 #3 degrees from #4 #5 center at #6 #7 {%
  \!angle=#3pt\relax
  \ifdim\!angle>\!zpt 
    \def\!sign{}
  \else 
    \def\!sign{-}\!angle=-\!angle
  \fi
  \!xxloc=\!M{#6}\!xunit
  \!yyloc=\!M{#7}\!yunit     
  \!xxS=\!M{#4}\!xunit
  \!yyS=\!M{#5}\!yunit
  \advance\!xxS -\!xxloc
  \advance\!yyS -\!yyloc
  \!divide\!xxS{#1pt}\!xxS 
  \!divide\!yyS{#2pt}\!yyS 
  \let\!MC=\!M
  \!setdimenmode
  \!xS=#1\!xxS  \advance\!xS\!xxloc
  \!yS=#2\!yyS  \advance\!yS\!yyloc
  \!start (\!xS,\!yS)%
  \!loop\ifdim\!angle>14.9999pt
    \!rotate(\!xxS,\!yyS)by(\!cos,\!sign\!sin)to(\!xxM,\!yyM) 
    \!rotate(\!xxM,\!yyM)by(\!cos,\!sign\!sin)to(\!xxE,\!yyE)
    \!xM=#1\!xxM  \advance\!xM\!xxloc  \!yM=#2\!yyM  \advance\!yM\!yyloc
    \!xE=#1\!xxE  \advance\!xE\!xxloc  \!yE=#2\!yyE  \advance\!yE\!yyloc
    \!qjoin (\!xM,\!yM) (\!xE,\!yE)
    \!xxS=\!xxE  \!yyS=\!yyE 
    \advance \!angle -15pt
  \repeat
  \ifdim\!angle>\!zpt
    \!angle=100.53096\!angle
    \divide \!angle 360 
    \!sinandcos\!angle\!!sin\!!cos
    \!rotate(\!xxS,\!yyS)by(\!!cos,\!sign\!!sin)to(\!xxM,\!yyM) 
    \!rotate(\!xxM,\!yyM)by(\!!cos,\!sign\!!sin)to(\!xxE,\!yyE)
    \!xM=#1\!xxM  \advance\!xM\!xxloc  \!yM=#2\!yyM  \advance\!yM\!yyloc
    \!xE=#1\!xxE  \advance\!xE\!xxloc  \!yE=#2\!yyE  \advance\!yE\!yyloc
    \!qjoin (\!xM,\!yM) (\!xE,\!yE)
  \fi
  \let\!M=\!MC
  \ignorespaces}

\def\!rotate(#1,#2)by(#3,#4)to(#5,#6){%
  \!dimenA=#3#1\advance \!dimenA -#4#2
  \!dimenB=#3#2\advance \!dimenB  #4#1
  \divide \!dimenA 32  \divide \!dimenB 32 
  #5=\!dimenA  #6=\!dimenB
  \ignorespaces}
\def\!sin{4.17684}
\def\!cos{31.72624}

\def\!sinandcos#1#2#3{%
 \!dimenD=#1
 \!dimenA=\!dimenD
 \!dimenB=32pt
 \!removept\!dimenD\!value
 \!dimenC=\!dimenD
 \!dimenC=\!value\!dimenC \divide\!dimenC by 64 
 \advance\!dimenB by -\!dimenC
 \!dimenC=\!value\!dimenC \divide\!dimenC by 96 
 \advance\!dimenA by -\!dimenC
 \!dimenC=\!value\!dimenC \divide\!dimenC by 128 
 \advance\!dimenB by \!dimenC%
 \!removept\!dimenA#2
 \!removept\!dimenB#3
 \ignorespaces}




\def\putrule#1from #2 #3 to #4 #5 {%
  \!xloc=\!M{#2}\!xunit  \!xxloc=\!M{#4}\!xunit%
  \!yloc=\!M{#3}\!yunit  \!yyloc=\!M{#5}\!yunit%
  \!dxpos=\!xxloc  \advance\!dxpos by -\!xloc
  \!dypos=\!yyloc  \advance\!dypos by -\!yloc
  \ifdim\!dypos=\!zpt
    \def\!!Line{\!puthline{#1}}\ignorespaces
  \else
    \ifdim\!dxpos=\!zpt
      \def\!!Line{\!putvline{#1}}\ignorespaces
    \else 
       \def\!!Line{}
    \fi
  \fi
  \let\!ML=\!M
  \!setdimenmode
  \!!Line%
  \let\!M=\!ML
  \ignorespaces}

\def\!putsolidhline#1{%
  \ifdim\!dxpos>\!zpt 
    \put{\!hline\!dxpos}#1[l] at {\!xloc} {\!yloc}
  \else 
    \put{\!hline{-\!dxpos}}#1[l] at {\!xxloc} {\!yyloc}
  \fi
  \ignorespaces}
 
\def\!putsolidvline#1{%
  \ifdim\!dypos>\!zpt 
    \put{\!vline\!dypos}#1[b] at {\!xloc} {\!yloc}
  \else 
    \put{\!vline{-\!dypos}}#1[b] at {\!xxloc} {\!yyloc}
  \fi
  \ignorespaces}
 
\def\!hline#1{\hbox to #1{\leaders \hrule height\linethickness\hfill}}
\def\!vline#1{\vbox to #1{\leaders \vrule width\linethickness\vfill}}

\def\!putdashedhline#1{%
  \ifdim\!dxpos>\!zpt 
    \!DLsetup\!Flist\!dxpos
    \put{\hbox to \!totalleaderlength{\!hleaders}\!hpartialpattern\!Rtrunc}
      #1[l] at {\!xloc} {\!yloc} 
  \else 
    \!DLsetup\!Blist{-\!dxpos}
    \put{\!hpartialpattern\!Ltrunc\hbox to \!totalleaderlength{\!hleaders}}
      #1[r] at {\!xloc} {\!yloc} 
  \fi
  \ignorespaces}
 
\def\!putdashedvline#1{%
  \!dypos=-\!dypos
  \ifdim\!dypos>\!zpt 
    \!DLsetup\!Flist\!dypos 
    \put{\vbox{\vbox to \!totalleaderlength{\!vleaders}
      \!vpartialpattern\!Rtrunc}}#1[t] at {\!xloc} {\!yloc} 
  \else 
    \!DLsetup\!Blist{-\!dypos}
    \put{\vbox{\!vpartialpattern\!Ltrunc
      \vbox to \!totalleaderlength{\!vleaders}}}#1[b] at {\!xloc} {\!yloc} 
  \fi
  \ignorespaces}

\def\!DLsetup#1#2{
  \let\!RSlist=#1
  \!countB=#2
  \!countA=\!leaderlength
  \divide\!countB by \!countA
  \!totalleaderlength=\!countB\!leaderlength
  \!Rresiduallength=#2%
  \advance \!Rresiduallength by -\!totalleaderlength
  \!Lresiduallength=\!leaderlength
  \advance \!Lresiduallength by -\!Rresiduallength
  \ignorespaces}
 
\def\!hleaders{%
  \def\!Rule##1{\vrule height\linethickness width##1}%
  \def\!Skip##1{\hskip##1}%
  \leaders\hbox{\!RSlist}\hfill}
 
\def\!hpartialpattern#1{%
  \!dimenA=\!zpt \!dimenB=\!zpt 
  \def\!Rule##1{#1{##1}\vrule height\linethickness width\!dimenD}%
  \def\!Skip##1{#1{##1}\hskip\!dimenD}%
  \!RSlist}
 
\def\!vleaders{%
  \def\!Rule##1{\hrule width\linethickness height##1}%
  \def\!Skip##1{\vskip##1}%
  \leaders\vbox{\!RSlist}\vfill}
 
\def\!vpartialpattern#1{%
  \!dimenA=\!zpt \!dimenB=\!zpt 
  \def\!Rule##1{#1{##1}\hrule width\linethickness height\!dimenD}%
  \def\!Skip##1{#1{##1}\vskip\!dimenD}%
  \!RSlist}
 
\def\!Rtrunc#1{\!trunc{#1}>\!Rresiduallength}
\def\!Ltrunc#1{\!trunc{#1}<\!Lresiduallength}
 
\def\!trunc#1#2#3{%
  \!dimenA=\!dimenB         
  \advance\!dimenB by #1%
  \!dimenD=\!dimenB  \ifdim\!dimenD#2#3\!dimenD=#3\fi
  \!dimenC=\!dimenA  \ifdim\!dimenC#2#3\!dimenC=#3\fi
  \advance \!dimenD by -\!dimenC}

\def\!start (#1,#2){%
  \!plotxorigin=\!xorigin  \advance \!plotxorigin by \!plotsymbolxshift
  \!plotyorigin=\!yorigin  \advance \!plotyorigin by \!plotsymbolyshift
  \!xS=\!M{#1}\!xunit \!yS=\!M{#2}\!yunit
  \!rotateaboutpivot\!xS\!yS
  \!copylist\!UDlist\to\!!UDlist
  \!getnextvalueof\!downlength\from\!!UDlist
  \!distacross=\!zpt
  \!intervalno=0 
  \global\totalarclength=\!zpt
  \ignorespaces}

\def\!ljoin (#1,#2){%
  \advance\!intervalno by 1
  \!xE=\!M{#1}\!xunit \!yE=\!M{#2}\!yunit
  \!rotateaboutpivot\!xE\!yE
  \!xdiff=\!xE \advance \!xdiff by -\!xS
  \!ydiff=\!yE \advance \!ydiff by -\!yS
  \!Pythag\!xdiff\!ydiff\!arclength
  \global\advance \totalarclength by \!arclength%
  \!drawlinearsegment
  \!xS=\!xE \!yS=\!yE
  \ignorespaces}

\def\!linearsolid{%
  \!npoints=\!arclength
  \!countA=\plotsymbolspacing
  \divide\!npoints by \!countA
  \ifnum \!npoints<1 
    \!npoints=1 
  \fi
  \divide\!xdiff by \!npoints
  \divide\!ydiff by \!npoints
  \!xpos=\!xS \!ypos=\!yS
  \loop\ifnum\!npoints>-1
    \!plotifinbounds
    \advance \!xpos by \!xdiff
    \advance \!ypos by \!ydiff
    \advance \!npoints by -1
  \repeat
  \ignorespaces}

\def\!lineardashed{%
  \ifdim\!distacross>\!arclength
    \advance \!distacross by -\!arclength  
  \else
    \loop\ifdim\!distacross<\!arclength
      \!divide\!distacross\!arclength\!dimenA
      \!removept\!dimenA\!t
      \!xpos=\!t\!xdiff \advance \!xpos by \!xS
      \!ypos=\!t\!ydiff \advance \!ypos by \!yS
      \!plotifinbounds
      \advance\!distacross by \plotsymbolspacing
      \!advancedashing
    \repeat  
    \advance \!distacross by -\!arclength
  \fi
  \ignorespaces}

\def\!!advancedashing{%
  \advance\!downlength by -\plotsymbolspacing
  \ifdim \!downlength>\!zpt
  \else
    \advance\!distacross by \!downlength
    \!getnextvalueof\!uplength\from\!!UDlist
    \advance\!distacross by \!uplength
    \!getnextvalueof\!downlength\from\!!UDlist
  \fi}

\def\inboundscheckoff{%
  \def\!plotifinbounds{\!plot(\!xpos,\!ypos)}%
  \def\!initinboundscheck{\relax}\ignorespaces}
 
\inboundscheckoff
 
\def\!!plotifinbounds{%
  \ifdim \!xpos<\!checkleft
  \else
    \ifdim \!xpos>\!checkright
    \else
      \ifdim \!ypos<\!checkbot
      \else
         \ifdim \!ypos>\!checktop
         \else
           \!plot(\!xpos,\!ypos)
         \fi 
      \fi
    \fi
  \fi}

\def\!!initinboundscheck{%
  \!checkleft=\!arealloc     \advance\!checkleft by \!xorigin
  \!checkright=\!arearloc    \advance\!checkright by \!xorigin
  \!checkbot=\!areabloc      \advance\!checkbot by \!yorigin
  \!checktop=\!areatloc      \advance\!checktop by \!yorigin}

%


\def\!logten#1#2{%
  \expandafter\!!logten#1\!nil
  \!removept\!dimenF#2%
  \ignorespaces}

\def\!!logten#1#2\!nil{%
  \if -#1%
    \!dimenF=\!zpt
    \def\!next{\ignorespaces}%
  \else
    \if +#1%
      \def\!next{\!!logten#2\!nil}%
    \else
      \if .#1%
        \def\!next{\!!logten0.#2\!nil}%
      \else
        \def\!next{\!!!logten#1#2..\!nil}%
      \fi
    \fi
  \fi
  \!next}

\def\!!!logten#1#2.#3.#4\!nil{%
  \!dimenF=1pt 
  \if 0#1%
    \!!logshift#3pt 
  \else 
    \!logshift#2/
    \!dimenE=#1.#2#3pt 
  \fi 
  \ifdim \!dimenE<\!rootten
    \multiply \!dimenE 10 
    \advance  \!dimenF -1pt
  \fi
  \!dimenG=\!dimenE
    \advance\!dimenG 10pt
  \advance\!dimenE -10pt 
  \multiply\!dimenE 10 
  \!divide\!dimenE\!dimenG\!dimenE
  \!removept\!dimenE\!t
  \!dimenG=\!t\!dimenE
  \!removept\!dimenG\!tt
  \!dimenH=\!tt\!tenAe
    \divide\!dimenH 100
  \advance\!dimenH \!tenAc
  \!dimenH=\!tt\!dimenH
    \divide\!dimenH 100   
  \advance\!dimenH \!tenAa
  \!dimenH=\!t\!dimenH
    \divide\!dimenH 100 
  \advance\!dimenF \!dimenH}

\def\!logshift#1{%
  \if #1/%
    \def\!next{\ignorespaces}%
  \else
    \advance\!dimenF 1pt 
    \def\!next{\!logshift}%
  \fi 
  \!next}
 
 \def\!!logshift#1{%
   \advance\!dimenF -1pt
   \if 0#1%
     \def\!next{\!!logshift}%
   \else
     \if p#1%
       \!dimenF=1pt
       \def\!next{\!dimenE=1p}%
     \else
       \def\!next{\!dimenE=#1.}%
     \fi
   \fi
   \!next}

\def\beginpicture{%
  \setbox\!picbox=\hbox\bgroup%
  \!xleft=\maxdimen  
  \!xright=-\maxdimen
  \!ybot=\maxdimen
  \!ytop=-\maxdimen}
 
\def\endpicture{%
  \ifdim\!xleft=\maxdimen
    \!xleft=\!zpt \!xright=\!zpt \!ybot=\!zpt \!ytop=\!zpt 
  \fi
  \global\!Xleft=\!xleft \global\!Xright=\!xright
  \global\!Ybot=\!ybot \global\!Ytop=\!ytop
  \egroup%
  \ht\!picbox=\!Ytop  \dp\!picbox=-\!Ybot
  \ifdim\!Ybot>\!zpt
  \else 
    \ifdim\!Ytop<\!zpt
      \!Ybot=\!Ytop
    \else
      \!Ybot=\!zpt
    \fi
  \fi
  \hbox{\kern-\!Xleft\lower\!Ybot\box\!picbox\kern\!Xright}}
 
\def\endpicturesave <#1,#2>{%
  \endpicture \global #1=\!Xleft \global #2=\!Ybot \ignorespaces}

\def\setcoordinatesystem{%
  \!ifnextchar{u}{\!getlengths }
    {\!getlengths units <\!xunit,\!yunit>}}
\def\!getlengths units <#1,#2>{%
  \!xunit=#1\relax
  \!yunit=#2\relax
  \!ifcoordmode 
    \let\!SCnext=\!SCccheckforRP
  \else
    \let\!SCnext=\!SCdcheckforRP
  \fi
  \!SCnext}
\def\!SCccheckforRP{%
  \!ifnextchar{p}{\!cgetreference }
    {\!cgetreference point at {\!xref} {\!yref} }}
\def\!cgetreference point at #1 #2 {%
  \edef\!xref{#1}\edef\!yref{#2}%
  \!xorigin=\!xref\!xunit  \!yorigin=\!yref\!yunit  
  \!initinboundscheck 
  \ignorespaces}
\def\!SCdcheckforRP{%
  \!ifnextchar{p}{\!dgetreference}%
    {\ignorespaces}}
\def\!dgetreference point at #1 #2 {%
  \!xorigin=#1\relax  \!yorigin=#2\relax
  \ignorespaces}

\long\def\put#1#2 at #3 #4 {%
  \!setputobject{#1}{#2}%
  \!xpos=\!M{#3}\!xunit  \!ypos=\!M{#4}\!yunit  
  \!rotateaboutpivot\!xpos\!ypos%
  \advance\!xpos -\!xorigin  \advance\!xpos -\!xshift
  \advance\!ypos -\!yorigin  \advance\!ypos -\!yshift
  \kern\!xpos\raise\!ypos\box\!putobject\kern-\!xpos%
  \!doaccounting\ignorespaces}
 
\long\def\multiput #1#2 at {%
  \!setputobject{#1}{#2}%
  \!ifnextchar"{\!putfromfile}{\!multiput}}
\def\!putfromfile"#1"{%
  \expandafter\!multiput \input #1 /}
\def\!multiput{%
  \futurelet\!nextchar\!!multiput}
\def\!!multiput{%
  \if *\!nextchar
    \def\!nextput{\!alsoby}%
  \else
    \if /\!nextchar
      \def\!nextput{\!finishmultiput}%
    \else
      \def\!nextput{\!alsoat}%
    \fi
  \fi
  \!nextput}
\def\!finishmultiput/{%
  \setbox\!putobject=\hbox{}%
  \ignorespaces}
 
\def\!alsoat#1 #2 {%
  \!xpos=\!M{#1}\!xunit  \!ypos=\!M{#2}\!yunit  
  \!rotateaboutpivot\!xpos\!ypos%
  \advance\!xpos -\!xorigin  \advance\!xpos -\!xshift
  \advance\!ypos -\!yorigin  \advance\!ypos -\!yshift
  \kern\!xpos\raise\!ypos\copy\!putobject\kern-\!xpos%
  \!doaccounting
  \!multiput}
 
\def\!alsoby*#1 #2 #3 {%
  \!dxpos=\!M{#2}\!xunit \!dypos=\!M{#3}\!yunit 
  \!rotateonly\!dxpos\!dypos
  \!ntemp=#1%
  \!!loop\ifnum\!ntemp>0
    \advance\!xpos by \!dxpos  \advance\!ypos by \!dypos
    \kern\!xpos\raise\!ypos\copy\!putobject\kern-\!xpos%
    \advance\!ntemp by -1
  \repeat
  \!doaccounting 
  \!multiput}
 
\def\accountingon{\def\!doaccounting{\!!doaccounting}\ignorespaces}

\accountingon
\def\!!doaccounting{%
  \!xtemp=\!xpos  
  \!ytemp=\!ypos
  \ifdim\!xtemp<\!xleft 
     \!xleft=\!xtemp 
  \fi
  \advance\!xtemp by  \!wd 
  \ifdim\!xright<\!xtemp 
    \!xright=\!xtemp
  \fi
  \advance\!ytemp by -\!dp
  \ifdim\!ytemp<\!ybot  
    \!ybot=\!ytemp
  \fi
  \advance\!ytemp by  \!dp
  \advance\!ytemp by  \!ht 
  \ifdim\!ytemp>\!ytop  
    \!ytop=\!ytemp  
  \fi}
 
\long\def\!setputobject#1#2{%
  \setbox\!putobject=\hbox{#1}%
  \!ht=\ht\!putobject  \!dp=\dp\!putobject  \!wd=\wd\!putobject
  \wd\!putobject=\!zpt
  \!xshift=.5\!wd   \!yshift=.5\!ht   \advance\!yshift by -.5\!dp
  \edef\!putorientation{#2}%
  \expandafter\!SPOreadA\!putorientation[]\!nil%
  \expandafter\!SPOreadB\!putorientation<\!zpt,\!zpt>\!nil\ignorespaces}
 
\def\!SPOreadA#1[#2]#3\!nil{\!etfor\!orientation:=#2\do\!SPOreviseshift}
 
\def\!SPOreadB#1<#2,#3>#4\!nil{\advance\!xshift by -#2\advance\!yshift by -#3}
 
\def\!SPOreviseshift{%
  \if l\!orientation 
    \!xshift=\!zpt
  \else 
    \if r\!orientation 
      \!xshift=\!wd
    \else 
      \if b\!orientation
        \!yshift=-\!dp
      \else 
        \if B\!orientation 
          \!yshift=\!zpt
        \else 
          \if t\!orientation 
            \!yshift=\!ht
          \fi 
        \fi
      \fi
    \fi
  \fi}

\long\def\!dimenput#1#2(#3,#4){%
  \!setputobject{#1}{#2}%
  \!xpos=#3\advance\!xpos by -\!xshift
  \!ypos=#4\advance\!ypos by -\!yshift
  \kern\!xpos\raise\!ypos\box\!putobject\kern-\!xpos%
  \!doaccounting\ignorespaces}

\def\!setdimenmode{%
  \let\!M=\!M!!\ignorespaces}
\def\!setcoordmode{%
  \let\!M=\!M!\ignorespaces}
\def\!ifcoordmode{%
  \ifx \!M \!M!}
\def\!ifdimenmode{%
  \ifx \!M \!M!!}
\def\!M!#1#2{#1#2} 
\def\!M!!#1#2{#1}
\!setcoordmode
\let\setdimensionmode=\!setdimenmode
\let\setcoordinatemode=\!setcoordmode




\def\!stack[#1]{%
  \let\!lglue=\hfill \let\!rglue=\hfill
  \expandafter\let\csname !#1glue\endcsname=\relax
  \!ifnextchar<{\!!stack}{\!!stack<\stackleading>}}
\def\!!stack<#1>#2{%
  \vbox{\def\!valueslist{}\!ecfor\!value:=#2\do{%
    \expandafter\!rightappend\!value\withCS{\\}\to\!valueslist}%
    \!lop\!valueslist\to\!value
    \let\\=\cr\lineskiplimit=\maxdimen\lineskip=#1%
    \baselineskip=-1000pt\halign{\!lglue##\!rglue\cr \!value\!valueslist\cr}}%
  \ignorespaces}


\def\!lines[#1]#2{%
  \let\!lglue=\hfill \let\!rglue=\hfill
  \expandafter\let\csname !#1glue\endcsname=\relax
  \vbox{\halign{\!lglue##\!rglue\cr #2\crcr}}%
  \ignorespaces}


\def\!Lines[#1]#2{%
  \let\!lglue=\hfill \let\!rglue=\hfill
  \expandafter\let\csname !#1glue\endcsname=\relax
  \vtop{\halign{\!lglue##\!rglue\cr #2\crcr}}%
  \ignorespaces}

 
 
 
\def\setplotsymbol(#1#2){%
  \!setputobject{#1}{#2}
  \setbox\!plotsymbol=\box\!putobject%
  \!plotsymbolxshift=\!xshift 
  \!plotsymbolyshift=\!yshift 
  \ignorespaces}
 
\setplotsymbol({\fiverm .})

 
\def\!!plot(#1,#2){%
  \!dimenA=-\!plotxorigin \advance \!dimenA by #1
  \!dimenB=-\!plotyorigin \advance \!dimenB by #2
  \kern\!dimenA\raise\!dimenB\copy\!plotsymbol\kern-\!dimenA%
  \ignorespaces}
 
\def\!!!plot(#1,#2){%
  \!dimenA=-\!plotxorigin \advance \!dimenA by #1
  \!dimenB=-\!plotyorigin \advance \!dimenB by #2
  \kern\!dimenA\raise\!dimenB\copy\!plotsymbol\kern-\!dimenA%
  \!countE=\!dimenA
  \!countF=\!dimenB
  \immediate\write\!replotfile{\the\!countE,\the\!countF.}%
  \ignorespaces}

\def\savelinesandcurves on "#1" {%
  \immediate\closeout\!replotfile
  \immediate\openout\!replotfile=#1%
  \let\!plot=\!!!plot}

\def\dontsavelinesandcurves {%
  \let\!plot=\!!plot}
\dontsavelinesandcurves

{\catcode`\%=11\xdef\!Commentsignal{
\def\writesavefile#1 {%
  \immediate\write\!replotfile{\!Commentsignal #1}%
  \ignorespaces}

\def\replot"#1" {%
  \expandafter\!replot\input #1 /}
\def\!replot#1,#2. {%
  \!dimenA=#1sp
  \kern\!dimenA\raise#2sp\copy\!plotsymbol\kern-\!dimenA
  \futurelet\!nextchar\!!replot}
\def\!!replot{%
  \if /\!nextchar 
    \def\!next{\!finish}%
  \else
    \def\!next{\!replot}%
  \fi
  \!next}


 
 
\def\!Pythag#1#2#3{%
  \!dimenE=#1\relax                                     
  \ifdim\!dimenE<\!zpt 
    \!dimenE=-\!dimenE 
  \fi
  \!dimenF=#2\relax
  \ifdim\!dimenF<\!zpt 
    \!dimenF=-\!dimenF 
  \fi
  \advance \!dimenF by \!dimenE
  \ifdim\!dimenF=\!zpt 
    \!dimenG=\!zpt
  \else 
    \!divide{8\!dimenE}\!dimenF\!dimenE
    \advance\!dimenE by -4pt
      \!dimenE=2\!dimenE
    \!removept\!dimenE\!!t
    \!dimenE=\!!t\!dimenE
    \advance\!dimenE by 64pt
    \divide \!dimenE by 2
    \!dimenH=7pt
    \!!Pythag\!!Pythag\!!Pythag
    \!removept\!dimenH\!!t
    \!dimenG=\!!t\!dimenF
    \divide\!dimenG by 8
  \fi
  #3=\!dimenG
  \ignorespaces}

\def\!!Pythag{
  \!divide\!dimenE\!dimenH\!dimenI
  \advance\!dimenH by \!dimenI
    \divide\!dimenH by 2}

\def\placehypotenuse for <#1> and <#2> in <#3> {%
  \!Pythag{#1}{#2}{#3}}

 
 
 
\def\!qjoin (#1,#2) (#3,#4){%
  \advance\!intervalno by 1
  \!ifcoordmode
    \edef\!xmidpt{#1}\edef\!ymidpt{#2}%
  \else
    \!dimenA=#1\relax \edef\!xmidpt{\the\!dimenA}%
    \!dimenA=#2\relax \edef\!xmidpt{\the\!dimenA}%
  \fi
  \!xM=\!M{#1}\!xunit  \!yM=\!M{#2}\!yunit   \!rotateaboutpivot\!xM\!yM
  \!xE=\!M{#3}\!xunit  \!yE=\!M{#4}\!yunit   \!rotateaboutpivot\!xE\!yE
%
  \!dimenA=\!xM  \advance \!dimenA by -\!xS
  \!dimenB=\!xE  \advance \!dimenB by -\!xM
  \!xB=3\!dimenA \advance \!xB by -\!dimenB
  \!xC=2\!dimenB \advance \!xC by -2\!dimenA
%
  \!dimenA=\!yM  \advance \!dimenA by -\!yS%
  \!dimenB=\!yE  \advance \!dimenB by -\!yM%
  \!yB=3\!dimenA \advance \!yB by -\!dimenB%
  \!yC=2\!dimenB \advance \!yC by -2\!dimenA%
%
  \!xprime=\!xB  \!yprime=\!yB
  \!dxprime=.5\!xC  \!dyprime=.5\!yC
  \!getf \!midarclength=\!dimenA
  \!getf \advance \!midarclength by 4\!dimenA
  \!getf \advance \!midarclength by \!dimenA
  \divide \!midarclength by 12
%
  \!arclength=\!dimenA
  \!getf \advance \!arclength by 4\!dimenA
  \!getf \advance \!arclength by \!dimenA
  \divide \!arclength by 12
  \advance \!arclength by \!midarclength
  \global\advance \totalarclength by \!arclength
%
%
  \ifdim\!distacross>\!arclength 
    \advance \!distacross by -\!arclength
  \else
    \!initinverseinterp
    \loop\ifdim\!distacross<\!arclength
      \!inverseinterp
      \!xpos=\!t\!xC \advance\!xpos by \!xB
        \!xpos=\!t\!xpos \advance \!xpos by \!xS
      \!ypos=\!t\!yC \advance\!ypos by \!yB
        \!ypos=\!t\!ypos \advance \!ypos by \!yS
      \!plotifinbounds
      \advance\!distacross \plotsymbolspacing
      \!advancedashing
    \repeat  
    \advance \!distacross by -\!arclength
  \fi
  \!xS=\!xE
  \!yS=\!yE
  \ignorespaces}

\def\!getf{\!Pythag\!xprime\!yprime\!dimenA%
  \advance\!xprime by \!dxprime
  \advance\!yprime by \!dyprime}

\def\!initinverseinterp{%
  \ifdim\!arclength>\!zpt
    \!divide{8\!midarclength}\!arclength\!dimenE
    \ifdim\!dimenE<\!wmin \!setinverselinear
    \else 
      \ifdim\!dimenE>\!wmax \!setinverselinear
      \else
        \def\!inverseinterp{\!inversequad}\ignorespaces
%
%
         \!removept\!dimenE\!Ew
         \!dimenF=-\!Ew\!dimenE
         \advance\!dimenF by 32pt
         \!dimenG=8pt 
         \advance\!dimenG by -\!dimenE
         \!dimenG=\!Ew\!dimenG
         \!divide\!dimenF\!dimenG\!beta
         \!gamma=1pt
         \advance \!gamma by -\!beta
      \fi
    \fi
  \fi
  \ignorespaces}

\def\!inversequad{%
  \!divide\!distacross\!arclength\!dimenG
  \!removept\!dimenG\!v
  \!dimenG=\!v\!gamma
  \advance\!dimenG by \!beta
  \!dimenG=\!v\!dimenG
  \!removept\!dimenG\!t}

\def\!setinverselinear{%
  \def\!inverseinterp{\!inverselinear}%
  \divide\!dimenE by 8 \!removept\!dimenE\!t
  \!countC=\!intervalno \multiply \!countC 2
  \!countB=\!countC     \advance \!countB -1
  \!countA=\!countB     \advance \!countA -1
  \wlog{\the\!countB th point (\!xmidpt,\!ymidpt) being plotted 
    doesn't lie in the}%
  \wlog{ middle third of the arc between the \the\!countA th 
    and \the\!countC th points:}%
  \wlog{ [arc length \the\!countA\space to \the\!countB]/[arc length 
    \the \!countA\space to \the\!countC]=\!t.}%
  \ignorespaces}
 
\def\!inverselinear{%
  \!divide\!distacross\!arclength\!dimenG
  \!removept\!dimenG\!t}

 

\def\startrotation{%
  \let\!rotateaboutpivot=\!!rotateaboutpivot
  \let\!rotateonly=\!!rotateonly
  \!ifnextchar{b}{\!getsincos }%
    {\!getsincos by {\!cosrotationangle} {\!sinrotationangle} }}
\def\!getsincos by #1 #2 {%
  \edef\!cosrotationangle{#1}%
  \edef\!sinrotationangle{#2}%
  \!ifcoordmode 
    \let\!ROnext=\!ccheckforpivot
  \else
    \let\!ROnext=\!dcheckforpivot
  \fi
  \!ROnext}
\def\!ccheckforpivot{%
  \!ifnextchar{a}{\!cgetpivot}%
    {\!cgetpivot about {\!xpivotcoord} {\!ypivotcoord} }}
\def\!cgetpivot about #1 #2 {%
  \edef\!xpivotcoord{#1}%
  \edef\!ypivotcoord{#2}%
  \!xpivot=#1\!xunit  \!ypivot=#2\!yunit
  \ignorespaces}
\def\!dcheckforpivot{%
  \!ifnextchar{a}{\!dgetpivot}{\ignorespaces}}
\def\!dgetpivot about #1 #2 {%
  \!xpivot=#1\relax  \!ypivot=#2\relax
  \ignorespaces}

\def\stoprotation{%
  \let\!rotateaboutpivot=\!!!rotateaboutpivot
  \let\!rotateonly=\!!!rotateonly
  \ignorespaces}
 
\def\!!rotateaboutpivot#1#2{%
  \!dimenA=#1\relax  \advance\!dimenA -\!xpivot
  \!dimenB=#2\relax  \advance\!dimenB -\!ypivot
  \!dimenC=\!cosrotationangle\!dimenA
    \advance \!dimenC -\!sinrotationangle\!dimenB
  \!dimenD=\!cosrotationangle\!dimenB
    \advance \!dimenD  \!sinrotationangle\!dimenA
  \advance\!dimenC \!xpivot  \advance\!dimenD \!ypivot
  #1=\!dimenC  #2=\!dimenD
  \ignorespaces}

\def\!!rotateonly#1#2{%
  \!dimenA=#1\relax  \!dimenB=#2\relax 
  \!dimenC=\!cosrotationangle\!dimenA
    \advance \!dimenC -\!rotsign\!sinrotationangle\!dimenB
  \!dimenD=\!cosrotationangle\!dimenB
    \advance \!dimenD  \!rotsign\!sinrotationangle\!dimenA
  #1=\!dimenC  #2=\!dimenD
  \ignorespaces}
\def\!rotsign{}
\def\!!!rotateaboutpivot#1#2{\relax}
\def\!!!rotateonly#1#2{\relax}
\stoprotation

\def\!reverserotateonly#1#2{%
  \def\!rotsign{-}%
  \!rotateonly{#1}{#2}%
  \def\!rotsign{}%
  \ignorespaces}

\def\!getspan span <#1>{%
  \!dshade=#1\relax
  \!ifcoordmode 
    \let\!GRnext=\!GRccheckforAP
  \else
    \let\!GRnext=\!GRdcheckforAP
  \fi
  \!GRnext}
\def\!GRccheckforAP{%
  \!ifnextchar{p}{\!cgetanchor }
    {\!cgetanchor point at {\!xshadesave} {\!yshadesave} }}
\def\!cgetanchor point at #1 #2 {%
  \edef\!xshadesave{#1}\edef\!yshadesave{#2}%
  \!xshade=\!xshadesave\!xunit  \!yshade=\!yshadesave\!yunit
  \ignorespaces}
\def\!GRdcheckforAP{%
  \!ifnextchar{p}{\!dgetanchor}%
    {\ignorespaces}}
\def\!dgetanchor point at #1 #2 {%
  \!xshade=#1\relax  \!yshade=#2\relax
  \ignorespaces}

\def\setshadesymbol{%
  \!ifnextchar<{\!setshadesymbol}{\!setshadesymbol<,,,> }}

\def\!setshadesymbol <#1,#2,#3,#4> (#5#6){%
  \!setputobject{#5}{#6}%
  \setbox\!shadesymbol=\box\!putobject%
  \!shadesymbolxshift=\!xshift \!shadesymbolyshift=\!yshift
%
  \!dimenA=\!xshift \advance\!dimenA \!smidge
  \!override\!dimenA{#1}\!lshrinkage%
  \!dimenA=\!wd \advance \!dimenA -\!xshift
    \advance\!dimenA \!smidge
    \!override\!dimenA{#2}\!rshrinkage
  \!dimenA=\!dp \advance \!dimenA \!yshift
    \advance\!dimenA \!smidge
    \!override\!dimenA{#3}\!bshrinkage
  \!dimenA=\!ht \advance \!dimenA -\!yshift
    \advance\!dimenA \!smidge
    \!override\!dimenA{#4}\!tshrinkage
  \ignorespaces}
\def\!smidge{-.2pt}%

\def\!override#1#2#3{%
  \edef\!!override{#2}%
  \ifx \!!override\empty
    #3=#1\relax
  \else
    \if z\!!override
      #3=\!zpt
    \else
      \ifx \!!override\!blankz
        #3=\!zpt
      \else
        #3=#2\relax
      \fi
    \fi
  \fi
  \ignorespaces}
\def\!blankz{ z}

\setshadesymbol ({\fiverm .})

\def\!startvshade#1(#2,#3,#4){%
  \let\!!xunit=\!xunit%
  \let\!!yunit=\!yunit%
  \let\!!xshade=\!xshade%
  \let\!!yshade=\!yshade%
  \def\!getshrinkages{\!vgetshrinkages}%
  \let\!setshadelocation=\!vsetshadelocation%
  \!xS=\!M{#2}\!!xunit
  \!ybS=\!M{#3}\!!yunit
  \!ytS=\!M{#4}\!!yunit
  \!shadexorigin=\!xorigin  \advance \!shadexorigin \!shadesymbolxshift
  \!shadeyorigin=\!yorigin  \advance \!shadeyorigin \!shadesymbolyshift
  \ignorespaces}
 
\def\!starthshade#1(#2,#3,#4){%
  \let\!!xunit=\!yunit%
  \let\!!yunit=\!xunit%
  \let\!!xshade=\!yshade%
  \let\!!yshade=\!xshade%
  \def\!getshrinkages{\!hgetshrinkages}%
  \let\!setshadelocation=\!hsetshadelocation%
  \!xS=\!M{#2}\!!xunit
  \!ybS=\!M{#3}\!!yunit
  \!ytS=\!M{#4}\!!yunit
  \!shadexorigin=\!xorigin  \advance \!shadexorigin \!shadesymbolxshift
  \!shadeyorigin=\!yorigin  \advance \!shadeyorigin \!shadesymbolyshift
  \ignorespaces}

\def\!lattice#1#2#3#4#5{%
  \!dimenA=#1
  \!dimenB=#2
  \!countB=\!dimenB
%
  \!dimenC=#3
  \advance\!dimenC -\!dimenA
  \!countA=\!dimenC
  \divide\!countA \!countB
  \ifdim\!dimenC>\!zpt
    \!dimenD=\!countA\!dimenB
    \ifdim\!dimenD<\!dimenC
      \advance\!countA 1 
    \fi
  \fi
  \!dimenC=\!countA\!dimenB
    \advance\!dimenC \!dimenA
  #4=\!countA
  #5=\!dimenC
  \ignorespaces}

\def\!qshade#1(#2,#3,#4)#5(#6,#7,#8){%
  \!xM=\!M{#2}\!!xunit
  \!ybM=\!M{#3}\!!yunit
  \!ytM=\!M{#4}\!!yunit
  \!xE=\!M{#6}\!!xunit
  \!ybE=\!M{#7}\!!yunit
  \!ytE=\!M{#8}\!!yunit
  \!getcoeffs\!xS\!ybS\!xM\!ybM\!xE\!ybE\!ybB\!ybC
  \!getcoeffs\!xS\!ytS\!xM\!ytM\!xE\!ytE\!ytB\!ytC
  \def\!getylimits{\!qgetylimits}%
  \!shade{#1}\ignorespaces}
 
\def\!lshade#1(#2,#3,#4){%
  \!xE=\!M{#2}\!!xunit
  \!ybE=\!M{#3}\!!yunit
  \!ytE=\!M{#4}\!!yunit
  \!dimenE=\!xE  \advance \!dimenE -\!xS
  \!dimenC=\!ytE \advance \!dimenC -\!ytS
  \!divide\!dimenC\!dimenE\!ytB
  \!dimenC=\!ybE \advance \!dimenC -\!ybS
  \!divide\!dimenC\!dimenE\!ybB
  \def\!getylimits{\!lgetylimits}%
  \!shade{#1}\ignorespaces}
 
\def\!getcoeffs#1#2#3#4#5#6#7#8{%
  \!dimenC=#4\advance \!dimenC -#2
  \!dimenE=#3\advance \!dimenE -#1
  \!divide\!dimenC\!dimenE\!dimenF
  \!dimenC=#6\advance \!dimenC -#4
  \!dimenH=#5\advance \!dimenH -#3
  \!divide\!dimenC\!dimenH\!dimenG
  \advance\!dimenG -\!dimenF
  \advance \!dimenH \!dimenE
  \!divide\!dimenG\!dimenH#8
  \!removept#8\!t
  #7=-\!t\!dimenE
  \advance #7\!dimenF
  \ignorespaces}

\def\!shade#1{%
  \!getshrinkages#1<,,,>\!nil
  \advance \!dimenE \!xS
  \!lattice\!!xshade\!dshade\!dimenE
    \!parity\!xpos
  \!dimenF=-\!dimenF
    \advance\!dimenF \!xE
  \!loop\!not{\ifdim\!xpos>\!dimenF}
    \!shadecolumn%
    \advance\!xpos \!dshade
    \advance\!parity 1
  \repeat
  \!xS=\!xE
  \!ybS=\!ybE
  \!ytS=\!ytE
  \ignorespaces}

\def\!vgetshrinkages#1<#2,#3,#4,#5>#6\!nil{%
  \!override\!lshrinkage{#2}\!dimenE
  \!override\!rshrinkage{#3}\!dimenF
  \!override\!bshrinkage{#4}\!dimenG
  \!override\!tshrinkage{#5}\!dimenH
  \ignorespaces}
\def\!hgetshrinkages#1<#2,#3,#4,#5>#6\!nil{%
  \!override\!lshrinkage{#2}\!dimenG
  \!override\!rshrinkage{#3}\!dimenH
  \!override\!bshrinkage{#4}\!dimenE
  \!override\!tshrinkage{#5}\!dimenF
  \ignorespaces}

\def\!shadecolumn{%
  \!dxpos=\!xpos
  \advance\!dxpos -\!xS
  \!removept\!dxpos\!dx
  \!getylimits
  \advance\!ytpos -\!dimenH
  \advance\!ybpos \!dimenG
  \!yloc=\!!yshade
  \ifodd\!parity 
     \advance\!yloc \!dshade
  \fi
  \!lattice\!yloc{2\!dshade}\!ybpos%
    \!countA\!ypos
  \!dimenA=-\!shadexorigin \advance \!dimenA \!xpos
  \loop\!not{\ifdim\!ypos>\!ytpos}
    \!setshadelocation
    \!rotateaboutpivot\!xloc\!yloc%
    \!dimenA=-\!shadexorigin \advance \!dimenA \!xloc
    \!dimenB=-\!shadeyorigin \advance \!dimenB \!yloc
    \kern\!dimenA \raise\!dimenB\copy\!shadesymbol \kern-\!dimenA
    \advance\!ypos 2\!dshade
  \repeat
  \ignorespaces}
 
\def\!qgetylimits{%
  \!dimenA=\!dx\!ytC              
  \advance\!dimenA \!ytB
  \!ytpos=\!dx\!dimenA
  \advance\!ytpos \!ytS
  \!dimenA=\!dx\!ybC              
  \advance\!dimenA \!ybB
  \!ybpos=\!dx\!dimenA
  \advance\!ybpos \!ybS}
 
\def\!lgetylimits{%
  \!ytpos=\!dx\!ytB
  \advance\!ytpos \!ytS
  \!ybpos=\!dx\!ybB
  \advance\!ybpos \!ybS}
 
\def\!vsetshadelocation{
  \!xloc=\!xpos
  \!yloc=\!ypos}
\def\!hsetshadelocation{
  \!xloc=\!ypos
  \!yloc=\!xpos}





\def\!axisticks {%
  \def\!nextkeyword##1 {%
    \expandafter\ifx\csname !ticks##1\endcsname \relax
      \def\!next{\!fixkeyword{##1}}%
    \else
      \def\!next{\csname !ticks##1\endcsname}%
    \fi
    \!next}%
  \!axissetup
    \def\!axissetup{\relax}%
  \edef\!ticksinoutsign{\!ticksinoutSign}%
  \!ticklength=\longticklength
  \!tickwidth=\linethickness
  \!gridlinestatus
  \!setticktransform
  \!maketick
  \!tickcase=0
  \def\!LTlist{}%
  \!nextkeyword}

\def\ticksout{%
  \def\!ticksinoutSign{+}}

\ticksout

\def\nogridlines{%
  \def\!gridlinestatus{\!gridlinestoofalse}}
\nogridlines

\def\loggedticks{%
  \def\!setticktransform{\let\!ticktransform=\!logten}}
\def\unloggedticks{%
  \def\!setticktransform{\let\!ticktransform=\!donothing}}
\def\!donothing#1#2{\def#2{#1}}
\unloggedticks

\expandafter\def\csname !ticks/\endcsname{%
  \!not {\ifx \!LTlist\empty}
    \!placetickvalues
  \fi
  \def\!tickvalueslist{}%
  \def\!LTlist{}%
  \expandafter\csname !axis/\endcsname}

\def\!maketick{%
  \setbox\!boxA=\hbox{%
    \beginpicture
      \!setdimenmode
      \setcoordinatesystem point at {\!zpt} {\!zpt}   
      \linethickness=\!tickwidth
      \ifdim\!ticklength>\!zpt
        \putrule from {\!zpt} {\!zpt} to
          {\!ticksinoutsign\!tickxsign\!ticklength}
          {\!ticksinoutsign\!tickysign\!ticklength}
      \fi
      \if!gridlinestoo
        \putrule from {\!zpt} {\!zpt} to
          {-\!tickxsign\!xaxislength} {-\!tickysign\!yaxislength}
      \fi
    \endpicturesave <\!Xsave,\!Ysave>}%
    \wd\!boxA=\!zpt}
  
\def\!ticksin{%
  \def\!ticksinoutsign{-}%
  \!maketick
  \!nextkeyword}

\def\!ticksout{%
  \def\!ticksinoutsign{+}%
  \!maketick
  \!nextkeyword}

\def\!tickslength<#1> {%
  \!ticklength=#1\relax
  \!maketick
  \!nextkeyword}

\def\!tickslong{%
  \!tickslength<\longticklength> }

\def\!ticksshort{%
  \!tickslength<\shortticklength> }

\def\!tickswidth<#1> {%
  \!tickwidth=#1\relax
  \!maketick
  \!nextkeyword}

\def\!ticksandacross{%
  \!gridlinestootrue
  \!maketick
  \!nextkeyword}

\def\!ticksbutnotacross{%
  \!gridlinestoofalse
  \!maketick
  \!nextkeyword}

\def\!tickslogged{%
  \let\!ticktransform=\!logten
  \!nextkeyword}

\def\!ticksunlogged{%
  \let\!ticktransform=\!donothing
  \!nextkeyword}

\def\!ticksunlabeled{%
  \!tickcase=0
  \!nextkeyword}

\def\!ticksnumbered{%
  \!tickcase=1
  \!nextkeyword}

\def\!tickswithvalues#1/ {%
  \edef\!tickvalueslist{#1! /}%
  \!tickcase=2
  \!nextkeyword}

\def\!ticksquantity#1 {%
  \ifnum #1>1
    \!updatetickoffset
    \!countA=#1\relax
    \advance \!countA -1
    \!ticklocationincr=\!axisLength
      \divide \!ticklocationincr \!countA
    \!ticklocation=\!axisstart
    \loop \!not{\ifdim \!ticklocation>\!axisend}
      \!placetick\!ticklocation
      \ifcase\!tickcase
          \relax 
        \or
          \relax 
        \or
          \expandafter\!gettickvaluefrom\!tickvalueslist
          \edef\!tickfield{{\the\!ticklocation}{\!value}}%
          \expandafter\!listaddon\expandafter{\!tickfield}\!LTlist%
      \fi
      \advance \!ticklocation \!ticklocationincr
    \repeat
  \fi
  \!nextkeyword}

\def\!ticksat#1 {%
  \!updatetickoffset
  \edef\!Loc{#1}%
  \if /\!Loc
    \def\next{\!nextkeyword}%
  \else
    \!ticksincommon
    \def\next{\!ticksat}%
  \fi
  \next}    
      
\def\!ticksfrom#1 to #2 by #3 {%
  \!updatetickoffset
  \edef\!arg{#3}%
  \expandafter\!separate\!arg\!nil
  \!scalefactor=1
  \expandafter\!countfigures\!arg/
  \edef\!arg{#1}%
  \!scaleup\!arg by\!scalefactor to\!countE
  \edef\!arg{#2}%
  \!scaleup\!arg by\!scalefactor to\!countF
  \edef\!arg{#3}%
  \!scaleup\!arg by\!scalefactor to\!countG
  \loop \!not{\ifnum\!countE>\!countF}
    \ifnum\!scalefactor=1
      \edef\!Loc{\the\!countE}%
    \else
      \!scaledown\!countE by\!scalefactor to\!Loc
    \fi
    \!ticksincommon
    \advance \!countE \!countG
  \repeat
  \!nextkeyword}

\def\!updatetickoffset{%
  \!dimenA=\!ticksinoutsign\!ticklength
  \ifdim \!dimenA>\!offset
    \!offset=\!dimenA
  \fi}

\def\!placetick#1{%
  \if!xswitch
    \!xpos=#1\relax
    \!ypos=\!axisylevel
  \else
    \!xpos=\!axisxlevel
    \!ypos=#1\relax
  \fi
  \advance\!xpos \!Xsave
  \advance\!ypos \!Ysave
  \kern\!xpos\raise\!ypos\copy\!boxA\kern-\!xpos
  \ignorespaces}

\def\!gettickvaluefrom#1 #2 /{%
  \edef\!value{#1}%
  \edef\!tickvalueslist{#2 /}%
  \ifx \!tickvalueslist\!endtickvaluelist
    \!tickcase=0
  \fi}
\def\!endtickvaluelist{! /}

\def\!ticksincommon{%
  \!ticktransform\!Loc\!t
  \!ticklocation=\!t\!!unit
  \advance\!ticklocation -\!!origin
  \!placetick\!ticklocation
  \ifcase\!tickcase
    \relax 
  \or 
    \ifdim\!ticklocation<-\!!origin
      \edef\!Loc{$\!Loc$}%
    \fi
    \edef\!tickfield{{\the\!ticklocation}{\!Loc}}%
    \expandafter\!listaddon\expandafter{\!tickfield}\!LTlist%
  \or 
    \expandafter\!gettickvaluefrom\!tickvalueslist
    \edef\!tickfield{{\the\!ticklocation}{\!value}}%
    \expandafter\!listaddon\expandafter{\!tickfield}\!LTlist%
  \fi}

\def\!separate#1\!nil{%
  \!ifnextchar{-}{\!!separate}{\!!!separate}#1\!nil}
\def\!!separate-#1\!nil{%
  \def\!sign{-}%
  \!!!!separate#1..\!nil}
\def\!!!separate#1\!nil{%
  \def\!sign{+}%
  \!!!!separate#1..\!nil}
\def\!!!!separate#1.#2.#3\!nil{%
  \def\!arg{#1}%
  \ifx\!arg\!empty
    \!countA=0
  \else
    \!countA=\!arg
  \fi
  \def\!arg{#2}%
  \ifx\!arg\!empty
    \!countB=0
  \else
    \!countB=\!arg
  \fi}
 
\def\!countfigures#1{%
  \if #1/%
    \def\!next{\ignorespaces}%
  \else
    \multiply\!scalefactor 10
    \def\!next{\!countfigures}%
  \fi
  \!next}

\def\!scaleup#1by#2to#3{%
  \expandafter\!separate#1\!nil
  \multiply\!countA #2\relax
  \advance\!countA \!countB
  \if -\!sign
    \!countA=-\!countA
  \fi
  #3=\!countA
  \ignorespaces}

\def\!scaledown#1by#2to#3{%
  \!countA=#1\relax
  \ifnum \!countA<0 
    \def\!sign{-}
    \!countA=-\!countA
  \else
    \def\!sign{}%
  \fi
  \!countB=\!countA
  \divide\!countB #2\relax
  \!countC=\!countB
    \multiply\!countC #2\relax
  \advance \!countA -\!countC
  \edef#3{\!sign\the\!countB.}
  \!countC=\!countA 
  \ifnum\!countC=0 
    \!countC=1
  \fi
  \multiply\!countC 10
  \!loop \ifnum #2>\!countC
    \edef#3{#3\!zero}%
    \multiply\!countC 10
  \repeat
  \edef#3{#3\the\!countA}
  \ignorespaces}

\def\!placetickvalues{%
  \advance\!offset \tickstovaluesleading
  \if!xswitch
    \setbox\!boxA=\hbox{%
      \def\\##1##2{%
        \!dimenput {##2} [B] (##1,\!axisylevel)}%
      \beginpicture 
        \!LTlist
      \endpicturesave <\!Xsave,\!Ysave>}%
    \!dimenA=\!axisylevel
      \advance\!dimenA -\!Ysave
      \advance\!dimenA \!tickysign\!offset
      \if -\!tickysign
        \advance\!dimenA -\ht\!boxA
      \else
        \advance\!dimenA  \dp\!boxA
      \fi
    \advance\!offset \ht\!boxA 
      \advance\!offset \dp\!boxA
    \!dimenput {\box\!boxA} [Bl] <\!Xsave,\!Ysave> (\!zpt,\!dimenA)
  \else
    \setbox\!boxA=\hbox{%
      \def\\##1##2{%
        \!dimenput {##2} [r] (\!axisxlevel,##1)}%
      \beginpicture 
        \!LTlist
      \endpicturesave <\!Xsave,\!Ysave>}%
    \!dimenA=\!axisxlevel
      \advance\!dimenA -\!Xsave
      \advance\!dimenA \!tickxsign\!offset
      \if -\!tickxsign
        \advance\!dimenA -\wd\!boxA
      \fi
    \advance\!offset \wd\!boxA
    \!dimenput {\box\!boxA} [Bl] <\!Xsave,\!Ysave> (\!dimenA,\!zpt)
  \fi}

\normalgraphs
\catcode`!=12 
\newread\epsffilein    
\newif\ifepsffileok    
\newif\ifepsfbbfound   
\newif\ifepsfverbose   
\newif\ifepsfdraft     
\newdimen\epsfxsize    
\newdimen\epsfysize    
\newdimen\epsftsize    
\newdimen\epsfrsize    
\newdimen\epsftmp      
\newdimen\pspoints     
\pspoints=1bp          
\epsfxsize=0pt         
\epsfysize=0pt         
\def\epsfbox#1{\global\def\epsfllx{72}\global\def\epsflly{72}%
   \global\def\epsfurx{540}\global\def\epsfury{720}%
   \def\lbracket{[}\def\testit{#1}\ifx\testit\lbracket
   \let\next=\epsfgetlitbb\else\let\next=\epsfnormal\fi\next{#1}}%
\def\epsfgetlitbb#1#2 #3 #4 #5]#6{\epsfgrab #2 #3 #4 #5 .\\%
   \epsfsetgraph{#6}}%
\def\epsfnormal#1{\epsfgetbb{#1}\epsfsetgraph{#1}}%
\def\epsfgetbb#1{%
%
%
\openin\epsffilein=#1
\ifeof\epsffilein\errmessage{I couldn't open #1, will ignore it}\else
%
%
   {\epsffileoktrue \chardef\other=12
    \def\do##1{\catcode`##1=\other}\dospecials \catcode`\ =10
    \loop
       \read\epsffilein to \epsffileline
       \ifeof\epsffilein\epsffileokfalse\else
%
%
          \expandafter\epsfaux\epsffileline:. \\%
       \fi
   \ifepsffileok\repeat
   \ifepsfbbfound\else
    \ifepsfverbose\message{No bounding box comment in #1; using defaults}\fi\fi
   }\closein\epsffilein\fi}%
%
%
%
\def\epsfclipoff{\def\epsfclipstring{\ifepsfdraft\space clip\fi}}%
\epsfclipoff
\def\epsfsetgraph#1{%
   \epsfrsize=\epsfury\pspoints
   \advance\epsfrsize by-\epsflly\pspoints
   \epsftsize=\epsfurx\pspoints
   \advance\epsftsize by-\epsfllx\pspoints
%
%
   \epsfxsize\epsfsize\epsftsize\epsfrsize
   \ifnum\epsfxsize=0 \ifnum\epsfysize=0
      \epsfxsize=\epsftsize \epsfysize=\epsfrsize
      \epsfrsize=0pt
%
%
     \else\epsftmp=\epsftsize \divide\epsftmp\epsfrsize
       \epsfxsize=\epsfysize \multiply\epsfxsize\epsftmp
       \multiply\epsftmp\epsfrsize \advance\epsftsize-\epsftmp
       \epsftmp=\epsfysize
       \loop \advance\epsftsize\epsftsize \divide\epsftmp 2
       \ifnum\epsftmp>0
          \ifnum\epsftsize<\epsfrsize\else
             \advance\epsftsize-\epsfrsize \advance\epsfxsize\epsftmp \fi
       \repeat
       \epsfrsize=0pt
     \fi
   \else \ifnum\epsfysize=0
     \epsftmp=\epsfrsize \divide\epsftmp\epsftsize
     \epsfysize=\epsfxsize \multiply\epsfysize\epsftmp   
     \multiply\epsftmp\epsftsize \advance\epsfrsize-\epsftmp
     \epsftmp=\epsfxsize
     \loop \advance\epsfrsize\epsfrsize \divide\epsftmp 2
     \ifnum\epsftmp>0
        \ifnum\epsfrsize<\epsftsize\else
           \advance\epsfrsize-\epsftsize \advance\epsfysize\epsftmp \fi
     \repeat
     \epsfrsize=0pt
    \else
     \epsfrsize=\epsfysize
    \fi
   \fi
%
%
   \ifepsfverbose\message{#1: width=\the\epsfxsize, height=\the\epsfysize}\fi
   \epsftmp=10\epsfxsize \divide\epsftmp\pspoints
   \vbox to\epsfysize{\vfil\hbox to\epsfxsize{%
      \ifnum\epsfrsize=0\relax
        \includegraphics{\ifepsfdraft}%
      \else
        \epsfrsize=10\epsfysize \divide\epsfrsize\pspoints
        \includegraphics{\ifepsfdraft}%
      \fi
      \hfil}}%
\global\epsfxsize=0pt\global\epsfysize=0pt}%
%
%
{\catcode`\%=12 \global\let\epsfpercent=
%
%
\long\def\epsfaux#1#2:#3\\{\ifx#1\epsfpercent
   \def\testit{#2}\ifx\testit\epsfbblit
      \epsfgrab #3 . . . \\%
      \epsffileokfalse
      \global\epsfbbfoundtrue
   \fi\else\ifx#1\par\else\epsffileokfalse\fi\fi}%
%
%
\def\epsfempty{}%
\def\epsfgrab #1 #2 #3 #4 #5\\{%
\global\def\epsfllx{#1}\ifx\epsfllx\epsfempty
      \epsfgrab #2 #3 #4 #5 .\\\else
   \global\def\epsflly{#2}%
   \global\def\epsfurx{#3}\global\def\epsfury{#4}\fi}%
%
%
\def\epsfsize#1#2{\epsfxsize}
%
%

%
%
%
%
%
%
%
\magnification=\magstephalf      
%
%
\vsize=7.5truein                 
\hsize=5.2truein                 
\newskip\stdskip                 
\stdskip=6pt plus3pt minus3pt    
\medskipamount=\stdskip          
\parindent=0pt                   
\parskip=\stdskip                
\abovedisplayskip=\stdskip       
\belowdisplayskip=\stdskip       
\mathsurround=0.75pt             
\overfullrule=0pt                
%
%
\def\ppar{\par\goodbreak\vskip 8pt plus 4pt minus 4pt}     
%
%
\def\stdspace{\hskip 0.75em plus 0.15em\ignorespaces}
%
%
%
%
%
%
%
\def\hexnumber#1{\ifcase#1 0\or 1\or 2\or 3\or 4\or 5\or 6\or 7\or 8\or
 9\or A\or B\or C\or D\or E\or F\fi}
%
%
\font\thirtnmsa=msam10 scaled 1315    
\font\tenmsa=msam10          \font\ninemsa=msam9
\font\sevenmsa=msam7         \font\sixmsa=msam6
\font\fivemsa=msam5
%
%
\newfam\msafam                  \textfont\msafam=\tenmsa
\scriptfont\msafam=\sevenmsa    \scriptscriptfont\msafam=\fivemsa
\edef\hexa{\hexnumber\msafam}        
\def\msa{\fam\msafam\tenmsa}         
%
%
\font\thirtnmsb=msbm10 scaled 1315   
\font\tenmsb=msbm10      \font\ninemsb=msbm9
\font\sevenmsb=msbm7     \font\sixmsb=msbm6
\font\fivemsb=msbm5
%
\newfam\msbfam                   \textfont\msbfam=\tenmsb       
\scriptfont\msbfam=\sevenmsb     \scriptscriptfont\msbfam=\fivemsb
\edef\hexb{\hexnumber\msbfam}    
\def\msb{\fam\msbfam\tenmsb}     
%
%
\font\thirtneufm=eufm10 scaled 1315   
\font\teneufm=eufm10                 \font\nineeufm=eufm9
\font\seveneufm=eufm7                \font\sixeufm=eufm6
\font\fiveeufm=eufm5
%
\newfam\eufmfam                    \textfont\eufmfam=\teneufm
\scriptfont\eufmfam=\seveneufm     \scriptscriptfont\eufmfam=\fiveeufm
\edef\hexf{\hexnumber\eufmfam}      
\def\frak{\fam\eufmfam\teneufm}     
%
%
%
\font\thirtnrm=cmr10 scaled 1315    
\font\ninerm=cmr9                   \font\sixrm=cmr6   
%
\font\thirtni=cmmi10 scaled 1315    
\font\ninei=cmmi9                   \font\sixi=cmmi6  
%
\font\thirtnsy=cmsy10 scaled 1315   
\font\ninesy=cmsy9                  \font\sixsy=cmsy6  
%
\font\thirtnbf=cmbx10 scaled 1315   
\font\ninebf=cmbx9                  \font\sixbf=cmbx6  
%
%
\font\thirtnex=cmex10 scaled 1315   
\font\nineex=cmex9                  
%
%
\font\thirtnit=cmti10 scaled 1315  
\font\nineit=cmti9                  
%
\font\thirtnsl=cmsl10 scaled 1315  
\font\ninesl=cmsl9                  
%
\font\thirtntt=cmtt10 scaled 1315  
\font\ninett=cmtt9                  
%
%
%
%
\def\small{%
%
%
\textfont0=\ninerm \scriptfont0=\sixrm \scriptscriptfont0=\fiverm
\def\rm{\fam0\ninerm}
%
%
\textfont1=\ninei \scriptfont1=\sixi \scriptscriptfont1=\fivei
%
%
\textfont2=\ninesy \scriptfont2=\sixsy \scriptscriptfont2=\fivesy
%
%
\textfont3=\nineex \scriptfont3=\nineex \scriptscriptfont3=\nineex
%
%
\textfont\bffam=\ninebf \scriptfont\bffam=\sixbf
\scriptscriptfont\bffam=\fivebf \def\bf{\fam\bffam\ninebf}%
%
%
\textfont\itfam=\nineit \def\it{\fam\itfam\nineit}%
\textfont\slfam=\ninesl \def\sl{\fam\slfam\ninesl}%
\textfont\ttfam=\ninett \def\tt{\fam\ttfam\ninett}%
%
%
%
\textfont\msafam=\ninemsa \scriptfont\msafam=\sixmsa
\scriptscriptfont\msafam=\fivemsa \def\msa{\fam\msafam\ninemsa}%
%
%
\textfont\msbfam=\ninemsb \scriptfont\msbfam=\sixmsb
\scriptscriptfont\msbfam=\fivemsb \def\msb{\fam\msbfam\ninemsb}%
%
%
\textfont\eufmfam=\nineeufm  \scriptfont\eufmfam=\sixeufm
\scriptscriptfont\eufmfam=\fiveeufm \def\frak{\fam\eufmfam\nineeufm}%
%
%
%
\normalbaselineskip=11pt%
\setbox\strutbox=\hbox{\vrule height8pt depth3pt width0pt}%
%
%
\normalbaselines\rm
%
%
\stdskip=4pt plus2pt minus2pt    
\medskipamount=\stdskip          
\parskip=\stdskip                
\abovedisplayskip=\stdskip       
\belowdisplayskip=\stdskip       
\def\ppar{\par\goodbreak\vskip 6pt plus 3pt minus 3pt}%
%
%
\def\section##1{\global\advance\sectionnumber by 1
\vskip-\lastskip\penalty-800\vskip 20pt plus10pt minus5pt 
\egroup{\bf\number\sectionnumber\quad##1}\bgroup\small         
\vskip 6pt plus3pt minus3pt
\nobreak\resultnumber=1}
}    
%
\def\beginsmall{\bgroup\small}
\let\endsmall\egroup
%
%
%
%
\def\large{%
\textfont0=\thirtnrm \scriptfont0=\ninerm \scriptscriptfont0=\sevenrm
\def\rm{\fam0\thirtnrm}%
\textfont1=\thirtni \scriptfont1=\ninei \scriptscriptfont1=\seveni
\textfont2=\thirtnsy \scriptfont2=\ninesy \scriptscriptfont2=\sevensy
\textfont3=\thirtnex \scriptfont3=\thirtnex \scriptscriptfont3=\thirtnex
\textfont\bffam=\thirtnbf \scriptfont\bffam=\ninebf
\scriptscriptfont\bffam=\sevenbf \def\bf{\fam\bffam\thirtnbf}%
\textfont\itfam=\thirtnit \def\it{\fam\itfam\thirtnit}%
\textfont\slfam=\thirtnsl \def\sl{\fam\slfam\thirtnsl}%
\textfont\ttfam=\thirtntt \def\tt{\fam\ttfam\thirtntt}%
\textfont\msafam=\thirtnmsa \scriptfont\msafam=\ninemsa
\scriptscriptfont\msafam=\sevenmsa \def\msa{\fam\msafam\thirtnmsa}%
\textfont\msbfam=\thirtnmsb \scriptfont\msbfam=\ninemsb
\scriptscriptfont\msbfam=\sevenmsb \def\msb{\fam\msbfam\thirtnmsb}%
\textfont\eufmfam=\thirtneufm  \scriptfont\eufmfam=\nineeufm
\scriptscriptfont\eufmfam=\seveneufm \def\frak{\fam\eufmfam\teneufm}%
\normalbaselineskip=16pt%
\setbox\strutbox=\hbox{\vrule height11.5pt depth4.5pt width0pt}%
\normalbaselines\rm}%
%
%
\def\Bbb#1{{\msb#1}}

%

\def\re{\Bbb R}
%
\mathchardef\plussquare="0\hexa01
\mathchardef\nge="3\hexb0B
\mathchardef\maltesecross="0\hexa7A
\mathchardef\del="0\hexf01
%
%
%
%
\font\sc=cmcsc10
%
%
%
%
\def\sqr#1#2{{\vcenter{\vbox{\hrule  height.#2truept
	\hbox{\vrule width.#2truept height#1truept 
	\kern#1truept \vrule width.#2truept}
	\hrule height.#2truept}}}}
\def\sq{\sqr55}    
%
%
%
%
\newcount\sectionnumber            
\newcount\resultnumber             
\sectionnumber=0\resultnumber=1    
%
%
%
\def\section#1{\global\advance\sectionnumber by 1
\xdef\nextkey{\number\sectionnumber}
\vskip-\lastskip\penalty-800\vskip 20pt plus10pt minus5pt 
{\large\bf\number\sectionnumber\quad#1}         
\vskip 8pt plus4pt minus4pt
\nobreak\resultnumber=1}      
%
%
%
%
%
         
%
%
%
%

%
\def\proc#1{\xdef\nextkey{\number\sectionnumber.\number\resultnumber}%
\vskip-\lastskip\ppar\bf%
\noindent#1\ \number\sectionnumber.\number\resultnumber
\stdspace\sl\global\advance\resultnumber by 1\ignorespaces}
\def\endproc{\rm\ppar} 
%
%
\def\endprf{\unskip\stdspace\hbox{}
\hfill$\sq$\par\medskip}                 
\def\proof#1{\vskip-\lastskip\ppar\noindent{\bf#1}%
\stdspace\rm\ignorespaces}        
\let\endproof\endprf              
%
%
%
%
%
%
%
%
\def\proclaim#1{\vskip-\lastskip\ppar\bf%
\noindent#1\stdspace\sl\ignorespaces} 
\let\endproclaim\endproc
%
%
%
%

%
%
%
%
%
%
\def\label{\xdef\nextkey{\number\sectionnumber.\number\resultnumber}%
\number\sectionnumber.\number\resultnumber
\global\advance\resultnumber by 1}
%
%
%
%
%
%
%
%
%
%
%
%
%
%
%
%
\newcount\refnumber              
\refnumber=1                     
\long\def\reflist#1\endreflist{%
\long\def\thereflist{#1}{\def\refkey##1##2\par{\xdef##1{\number\refnumber}%
\global\advance\refnumber by 1}%
\def\key##1##2\par{\expandafter\xdef%
\csname##1\endcsname{\number\refnumber}%
\global\advance\refnumber by 1}#1\par}}
\long\def\references{%
\penalty-800\vskip-\lastskip\vskip 15pt plus10pt minus5pt 
{\large\bf References}\ppar 
{\leftskip=25pt\frenchspacing    
\small\parskip=3pt plus2pt       
\def\refkey##1##2\par{\noindent  
\llap{[##1]\stdspace}\ignorespaces##2\par}         
\def\key##1##2\par{\noindent  
\llap{[\ref{##1}]\stdspace}\ignorespaces##2\par}  
\def\,{\thinspace}\thereflist\par}}
%
%
%
\newcount\footnotenumber         
\footnotenumber=1                
\def\fnote#1{\xdef\nextkey{\number\footnotenumber}%
{\small\ifnum\footnotenumber>9\parindent=14pt%
\else\parindent=10pt\fi\footnote{$^{\number\footnotenumber}$}%
{\hglue-5pt#1}\global\advance\footnotenumber by 1}}
%
%
%
%
%
%
%
\newcount\figurenumber          
\figurenumber=1                 
\def\caption#1{\xdef\nextkey{\number\figurenumber}%
\cl{\small Figure \number\figurenumber: #1}%
\global\advance\figurenumber by 1}
\def\figurelabel{\xdef\nextkey{\number\figurenumber}%
\cl{\small Figure \number\figurenumber}%
\global\advance\figurenumber by 1}
\long\def\figure#1\endfigure{{\xdef\nextkey{\number\figurenumber}%
\let\captiontext\relax\def\caption##1{\xdef\captiontext{##1}}%
\midinsert\cl{\ignorespaces#1\unskip\unskip\unskip\unskip}\vglue6pt\cl{\small 
Figure \number\figurenumber\ifx\captiontext\relax\else: \captiontext
\fi}\endinsert\global\advance\figurenumber by 1}}
%
%
%
%
%
%
%
\def\nextkey{??}   
%
\def\key#1{\expandafter\xdef\csname #1\endcsname{\nextkey}}
\def\ref#1{\expandafter\ifx\csname #1\endcsname\relax
\immediate\write16{Reference {#1} undefined}??\else
\csname #1\endcsname\fi}
%
%
%
%
%
%
%
\newread\gtinfile
\newwrite\gtreffile
\def\useforwardrefs{
\openin\gtinfile\jobname.ref
\ifeof\gtinfile
\closein\gtinfile
\immediate\write16{No file \jobname.ref}
\else
\closein\gtinfile
\fi
\immediate\openout\gtreffile \jobname.ref
%
%
\def\key##1{{\def\\{\noexpand}%
\expandafter\xdef\csname ##1\endcsname{\nextkey}%
\immediate\write\gtreffile{\\\expandafter\\\def\\\csname ##1\\\endcsname%
{\nextkey}}}}
%
%
\long\def\reflist##1\endreflist{%
\long\def\thereflist{##1}{\def\refkey####1####2\par{\xdef####1{%
\number\refnumber}{\def\\{\noexpand}\immediate\write\gtreffile
{\\\def\\####1{\number\refnumber}}}\global\advance\refnumber by 1}%
\def\key####1####2\par{\expandafter\xdef%
\csname####1\endcsname{\number\refnumber}%
{\def\\{\noexpand}\immediate\write\gtreffile
{\\\expandafter\\\def\\\csname ####1\\\endcsname{\number\refnumber}}}
\global\advance\refnumber by 1}##1\par}}
\long\def\biblio##1\endbiblio{\reflist##1\endreflist\references}%
%
%
\def\numkey##1{{\def\\{\noexpand}%
\xdef##1{\number\sectionnumber.\number\resultnumber}
\immediate\write\gtreffile{\\\def\\##1%
{\number\sectionnumber.\number\resultnumber}}}}
\def\seckey##1{{\def\\{\noexpand}\xdef##1{\number\sectionnumber}
\immediate\write\gtreffile{\\\def\\##1{\number\sectionnumber}}}}
\def\figkey##1{\xdef##1{\number\figurenumber}%
{\def\\{\noexpand}\immediate\write\gtreffile%
{\\\def\\##1{\number\figurenumber}}}
\number\figurenumber\global\advance\figurenumber by 1}
}   
%
%
%
%
\def\figkey#1{\xdef#1{\number\figurenumber}%
\number\figurenumber\global\advance\figurenumber by 1}
\def\fig#1#2\endfig{%
\midinsert\cl{#2}\vglue6pt\cl{\small Figure #1}\endinsert}
\def\newfig{\number\figurenumber\global\advance\figurenumber by 1}
\def\numkey#1{\xdef#1{\number\sectionnumber.\number\resultnumber}}
\def\seckey#1{\xdef#1{\number\sectionnumber}}
%
%
%
%
%
%
%
%
%
\def\verb{\catcode`\"=\active}       
\def\brev{\catcode`\"=12}            
\brev                                
\verb                                
{\obeyspaces\gdef {\ }}              
{\catcode`\`=\active\gdef`{\relax\lq}}
\def"{%
\begingroup\baselineskip=12pt\def\par{\leavevmode\endgraf}%
\tt\obeylines\obeyspaces\parskip=0pt\parindent=0pt%
\catcode`\$=12\catcode`\&=12\catcode`\^=12\catcode`\#=12%
\catcode`\_=12\catcode`\~=12%
\catcode`\{=12\catcode`\}=12\catcode`\%=12\catcode`\\=12%
\catcode`\`=\active\let"\endgroup}
\brev      
%
%
%
%
%
%
\def\item#1{\par\leavevmode\llap{#1\stdspace}%
\ignorespaces}                             
%
%

%
%
\def\np{\vfil\eject}         
\def\cl{\centerline}         
%
%
%

%
%
%
%
%
\def\title#1{\def\thetitle{#1}}

\def\author#1{\edef\previousauthors{\theauthors}
 \ifx\theauthors\relax\def\theauthors{#1}\else
 \def\theauthors{\previousauthors\par#1}\fi}

%
\def\address#1{\edef\previousaddresses{\theaddress}
 \ifx\theaddress\relax\def\theaddress{#1}\else
 \def\theaddress{\previousaddresses\par\vskip 2pt\par#1}\fi}
\def\secondaddress#1{\edef\previousaddresses{\theaddress}
 \ifx\theaddress\relax\def\theaddress{#1}\else
 \def\theaddress{\previousaddresses\par{\rm and}\par#1}\fi}   

\def\email#1{\edef\previousemails{\theemail}
 \ifx\theemail\relax\def\theemail{#1}\else
 \def\theemail{\previousemails\hskip 0.75em\relax#1}\fi}
\def\secondemail#1{\edef\previousemails{\theemail}
 \ifx\theemail\relax\def\theemail{#1}\else
 \def\theemail{\previousemails\hskip 0.75em{\rm and}\hskip 0.75em
 \relax#1}\fi}
\def\url#1{\edef\previousurls{\theurl}
 \ifx\theurl\relax\def\theurl{#1}\else
 \def\theurl{\previousurls\hskip 0.75em\relax#1}\fi}
\def\secondurl#1{\edef\previousurls{\theurl}
 \ifx\theurl\relax\def\theurl{#1}\else
 \def\theurl{\previousurls\hskip 0.75em{\rm and}\hskip 0.75em
 \relax#1}\fi}
\long\def\abstract#1\endabstract{\long\def\theabstract{#1}}
\def\primaryclass#1{\def\theprimaryclass{#1}}

\def\keywords#1{\def\thekeywords{#1}}
%
%
\let\\\par\let\thetitle\relax
\let\theauthors\relax
\let\theaddress\relax\let\theshortaddress\relax
\let\theemail\relax\let\theurl\relax
\let\theabstract\relax\let\theprimaryclass\relax
\let\thesecondaryclass\relax\let\thekeywords\relax
%
%
%
%
\long\def\maketitlepage{    

\vglue 0.2truein   

%
{\parskip=0pt\leftskip 0pt plus 1fil\def\\{\par\smallskip}{\large
\bf\thetitle}\par\medskip}   

\vglue 0.15truein 

%
{\parskip=0pt\leftskip 0pt plus 1fil\def\\{\par}{\sc\theauthors}
\par\medskip}%
 
\vglue 0.1truein 

%
{\small\parskip=0pt
{\leftskip 0pt plus 1fil\def\\{\par}{\sl\theaddress}\par}
\ifx\theemail\relax\else  
\vglue 5pt \def\\{\stdspace{\rm and}\stdspace} 
\cl{Email:\stdspace\tt\theemail}\fi
\ifx\theurl\relax\else    
\vglue 5pt \def\\{\stdspace{\rm and}\stdspace} 
\cl{URL:\stdspace\tt\theurl}\fi\par}

\vglue 7pt 

{\bf Abstract}

\vglue 5pt

\theabstract

\vglue 7pt 

{\bf AMS Classification numbers}\quad Primary:\quad \theprimaryclass\par

Secondary:\quad \thesecondaryclass

\vglue 5pt 

{\bf Keywords:}\quad \thekeywords

\np  

}    
%
%
\long\def\makeshorttitle{    


%
{\parskip=0pt\leftskip 0pt plus 1fil\def\\{\par\smallskip}{\large
\bf\thetitle}\par\medskip}   

\vglue 0.05truein 

%
{\parskip=0pt\leftskip 0pt plus 1fil\def\\{\par}{\sc\theauthors}
\par\medskip}%
 
\vglue 0.03truein 

%
{\small\parskip=0pt
{\leftskip 0pt plus 1fil\def\\{\par}{\sl\ifx\theshortaddress\relax
\theaddress\else\theshortaddress\fi}\par}
\ifx\theemail\relax\else  
\vglue 5pt \def\\{\stdspace{\rm and}\stdspace} 
\cl{Email:\stdspace\tt\theemail}\fi
\ifx\theurl\relax\else    
\vglue 5pt \def\\{\stdspace{\rm and}\stdspace} 
\cl{URL:\stdspace\tt\theurl}\fi\par}

\vglue 10pt 


{\small\leftskip 25pt\rightskip 25pt{\bf Abstract}\stdspace\theabstract

{\bf AMS Classification}\stdspace\theprimaryclass
\ifx\thesecondaryclass\relax\else; \thesecondaryclass\fi\par
{\bf Keywords}\stdspace \thekeywords\par}
\vglue 7pt
}    
\let\maketitle\makeshorttitle      
%
%

\parskip=0pt


\long\def\ifndef#1{\expandafter\ifx\csname#1\endcsname\relax}

\def\warning#1{\immediate\write16{l.\the\inputlineno: #1}}

\def\section#1#2\par{%
\vskip-\lastskip\penalty-800\vskip 20pt plus10pt minus5pt
{\noindent\large\bf #1\quad #2}%
\vskip 8pt plus4pt minus4pt
\nobreak\noindent}

\def\deftag#1#2{%
\ifndef{#2}\else\warning{tag `#2' defined more than once.}\fi 
\expandafter\def\csname #2\endcsname{#1}}
\def\tag#1{%
\ifndef{#1}\warning{tag `#1' undefined.}{\bf ?}%
\else\csname #1\endcsname\fi}
\def\Tag#1{\tag{#1}}
\def\eqtag#1{(\tag{#1})}
\def\eqTag#1{\hbox{\eqtag{#1}}}
\def\defcite#1#2{\deftag{#1}{#2}}
\def\cite#1{[\tag{#1}]}
\def\ecite#1#2{[\tag{#1}, #2]}

\def\beginproclaim#1. {\proclaim{#1.}}
\def\endproclaim{\endproc}

\def\beginproof#1{\proof{#1}}
\def\endproof{\endprf}

\def\bbb{\Bbb}

\def\Rom#1{{\rm\uppercase\expandafter{\romannumeral#1}}}
\def\a{\alpha}		\def\c{\gamma}
  \def\e{\varepsilon}

\def\bibitem#1{\key{\tag{#1}}}

\def\setminus{-}
\def\Z{{\bbb Z}} 
\def\R{{\bbb R}} 
\def\C{{\bbb C}} 
\let\sset=\subseteq		%
\def\setminus{\mathop{\bbb \char"72}\nolimits}

\let\del=\partial
\let\isomorphism=\cong		
\let\to=\rightarrow

\def\re{\mathop{\rm Re}\nolimits}

\def\({\left(} \def\){\right)}
\def\[{\left[} \def\]{\right]}

\def\braid{Z}
\def\artin{{\cal A}}
\def\artindn{\artin(D_n)}
\def\artinDn{\artin(\tilde{D}_n)}
\def\setminus{-}
\def\D{\Delta}
\def\nguyen{Nguy{\^e}{\~n}}
\def\Wbar{\overline{W}}

%
%
\deftag{1}{sec-intro}
\deftag{1.1}{eq-an-1-diagram}
\deftag{1.1}{tab-results}
\deftag{1.2}{eq-move-for-index-2-strand}
\deftag{1.3}{eq-labelled-dn}
\deftag{1.1}{thm-artindn}
\deftag{1.4}{eq-labelled-Dn}
\deftag{1.5}{eq-braid-generators-for-Dn}
\deftag{1.2}{thm-artinDn}
\deftag{2}{sec-dn}
\deftag{2.1}{eq-simple-roots-for-dn}
\deftag{2.2}{eq-def-of-tau}
\deftag{2.3}{eq-defn-of-w}
\deftag{2.1}{fig-picture-of-Delta}
\deftag{2.1}{thm-fundamental-elt}
\deftag{2.2}{fig-D-h1-Dinverse}
\deftag{2.3}{fig-braid-conjugation-odd}
\deftag{2.4}{fig-braid-conjugation-even}
\deftag{2.2}{thm-Bn(k)-for-odd-n}
\deftag{3}{sec-Dn}
\deftag{3.1}{thm-hyperplane-complement-for-Dn}
\deftag{4}{sec-other-artin-groups}
\defcite{1}{brieskorn:fundamentalgruppe}
\defcite{2}{brieskorn:groupes-de-tresses}
\defcite{3}{brieskorn:artin-gruppen-und-coxeter-gruppen}
\defcite{4}{charney:artin-groups-finite-type-biautomatic}
\defcite{5}{charney:finite-k-pi-1's-for-artin-groups}
\defcite{6}{ATLAS}
\defcite{7}{deligne:immuebles-des-groupes-de-tresses}
\defcite{8}{dbae:wd-proc-gps}
\defcite{9}{fox-neuwirth:braid-groups}
\defcite{10}{nguyen:regular-orbits-affine-weyl-groups}
\defcite{11}{salvetti:homotopy-type-of-artin-groups}
\defcite{12}{tom-dieck:symmetrische-brucken}
\defcite{13}{tom-dieck:rooted-cylinder-ribbons}
\defcite{14}{van-der-lek:extended-artin-groups}
\defcite{15}{van-der-lek:thesis}


\title	{Braid pictures for Artin groups}
\author {Daniel Allcock}
\address{Department of Mathematics, Harvard University, Cambridge, MA 02138}
\email  {allcock@math.harvard.edu}
\url	{http://www.math.harvard.edu/$\sim$allcock}

\abstract

We define the braid groups of a two-dimensional orbifold and introduce
conventions for drawing braid pictures. We use these to realize the
Artin groups associated to the spherical Coxeter diagrams $A_n$,
$B_n=C_n$ and $D_n$ and the affine diagrams $\tilde{A}_n$, $\tilde{B}_n$, $\tilde{C}_n$ and
$\tilde{D}_n$ as subgroups of the braid groups of various simple
orbifolds. The cases $D_n$, $\tilde{B}_n$, $\tilde{C}_n$ and $\tilde{D}_n$ are new.  In each
case the Artin group is a normal subgroup with abelian quotient; in
all cases except $\tilde{A}_n$ the quotient is finite. We illustrate the value
of our braid calculus by performing with pictures a nontrivial
calculation in the Artin groups of type $D_n$.

\endabstract
\primaryclass{20F36}
\keywords{braid groups, Artin groups}
\maketitlepage

\parindent=20pt

\section{\Tag{sec-intro}}{Introduction}

The purpose of this paper is to establish and explain a very close
connection between two quite different-seeming generalizations
of the classical braid group. One generalization is to the class 
of Artin groups. This generalization is very natural
when one studies the braid group from the point of view of singularity
theory, or Lie theory, or the theory of reflection groups. We
will be concerned mainly with the Artin groups $\artindn$ and $\artinDn$
associated to the spherical and affine Coxeter diagrams
$$
\beginpicture
\put {\epsfbox[27 676  207 724]{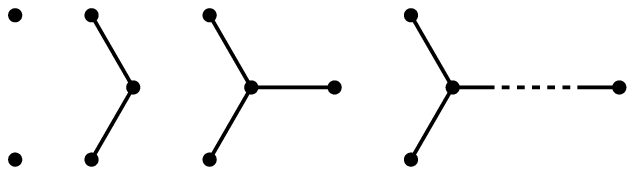}} [bl] at 0 0
\setcoordinatesystem units <1bp,1bp> point at 27 676
\put {\strut$D_2$} [t] <0pt,-1pt> at 30.0 677.215
\put {\strut$D_3$} [t] <0pt,-1pt> at 58.0 677.215
\put {\strut$D_4$} [t] <0pt,-1pt> at 104.0 677.215
\put {\strut$D_n$} [t] <0pt,-1pt> at 174.0 677.215
\endpicture
$$
$$
\beginpicture
\put {\epsfbox[13 676  261 724]{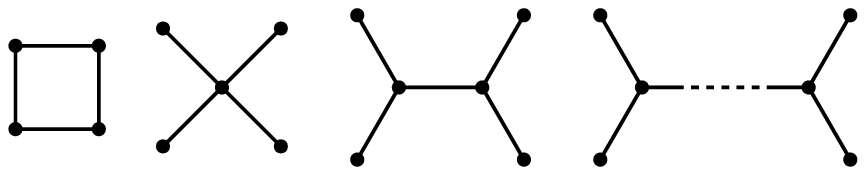}} [bl] at 0 0
\setcoordinatesystem units <1bp,1bp> point at 13 676
\put {\strut$\tilde{D}_3$} [t] <0pt,0pt> at 30.0 677.215
\put {\strut$\tilde{D}_4$} [t] <0pt,0pt> at 77.4558 677.215
\put {\strut$\tilde{D}_5$} [t] <0pt,0pt> at 140.426 677.215
\put {\strut$\tilde{D}_n$} [t] <0pt,0pt> at 222.426 677.215
\endpicture
$$
but we will also consider the diagrams associated to the other
`classical' Coxeter diagrams. These are the spherical diagrams $A_n$
and $B_n=C_n$ and the affine diagrams $\tilde{A}_n$, $\tilde{B}_n$ and $\tilde{C}_n$. 
(Each diagram $x_n$ has $n$ nodes and each diagram $X_n$  has $n+1$. The
remaining diagrams are given in section~\tag{sec-other-artin-groups}.)
The most direct way to define $\artin(D_n)$ and $\artin(\tilde{D}_n)$ is to give
generators and relations: one takes one generator for each node
of the diagram and imposes the relations that two of these
generators commute (resp. braid) if the corresponding nodes are
unjoined (resp. joined). When we say that two group elements $x$
and $y$ braid, we mean that they satisfy $xyx=yxy$. One can
obtain the classical braid group on $n$ strands by applying this
construction to the Coxeter diagram $A_{n-1}$ given below.
$$
\beginpicture
\put {\epsfbox[27 697  105 703]{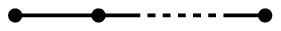}} [bl] at 0 0
\setcoordinatesystem units <1bp,1bp> point at 27 697
\endpicture
\eqno\eqTag{eq-an-1-diagram}
$$
It turns out that this procedure of assigning a group to a Coxeter
diagram is not just a random construction---the generators and
relations arise from natural geometric considerations and are closely
related to the associated Coxeter group. However, despite being
natural objects, the Artin groups are still somewhat mysterious. The
aim of this paper is to provide a very concrete way to understand
some of them.

The other generalization of the classical braid group is the braid
group of a two-dimensional orbifold. The braid groups of certain
two-orbifolds turn out to contain the Artin groups $\artindn$ and
$\artinDn$ as subgroups of very small index. The relevant orbifolds
are among the simplest possible ones: one is the plane with a single
cone point of order 2 and the other is the plane with two such cone
points. In sections~\tag{sec-dn} and~\tag{sec-Dn} we treat these
groups in detail. There are similar results, some already known, for
the Artin groups associated to the other classical Coxeter diagrams;
we summarize these in table~\tag{tab-results} and will discuss them in
section~\tag{sec-other-artin-groups}. 
The cases $D_n$, $\tilde{B}_n$, $\tilde{C}_n$
and $\tilde{D}_n$ are new,
the case $A_{n-1}$ is classical, $B_n$
is treated in \cite{brieskorn:groupes-de-tresses}, and the result for
$\tilde{A}_{n-1}$ appears in \cite{tom-dieck:rooted-cylinder-ribbons}. Our
interpretation of these Artin groups as braid groups leads to several
curious coincidences, which we will discuss in
section~\tag{sec-other-artin-groups}.

\topinsert
\hfil\vbox{\halign{%
\hfil#\hfil&\quad\hfil#\hfil&\quad\hfil#\hfil&\quad\hfil#\hfil\cr
{\bf diagram} & {\bf orbifold features}
& {\bf quotient group} & {\bf condition}\cr
\noalign{\vskip2pt}
$A_{n-1}$ & none & 1 & $n>1$ \cr
$B_n$, $C_n$ & 1 puncture & 1 & $n>1$ \cr
$D_n$ & 1 cone point & $\Z/2$ & $n>1$ \cr
$\tilde{A}_{n-1}$ & 1 puncture & $\Z$ & $n>2$ \cr
$\tilde{B}_n$ & 1 puncture, 1 cone point & $\Z/2$ & $n>2$ \cr
$\tilde{C}_n$ & 2 punctures & 1 & $n>1$ \cr
$\tilde{D}_n$ & 2 cone points & $\Z/2\times\Z/2$ & $n>2$\cr
}}\hfil
\medskip
\narrower\noindent
\small
Table~\Tag{tab-results}. Summary of results. Entries should be interpreted as
follows. For each row and each number $n$ satisfying the condition in
the last column, the Artin group associated to the given diagram is a
normal subgroup of the $n$-strand braid group of the orbifold which is
the plane equipped with the given features. All of the cone points
indicated have order 2. The quotient of the braid group by the Artin
group is given in the third column, and in each case the group
extension splits.
\endinsert

To define the braid groups, let $L$ be any two-dimensional
orbifold. Then its ($n$-strand) pure braid space is
$L^n\setminus\D_n$, where 
$$
\D_n=
\{\,(x_1,\ldots,x_n)\in L^n
\>|\>
\hbox{$x_i=x_j$ for some $i\neq j$}\,\}.
$$
The symmetric group $S_n$ acts on $L^n$ in the obvious way; the
action is free on $L^n\setminus\D_n$ and defines a
covering map to $X_n=(L^n\setminus\D_n)/S_n$, which we call the
($n$-strand) braid space of $L$. We fix a basepoint
$b=(b_1,\ldots,b_n)$ for $X_n$, chosen so that it does not lie
on the orbifold locus. Then the $n$-strand braid group
$\braid_n=\braid_n(L)$ is defined to be the orbifold fundamental group of
$X_n$.

We are interested in the case where the only orbifold
singularities of $L$ are cone points; in this case we may understand
elements of $\pi_1(X_n)$ as follows. Each element of $\braid_n$ may
be represented by a set of paths $\c_1,\ldots,\c_n:[0,1]\to L$
that miss the cone points of $L$ and have the following
properties. First, for some permutation $\pi$ of
$\{1,\ldots,n\}$, each $\c_i$ begins at $b_i$ and ends at
$b_{\pi(i)}$. Second, for each $t\in[0,1]$, the points $\c_i(t)$
are all distinct. We call any  $n$-tuple of paths satisfying
these conditions a braid. Then $\c:t\mapsto(\c_1(t),\ldots,\c_n(t))$ is
a path in $L^n\setminus\D_n$. If $\pi$ is nontrivial then $\c$
will not be a closed path, but its image in $X_n$ will be, and
hence define an element of $\braid_n$.

When $L$ is the complex plane with some punctures and/or cone
points then we can draw pictures of braids in essentially the
same manner as in the classical case. One simply draws a picture
of $L\times[0,1]$ and plots the curves
$t\mapsto(\c_i(t),t)$. One usually pictures the slices
$L\times\{t\}$ as horizontal planes with points having large
imaginary part appearing far from the viewer, and one regards
$t$ as increasing in the downward direction. If $L$ is just $\C$
with no cone points then we recover braid pictures in the usual
sense, for example:
$$
\beginpicture
\put {\epsfbox[8 331  115 392]{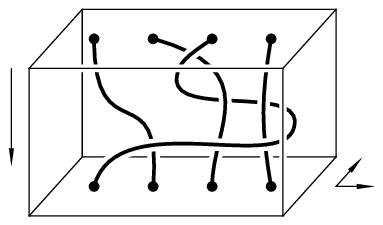}} [bl] at 0 0
\setcoordinatesystem units <1bp,1bp> point at 8 331
\put {$t$ } [r] <0pt,0pt> at 10.2 361.25
\put {\strut$\C$} [t] <0pt,0pt> at 103.7 340.0
\endpicture
$$
Of course we will typically draw less elaborate pictures. When
$L$ has cone points then we must indicate them in the picture of
$L\times[0,1]$. We do this by drawing each segment
$\{c\}\times[0,1]$, with $c$ a cone point, as a thick vertical
line with the order of the cone point indicated nearby. Here is
an example of a braid when $L$ has two cone points, of orders
$3$ and $4$:
$$
\beginpicture
\put {\epsfbox[29 334  107 422]{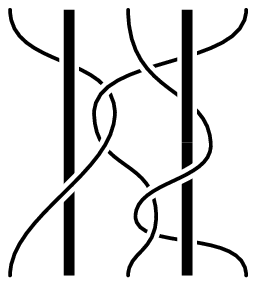}} [bl] at 0 0
\setcoordinatesystem units <1bp,1bp> point at 29 334
\put {\strut\lower1.5pt\hbox{$3$}} [t] <0pt,0pt> at 51.0 340.0
\put {\strut\lower1.5pt\hbox{$4$}} [t] <0pt,0pt> at 85.0 340.0
\endpicture
$$
Finally, one multiplies braids in the obvious way: if $B$ and
$B'$ are braids then we represent $BB'$ by simply placing a
picture of $B'$ below a picture of $B$.

As usual, two braids represent the same element of $\braid_n$ if one
may be deformed to the other through a continuum of
braids. Also, when $L$ has cone points then there is an
additional `move' that can be performed when a braid strand links
the `strand' representing a cone point of order $p$ exactly $p$ times:
$$
\beginpicture
\put {\epsfbox[29 334  138 405]{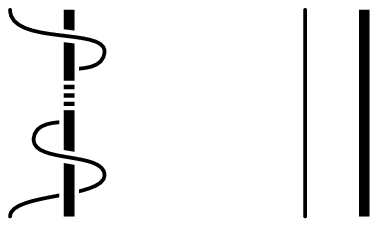}} [bl] at 0 0
\setcoordinatesystem units <1bp,1bp> point at 29 334
\put {\llap{($2p$ crossings)\quad}} [r] <0pt,0pt> at 34.0 369.75
\put {\strut$p$} [t] <0pt,0pt> at 51.0 340.0
\put {$=$} [c] <0pt,0pt> at 85.0 369.75
\put {\strut$p$} [t] <0pt,0pt> at 136.0 340.0
\put {\vrule width 0pt height 0pt depth 0pt{}} [c] <0pt,0pt> at 153.0 369.75
\endpicture
$$
These two braids represent the same element of $\pi_1(X_n)$
because a loop in $L$ encircling an order $p$ cone point exactly
$p$ times represents the trivial element of the orbifold
fundamental group of $L$. When $p=2$ the most useful formulation
of this move is as
$$
\vcenter{%
\beginpicture
\put {\epsfbox[29 334  150 379]{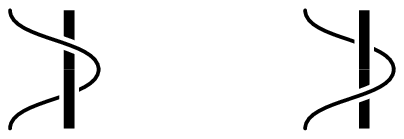}} [bl] at 0 0
\setcoordinatesystem units <1bp,1bp> point at 29 334
\put {$=$} [c] <0pt,0pt> at 85.0 357.0
\put {\vrule width 0pt height 0pt depth 0pt{}} [c] <0pt,0pt> at 153.0 357.0
\endpicture
}
\eqno\eqTag{eq-move-for-index-2-strand}
$$
We have suppressed the numeral `2' that would denote the order of the
cone point. Since all heavy lines up until section~\tag{sec-other-artin-groups} represent
cone points of order 2, we will continue to do this.

We can now state the main results of the paper. We take $n\geq2$
and label the standard generators of $\artindn$ by
$$
\vcenter{%
\beginpicture
\put {\epsfbox[12 671  93 729]{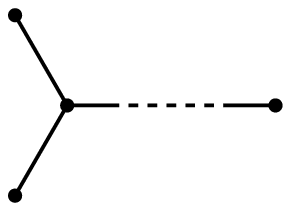}} [bl] at 0 0
\setcoordinatesystem units <1bp,1bp> point at 12 671
\put {$g_3\;$} [r] <0pt,0pt> at 30.0 700.0
\put {$g_2\,$} [r] <0pt,0pt> at 13.0 674.019
\put {$g_1\,$} [r] <0pt,0pt> at 13.0 725.981
\put {$\,g_n$} [l] <0pt,0pt> at 92.0 700.0
\endpicture
}
\eqno\eqTag{eq-labelled-dn}
$$
Let $k$ be the orbifold which is $\C$ with a single cone point of
order 2 and write $h_1,\ldots,h_n$ for the
$n$-strand braids in $k$ given by the following diagrams:
$$
\beginpicture
\put {\epsfbox[29 334  290 379]{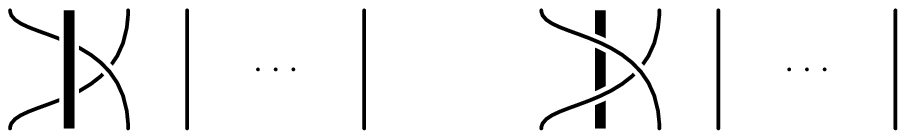}} [bl] at 0 0
\setcoordinatesystem units <1bp,1bp> point at 29 334
\put {\strut$h_1$} [t] <0pt,-1pt> at 85.0 340.0
\put {\vrule width 0pt height 0pt depth 0pt} [c] <0pt,0pt> at 17.0 357.0
\put {\strut$h_2$} [t] <0pt,-1pt> at 238.0 340.0
\put {\vrule width 0pt height 0pt depth 0pt} [c] <0pt,0pt> at 306.0 357.0
\endpicture
$$
$$
\beginpicture
\put {\epsfbox[33 334  294 379]{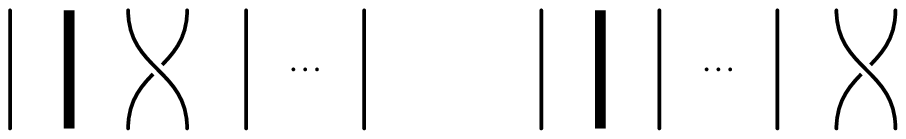}} [bl] at 0 0
\setcoordinatesystem units <1bp,1bp> point at 33 334
\put {\strut$h_3$} [t] <0pt,-1pt> at 85.0 340.0
\put {\vrule width 0pt height 0pt depth 0pt} [c] <0pt,0pt> at 17.0 357.0
\put {\strut$h_n$} [t] <0pt,-1pt> at 238.0 340.0
\put {\vrule width 0pt height 0pt depth 0pt} [c] <0pt,0pt> at 306.0 357.0
\endpicture
$$
It is easy to check that the map $g_i\mapsto h_i$ defines a
homomorphism $\artindn\to \braid_n(k)$. The verification of each
defining relation of $\artindn$ except for the commutativity of
$h_1$ and $h_2$ is easy and standard. The commutativity of
$h_1$ and $h_2$ follows from a double application of the
`orbifold move':
$$ 
\beginpicture
\put {\epsfbox[29 334  243 396]{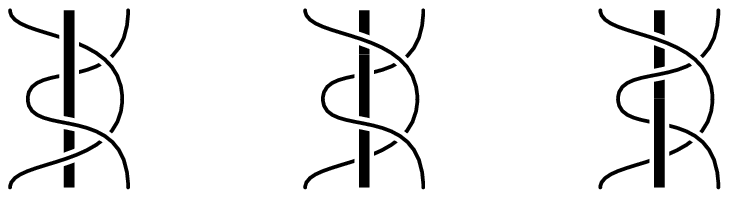}} [bl] at 0 0
\setcoordinatesystem units <1bp,1bp> point at 29 334
\put {$=$} [c] <0pt,0pt> at 93.5 365.5
\put {$=$} [c] <0pt,0pt> at 178.5 365.5
\endpicture
$$
In section~\tag{sec-dn} we prove the following theorem.

\beginproclaim 
Theorem~\Tag{thm-artindn}.  
For $n\geq2$, the map
$\artindn\to \braid_n(k)$ is an isomorphism onto its image, which
has index two in $\braid_n(k)$.  There is a complementary subgroup
$\Z/2$, whose nontrivial element acts on $\artindn$ by
exchanging $g_1$ with $g_2$ and fixing each of the remaining $g_i$.
\endproclaim

Something very similar happens for the affine artin group
$\artinDn$ for $n\geq3$. We label the standard generators by
$$
\vcenter{%
\beginpicture
\put {\epsfbox[10 671  286 729]{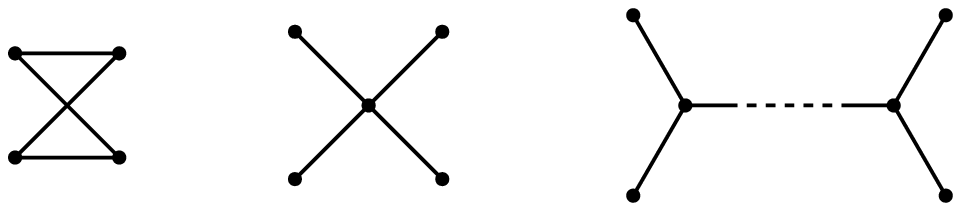}} [bl] at 0 0
\setcoordinatesystem units <1bp,1bp> point at 10 671
\put {$G_1\,$} [r] <0pt,0pt> at 13.0 715.0
\put {$G_2\,$} [r] <0pt,0pt> at 13.0 685.0
\put {$\,G_3$} [l] <0pt,0pt> at 47.0 715.0
\put {$\,G_4$} [l] <0pt,0pt> at 47.0 685.0
\put {$G_1\,$} [r] <0pt,0pt> at 93.6064 721.213
\put {$G_2\,$} [r] <0pt,0pt> at 93.6064 678.787
\put {$G_3\;$} [r] <0pt,0pt> at 114.82 700.0
\put {$\,G_4$} [l] <0pt,0pt> at 140.033 721.213
\put {$\,G_5$} [l] <0pt,0pt> at 140.033 678.787
\put {$G_1\,$} [r] <0pt,0pt> at 191.033 725.981
\put {$G_2\,$} [r] <0pt,0pt> at 191.033 674.019
\put {$G_3\;$} [r] <0pt,0pt> at 206.033 700.0
\put {$\;G_{n-1}$} [l] <0pt,0pt> at 270.033 700.0
\put {$\,G_n$} [l] <0pt,0pt> at 285.033 725.981
\put {$\,G_{n+1}$} [l] <0pt,0pt> at 285.033 674.019
\endpicture
}
\eqno\eqTag{eq-labelled-Dn}
$$
We take $K$ to be the orbifold which is $\C$ with two cone
points of order 2, and write $H_1,\ldots,H_{n+1}$ for the
$n$-strand braids in $K$ given by the diagrams
$$ 
\beginpicture
\put {\epsfbox[29 334  171 379]{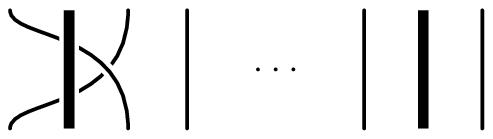}} [bl] at 0 0
\setcoordinatesystem units <1bp,1bp> point at 29 334
\put {\llap{$H_1=$\qquad}} [r] <0pt,0pt> at 34.0 357.0
\put {\vrule width 0pt height 0pt depth 0pt} [c] <0pt,0pt> at 17.0 357.0
\put {\vrule width 0pt height 0pt depth 0pt} [c] <0pt,0pt> at 187.0 357.0
\endpicture
$$
$$ 
\beginpicture
\put {\epsfbox[29 334  171 379]{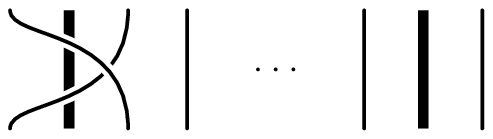}} [bl] at 0 0
\setcoordinatesystem units <1bp,1bp> point at 29 334
\put {\llap{$H_2=$\qquad}} [r] <0pt,0pt> at 34.0 357.0
\put {\vrule width 0pt height 0pt depth 0pt} [c] <0pt,0pt> at 17.0 357.0
\put {\vrule width 0pt height 0pt depth 0pt} [c] <0pt,0pt> at 187.0 357.0
\endpicture
$$
$$ 
\beginpicture
\put {\epsfbox[33 334  171 379]{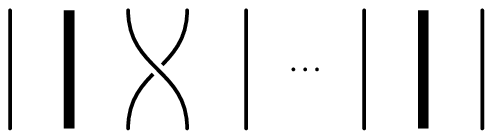}} [bl] at 0 0
\setcoordinatesystem units <1bp,1bp> point at 33 334
\put {\llap{$H_3=$\qquad}} [r] <0pt,0pt> at 34.0 357.0
\put {\vrule width 0pt height 0pt depth 0pt} [c] <0pt,0pt> at 17.0 357.0
\put {\vrule width 0pt height 0pt depth 0pt} [c] <0pt,0pt> at 187.0 357.0
\endpicture
$$
$$ 
\beginpicture
\put {\epsfbox[33 334  171 379]{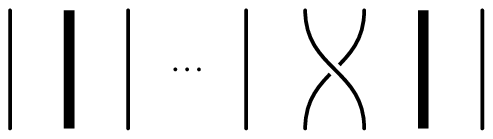}} [bl] at 0 0
\setcoordinatesystem units <1bp,1bp> point at 33 334
\put {\llap{$H_{n-1}=$\qquad}} [r] <0pt,0pt> at 34.0 357.0
\put {\vrule width 0pt height 0pt depth 0pt} [c] <0pt,0pt> at 17.0 357.0
\put {\vrule width 0pt height 0pt depth 0pt} [c] <0pt,0pt> at 187.0 357.0
\endpicture
$$
$$ 
\beginpicture
\put {\epsfbox[33 334  175 379]{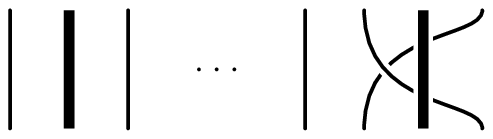}} [bl] at 0 0
\setcoordinatesystem units <1bp,1bp> point at 33 334
\put {\llap{$H_n=$\qquad}} [r] <0pt,0pt> at 34.0 357.0
\put {\vrule width 0pt height 0pt depth 0pt} [c] <0pt,0pt> at 17.0 357.0
\put {\vrule width 0pt height 0pt depth 0pt} [c] <0pt,0pt> at 187.0 357.0
\endpicture
$$
$$ 
\beginpicture
\put {\epsfbox[33 334  175 379]{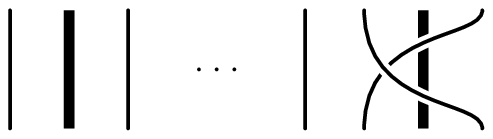}} [bl] at 0 0
\setcoordinatesystem units <1bp,1bp> point at 33 334
\put {\llap{$H_{n+1}=$\qquad}} [r] <0pt,0pt> at 34.0 357.0
\put {\vrule width 0pt height 0pt depth 0pt} [c] <0pt,0pt> at 17.0 357.0
\put {\vrule width 0pt height 0pt depth 0pt} [c] <0pt,0pt> at 187.0 357.0
\endpicture
$$
The third and fourth pictures should be ignored when $n=3$.
By the same argument as in the spherical case, the map
$G_i\mapsto H_i$ defines a homomorphism $\artinDn\to
\braid_n(K)$. 

\beginproclaim
Theorem~\Tag{thm-artinDn}.
For $n\geq3$, the map $\artinDn\to \braid_n(K)$ is an isomorphism onto
its image, which is normal in $\braid_n(K)$ and has index $4$. 
There is a complementary subgroup $\Z/2\times\Z/2$; one element of
this group  acts on
$\artinDn$ be exchanging $H_1$ with $H_2$, fixing each of the
remaining $H_i$, and another acts by exchanging $H_n$ with
$H_{n+1}$, fixing each of the remaining $H_i$.
\endproclaim

We treat the Artin groups of the other classical Coxeter diagrams in
section~\tag{sec-other-artin-groups}. The results and arguments are very similar to the
$D_n$ and $\tilde{D}_n$ cases, so our presentation there is much more
telegraphic. 

The author is supported by a National Science Foundation
Postdoctoral Fellowship.

\section{\Tag{sec-dn}}{The spherical Artin group $\artin(D_n)$}

In this section we will prove theorem~\tag{thm-artindn}. Our
earlier definition of the Artin group in terms of generators and relations
conceals the geometric origin of the group. The group $\artindn$ is
important because it turns out to be the fundamental group of the
space of conjugacy classes of regular semisimple elements of the Lie
algera ${\frak so}_{2n}\C$, and also of the discriminant
complement in the deformation space of the simple singularity
$D_n$. Both of these manifestations of $\artindn$ are closely related
to the appearance of $\artindn$ in the much simpler context of
finite reflection groups, to which we now turn.

Let $V^\R=\R^n$ be equipped
with the standard Euclidean metric, and let $V$ be the
complexification of $V^\R$. The $D_n$ root system is the set of vectors
or `roots' in $V^\R$ obtained from $(\pm1,\pm1,0,\dots,0)$ by
arbitrary permutation of coordinates. The Weyl group $W=W(D_n)$
is the group generated by the reflections in these roots, that is,
across the hyperplanes orthogonal to them. We write $\Sigma$
for the set of reflections in $W$, which are precisely the
reflections in the roots. For each $s\in\Sigma$ we let $H_s^\R$ be
the set of fixed points in $V^\R$ of $s$, which is just the
mirror of the reflection. Then the complexification $H_s$ of
$H_s^\R$ is the set of fixed points of $s$ in $V$. Because each
$H_s$ has real codimension two in $V$, the space
$V_0=V-\cup_{s\in\Sigma}H_s$ is connected. This is called the
pure braid space for $D_n$ and its fundamental group is called
the pure Artin group of type $D_n$. It is known that $W$ acts
freely on $V_0$, and the quotient manifold $V_0/W$ is called the
pure braid space for $D_n$.

It follows from work of Brieskorn \cite{brieskorn:fundamentalgruppe}
that $\pi_1(V_0/W)\isomorphism\artindn$. Indeed his work applies in
greater generality, where one replaces $W(D_n)$ by any other finite
reflection group $W$ and constructs the corresponding space $V_0/W$; then
$\pi_1(V_0/W)$ is the Artin group corresponding to the chosen
reflection group. For the constructions of Artin groups associated to
arbitrary reflection groups, see
\cite{nguyen:regular-orbits-affine-weyl-groups},
\cite{van-der-lek:extended-artin-groups} or
\cite{van-der-lek:thesis}. Also, we treat the classical
affine reflection groups in sections~\tag{sec-Dn}
and~\tag{sec-other-artin-groups}.

\beginproof{Proof of theorem~\tag{thm-artindn}:}
We take $k$ to be the orbifold $\C/(z\mapsto-z)$, and regard the map
$z\mapsto z^2$ as an orbifold covering map from $\C$ to $k$.
Recall
that $X_n=(k^n-\D_n)/S_n$. The idea of the proof is that there is a
two-fold orbifold covering map $V_0/W\to X_n$. We 
begin by describing $V_0$ explicitly:
$$ 
V_0=\{(x_1,\ldots,x_n)\in\C^n
|
\hbox{$x_i\neq\pm x_j$ if $i\neq j$}\}.
$$
Using ATLAS \cite{ATLAS} notation for group structures, the Weyl group
$W$ has structure $2^{n-1}{\,:\,}S_n$, where $S_n$ permutes the
coordinates and the group $2^{n-1}$ acts by changing the signs of any
even number of them. The latter group is a subgroup of the larger
group $2^n$, also normalized by $S_n$, consisting of all sign-changes
of coordinates.

We define an orbifold covering map $\sigma:V\to V/2^n=k^n$ by
$(x_1,\ldots,x_n)\mapsto(x_1^2,\ldots,x_n^2)$. The image of
$V_0$ in $k^n$ is exactly $k^n-\D_n$. Furthermore, $\sigma$
identifies the actions of $S_n$ on its domain and
range. Therefore the orbifolds $(V_0/2^n)/S_n$ and
$(k^n-\D_n)/S_n$ are isomorphic. The latter space is $X_n$ and
the former is $V_0/(2^n{\,:\,}S_n)$. Since $W=2^{n-1}{\,:\,}S_n$ has index
two in $2^n{\,:\,}S_n$, we see that $V_0/W$ is a double cover of
$V_0/(2^n{\,:\,}S_n)=X_n$. This proves that $\artindn$ has index two
in $\braid_n(k)$.

Next we show that this 2-fold covering map induces the stated map
on fundamental groups. For this we use Brieskorn's explicit
description \cite{brieskorn:fundamentalgruppe} of the standard Artin
generators. We choose a set of simple roots
$$
\vcenter{%
\beginpicture
\put {\epsfbox[12 671  93 729]{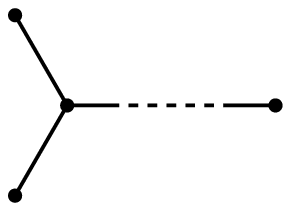}} [bl] at 0 0
\setcoordinatesystem units <1bp,1bp> point at 12 671
\put {\strut$(-1,1,0,\ldots,0)$} [b] <0pt,0pt> at 15.0 727.981
\put {\strut$(1,1,0,\ldots,0)$} [t] <0pt,0pt> at 15.0 672.019
\put {$(0,-1,1,0,\ldots,0)\;$} [r] <0pt,0pt> at 28.0 700.0
\put {$\,(0,\ldots,0,-1,1)$} [l] <0pt,0pt> at 92.0 700.0
\endpicture
}
\eqno\eqTag{eq-simple-roots-for-dn}
$$
for $D_n$ and let $C$ be the (open) Weyl chamber they define. This is
the set of all points in $V^\R$ with positive inner products with each
simple root. We fix a basepoint $p_0$ in $C$, and for each simple root
$s=1,\ldots,n$ with numberings as in \eqtag{eq-labelled-dn}, we let
$p_s$ be the image of $p_0$ under the reflection in that root. Then we
let $L_s$ be the complex line in $V$ containing $p_0$ and $p_s$, and
let $g_s:[0,1]\to L_s$ be the path from $p_0$ to $p_s$ that misses the
mirrors and has the following properties. On $[0,{1\over3}]$ and
$[{2\over3},1]$, $g_s$ is linear in the real line through $p_0$ and
$p_s$, and on $[{1\over3},{2\over3}]$ it is a positively-oriented
semicircle with center $(p_0+p_s)/2$ and very small radius. Then $g_s$
is not a loop but its image in $V_0/W$ is. We have
associated a loop in $V_0/W$ with each node of the diagram, and
Brieskorn's theorem says that we may take these to be the standard
generators.

To understand the image of each of these loops in $X_n$, it suffices
to compute their images under $\sigma:V_0\to k^n\setminus\D_n$, since
the image in $(k^n\setminus\D_n)/S_n$ is obtained simply by
`forgetting' in the usual way the ordering of the strands. First, we
observe that the point $(0,1,\ldots,n-1)$ lies in $C$, and that
Brieskorn's theorem obviously still holds if we choose $p_0$ to be a
perturbation of this point, say $p_0=(\e i,1,2,\ldots,n-1)$
where $\e$ is a small positive number. (Actually, there is no need to
restrict $\e$ to be small.)  Then the paths $g_1,\ldots,g_n$ in $V_0$
are given in by the following figure. For each $s$ we have shown the
union of the coordinate projections of $g_s$.
$$
\beginpicture
\put {\epsfbox[43 335  211 372]{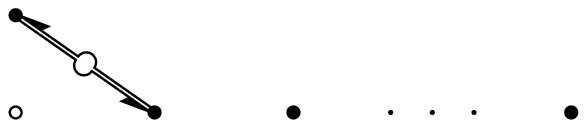}} [bl] at 0 0
\setcoordinatesystem units <1bp,1bp> point at 43 335
\put {\strut$0$} [t] <0pt,-1pt> at 48.0 338.0
\put {\strut$1$} [t] <0pt,-1pt> at 88.0 338.0
\put {\strut$2$} [t] <0pt,-1pt> at 128.0 338.0
\put {\strut$n-1$} [t] <0pt,-1pt> at 208.0 338.0
\put {$\,\e i$} [r] <-1pt,0pt> at 46.0 368.0
\put {$g_1$\qquad} [r] <0pt,0pt> at 8.0 340.0
\endpicture
$$
\medskip
$$
\beginpicture
\put {\epsfbox[5 308  211 371]{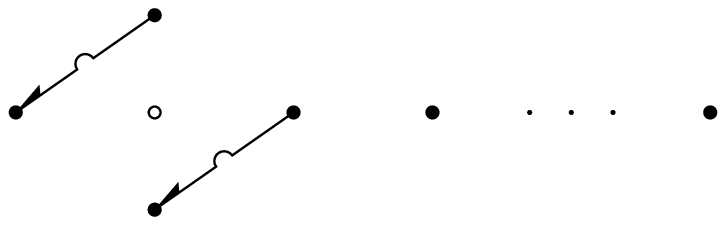}} [bl] at 0 0
\setcoordinatesystem units <1bp,1bp> point at 5 308
\put {\strut$-1$} [t] <0pt,-1pt> at 8.0 338.0
\put {$-\e i\,$} [r] <-1pt,0pt> at 46.0 312.0
\put {$g_2$\qquad} [r] <0pt,0pt> at 8.0 354.0
\put {\phantom{\strut$n-1$}} [b] <0pt,0pt> at 208.0 342.0
\endpicture
$$
\medskip
$$
\beginpicture
\put {\epsfbox[43 333  211 371]{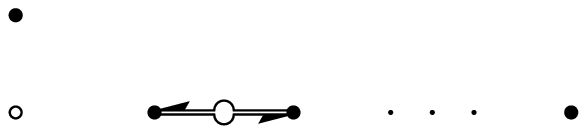}} [bl] at 0 0
\setcoordinatesystem units <1bp,1bp> point at 43 333
\put {$g_3$\qquad} [r] <0pt,0pt> at 8.0 354.0
\put {\phantom{\strut$n-1$}} [b] <0pt,0pt> at 208.0 342.0
\endpicture
$$
\medskip
$$
\beginpicture
\put {\epsfbox[43 333  211 371]{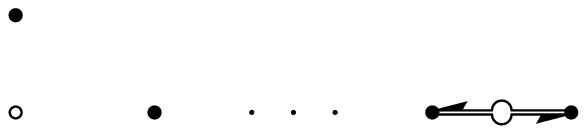}} [bl] at 0 0
\setcoordinatesystem units <1bp,1bp> point at 43 333
\put {$g_n$\qquad} [r] <0pt,0pt> at 8.0 354.0
\put {\strut$n-2$} [b] <0pt,0pt> at 168.0 342.0
\put {\strut$n-1$} [b] <0pt,0pt> at 208.0 342.0
\endpicture
$$
The images of the $g_s$
under $\sigma$ are given below; they are obtained by applying
the squaring map to the paths above.
$$ 
\beginpicture
\put {\epsfbox[7 445  261 461]{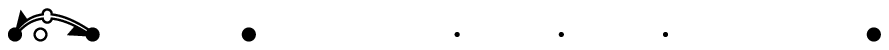}} [bl] at 0 0
\setcoordinatesystem units <1bp,1bp> point at 7 445
\put {$-\e^2\,$} [r] <0pt,0pt> at 8.65 450.0
\put {\strut$0$} [t] <0pt,-1pt> at 18.0 448.0
\put {\strut$1$} [t] <0pt,-1pt> at 33.0 448.0
\put {\strut$4$} [t] <0pt,-1pt> at 78.0 448.0
\put {\strut$(n-1)^2$} [t] <0pt,-1pt> at 258.0 448.0
\put {$g_1$\qquad\qquad} [r] <0pt,0pt> at 10.65 450.0
\endpicture
$$
$$ 
\beginpicture
\put {\epsfbox[7 439  261 455]{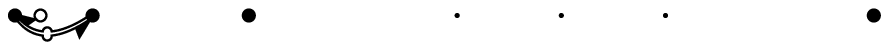}} [bl] at 0 0
\setcoordinatesystem units <1bp,1bp> point at 7 439
\put {$g_2$\qquad\qquad} [r] <0pt,0pt> at 10.65 450.0
\put {\phantom{\strut$(n-1)^2$}} [b] <0pt,0pt> at 258.0 452.0
\endpicture
$$
$$ 
\beginpicture
\put {\epsfbox[8 443  261 457]{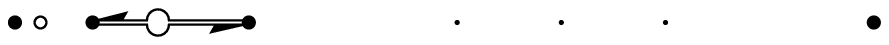}} [bl] at 0 0
\setcoordinatesystem units <1bp,1bp> point at 8 443
\put {$g_3$\qquad\qquad} [r] <0pt,0pt> at 10.65 450.0
\put {\phantom{\strut$(n-1)^2$}} [b] <0pt,0pt> at 258.0 452.0
\endpicture
$$
$$ 
\beginpicture
\put {\epsfbox[8 438  262 462]{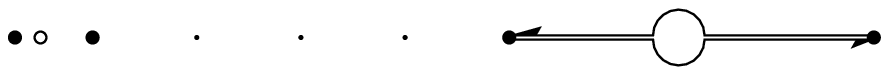}} [bl] at 0 0
\setcoordinatesystem units <1bp,1bp> point at 8 438
\put {$g_n$\qquad\qquad} [r] <0pt,0pt> at 10.65 450.0
\put {\strut$(n-2)^2$} [b] <0pt,0pt> at 153.0 452.0
\put {\strut$(n-1)^2$} [b] <0pt,0pt> at 258.0 452.0
\endpicture
$$

\smallskip\noindent Drawing these braids using the conventions of
section~\tag{sec-intro} shows that the images of the $g_i$ are the
$h_i$. Here, the basepoint $b$ for $k^n\setminus\D_n$ is
$(-\e^2,1,4,\ldots,(n-1)^2)$.

Now we show that $\braid_n(k)$ has the claimed semidirect product
structure. First we observe that the braid
$$ 
\vcenter{%
\beginpicture
\put {\epsfbox[29 334  107 379]{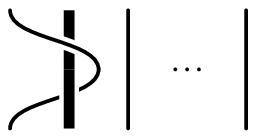}} [bl] at 0 0
\setcoordinatesystem units <1bp,1bp> point at 29 334
\put {$\tau=$} [r] <0pt,0pt> at 17.0 357.0
\endpicture
}
\eqno\eqTag{eq-def-of-tau}
$$
has order 2. Furthermore, it does not lie in the group $\artindn$
generated by the $h_i$. This follows from the fact that $V_0/W$ is
aspherical (\cite{brieskorn:groupes-de-tresses},
\cite{deligne:immuebles-des-groupes-de-tresses}) and hence its
fundamental group $\artindn$ is torsion-free. A more elementary way to
see that $\tau\notin\artindn$ is to note that its lift to $V_0$ joins
the basepoint $p_0$ to $(-\e i,1,2,\ldots,n-1)$ and that these two
points are inequivalent under $W$. It is easy to check that
conjugation by $\tau$ induces the stated automorphism of $\artindn$.
\endproof

Theorem~\tag{thm-artindn} justifies the use of our braid calculus to
perform calculations in $\artindn$. As an example of a nontrivial
application of this calculus, we will discuss the
`fundamental element' $\D$ described in
\cite{charney:artin-groups-finite-type-biautomatic} and
\cite{brieskorn:artin-gruppen-und-coxeter-gruppen}. From our
perspective this is the braid $w^{n-1}$, where $w=h_1h_2\cdots
h_n$ is given in pictures by
$$
\vcenter{%
\beginpicture
\put {\epsfbox[29 334  260 413]{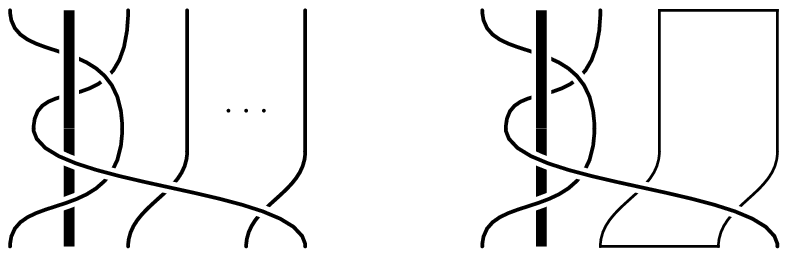}} [bl] at 0 0
\setcoordinatesystem units <1bp,1bp> point at 29 334
\put {$w\;=$} [r] <0pt,0pt> at 23.8 374.0
\put {$=$} [c] <0pt,0pt> at 147.9 374.0
\endpicture
}
\eqno\eqTag{eq-defn-of-w}
$$
In the last picture we have used a ribbon to indicate some
number of strands moving in parallel. This makes some
pictures easier to understand. By stacking together $n-1$
copies of $w$, one sees that $\D$ is represented
by the braid of fig.~\tag{fig-picture-of-Delta}.

\topinsert
\hfil
\beginpicture
\put {\epsfbox[23 334  128 473]{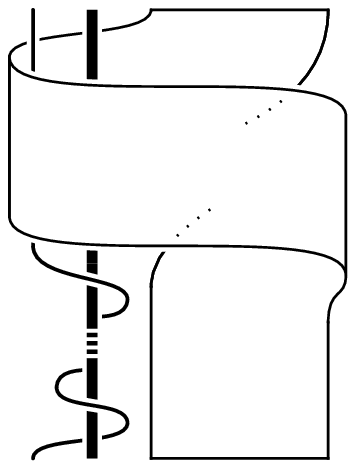}} [bl] at 0 0
\setcoordinatesystem units <1bp,1bp> point at 23 334
\put {$\Delta\;=$} [l] <0pt,0pt> at -51.0 404.6
\put {($(2n-2)$ crossings)} [r] <0pt,0pt> at 37.4 370.6
\endpicture
\hfil

\line{\small\hfil Figure~\Tag{fig-picture-of-Delta}. 
The fundamental element $\Delta=w^{n-1}=(h_1\cdots h_n)^{n-1}$.\hfil}
\endinsert

\beginproclaim
Theorem~\Tag{thm-fundamental-elt}.
Conjugation by $\D$ fixes each of $h_3,\ldots,h_n$, and either
swaps $h_1$ with $h_2$ (if $n$ is odd) or fixes each of them (if
$n$ is even). Furthermore, $\D^2$ (resp. $\D$) is central in
$\braid_n(k)$ if $n$ is odd (resp. even).
\endproclaim

\beginproof{Proof:}
It follows by inspection of fig.~\tag{fig-picture-of-Delta} that $\D$
commutes with $h_i$ for all $i>2$. Next we observe that $\D
h_1\D^{-1}$ is given in fig.~\tag{fig-D-h1-Dinverse}.
If $n$ is odd then this braid is equivalent to that of fig.~\tag{fig-braid-conjugation-odd},
where we have suppressed all but the first two braid strands.
Similarly, if $n$ is even then fig.~\tag{fig-D-h1-Dinverse} reduces to
fig.~\tag{fig-braid-conjugation-even}.
Therefore $\D h_1\D^{-1}=h_1$ if $n$ is even and $\D h_1\D^{-1}=h_2$
if $n$ is odd. A similar calculation shows that $\D h_2\D^{-1}=h_2$ is
$n$ is even and $\D h_2\D^{-1}=h_1$ if $n$ is odd. Since $\braid_n(k)$ is
generated by  $\tau$ and the $h_i$, to finish the proof it suffices to
check that $\D$ and $\tau$ commute. This follows by inspection of
fig.~\tag{fig-picture-of-Delta}. Indeed $w$ and $\tau$ also commute, by
inspection of \eqtag{eq-defn-of-w}.
\pageinsert
\vfil
\hfil
\beginpicture
\put {\epsfbox[22 120  273 432]{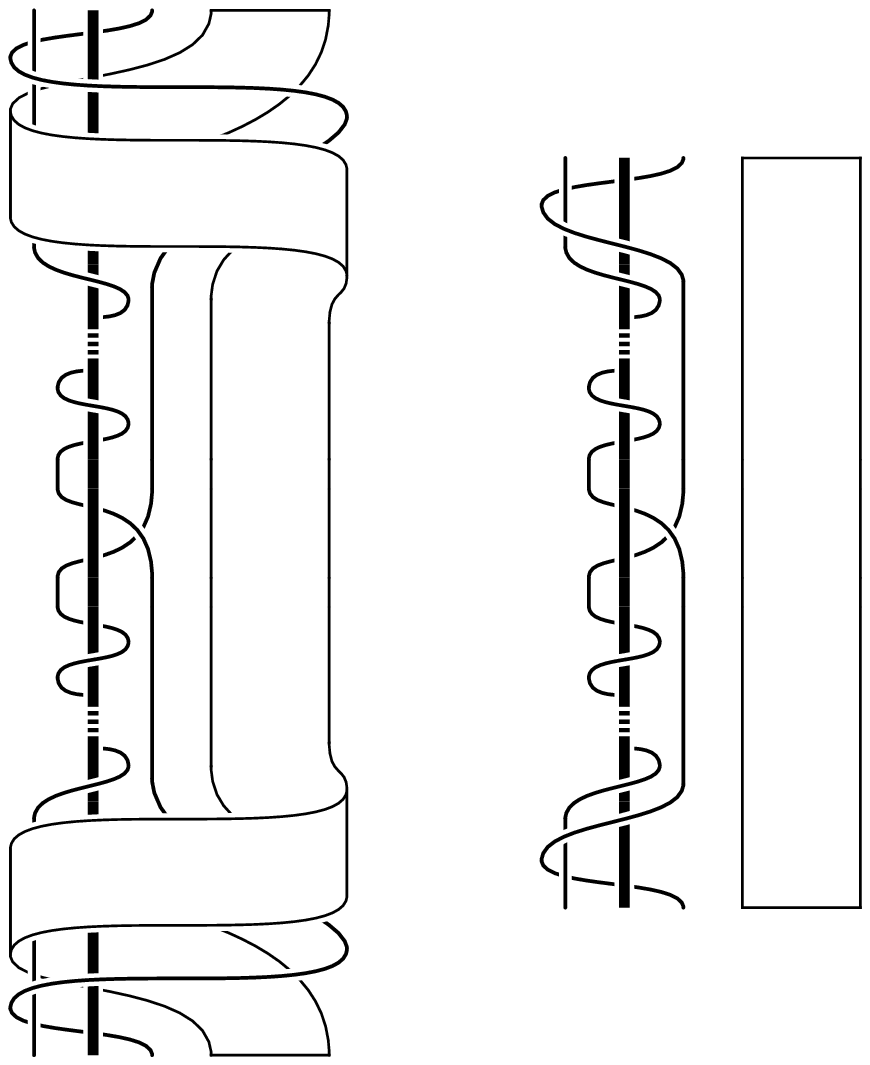}} [bl] at 0 0
\setcoordinatesystem units <1bp,1bp> point at 22 120
\put {$=$} [c] <0pt,0pt> at 158.1 276.25
\endpicture
\hfil
\par
\line{\small\hfil 
Figure~\Tag{fig-D-h1-Dinverse}. A picture of $\D h_1 \D^{-1}$.\hfil}
\vfil
\hfil
\beginpicture
\put {\epsfbox[22 332  311 425]{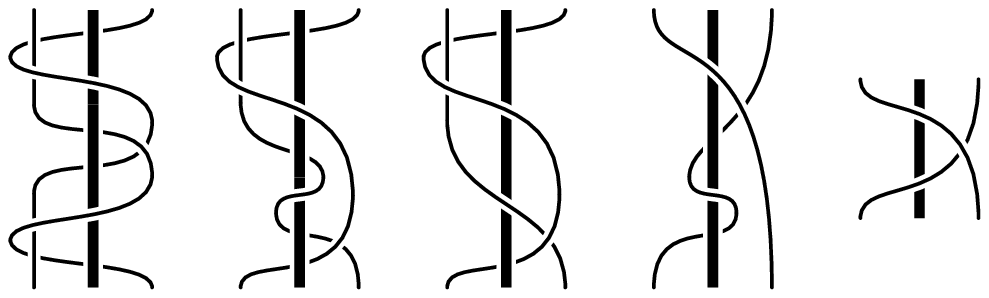}} [bl] at 0 0
\setcoordinatesystem units <1bp,1bp> point at 22 332
\put {$=$} [c] <0pt,0pt> at 80.75 379.95
\put {$=$} [c] <0pt,0pt> at 140.25 379.95
\put {$=$} [c] <0pt,0pt> at 199.75 379.95
\put {$=$} [c] <0pt,0pt> at 259.25 379.95
\endpicture
\hfil
\par
\line{\small\hfil
Figure~\Tag{fig-braid-conjugation-odd}. 
$\D h_1\D^{-1}=h_2$ if $n$ is odd.\hfil}
\vfil
\endinsert
\topinsert
\hfil
\beginpicture
\put {\epsfbox[22 334  209 464]{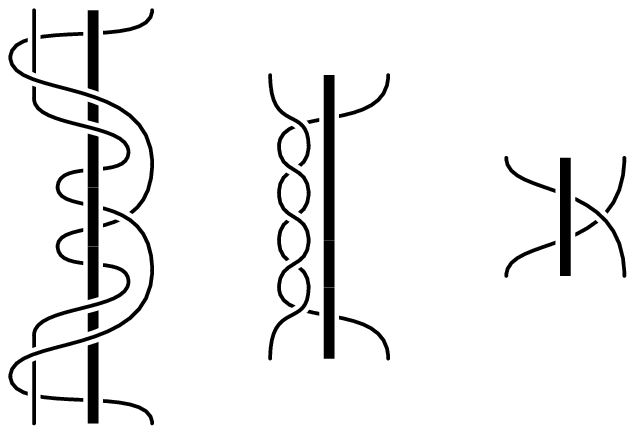}} [bl] at 0 0
\setcoordinatesystem units <1bp,1bp> point at 22 334
\put {$=$} [c] <0pt,0pt> at 85.0 399.5
\put {$=$} [c] <0pt,0pt> at 153.0 399.5
\endpicture
\hfil
\par
\line{\small\hfil
Figure~\Tag{fig-braid-conjugation-even}. 
$\D h_1\D^{-1}=h_1$ if $n$ is even.\hfil}
\endinsert
\endproof

This leads to another description of $\braid_n(k)$ if $n$ is odd:

\beginproclaim
Theorem~\Tag{thm-Bn(k)-for-odd-n}.
If $n$ is odd then a presentation for $\braid_n(k)$ may be obtained
from that of $\artindn$ by adjoining a new central element $z$,
subject to the relation that $z^2=\D^2$.
\endproclaim

\beginproof{Proof:}
Take $z=\tau\D$. By theorems~\tag{thm-artindn} and
\tag{thm-fundamental-elt}, the conjugation maps of $\D$ and
$\tau=\tau^{-1}$ coincide, so $z$ is central in $\braid_n(k)$. Since
$\D$ lies in $\artindn$ but $\tau$ does not,
$z\notin\artindn$. Since $\artindn$ has index two in $\braid_n(k)$,
the description of $\braid_n(k)$ is completed by computing the square
of $z$. In light of $\D\tau=\tau\D$, we have
$z^2=\tau\D\tau\D=\D^2$.
\endproof

We now identify our fundamental element with those of
Charney~\cite{charney:artin-groups-finite-type-biautomatic} and
of Brieskorn and
Saito~\cite{brieskorn:artin-gruppen-und-coxeter-gruppen}. It is easy
to find the lift to $V_0$ of the path in $(k^n-\D_n)/S_n$ represented
by the braid of fig.~\tag{fig-picture-of-Delta}. All we need is that
the lift begins at $p_0=(\e i,1,2,\ldots,n-1)$ and ends at
$p_0'=((-1)^{n-1}\e i,-1,-2,.\ldots,-(n-1))$. One can even read this
fact directly from fig.~\tag{fig-picture-of-Delta}: following the
strands around gives the induced permutation of coordinates, with a
sign-flip being present if the strand passes behind the `orbifold
strand' an odd number of times.  Now, the real part of $p_0'$ lies in
the Weyl chamber opposite to $C$, because it has negative inner
product with each simple root. Therefore the standard map from
$\artindn$ to $W(D_n)$ carries $\D$ to the element $\a$ of $W$ that
exchanges $C$ with $-C$. It is known that $\a$ is the unique longest
element of $W$ with respect to the standard generating set, with
length $n(n-1)$. Since $\D$ is defined as a word of length $n(n-1)$ in
the $h_i$, it is a minimal element of $\artindn$ mapping to $\a$. The
description of the fundamental element in the proof of lemma~2.3 of
\cite{charney:artin-groups-finite-type-biautomatic} shows that our
$\D$ coincides with Charney's.  (Our $\a$ is her $g_0$.)  Brieskorn
and Saito
\ecite{brieskorn:artin-gruppen-und-coxeter-gruppen}{section~5.8}
describe their fundamental element as a specific word of
length $n(n-1)$ in the standard Artin generators that maps to
$\a$. This implies that their description coincides with Charney's and
hence with ours. It is interesting to note that their explicit
expression for $\D$ differs from ours.

\section{\Tag{sec-Dn}}{The affine Artin group $\artin(\tilde{D}_n)$}

In this section we will prove theorem~\tag{thm-artinDn}. The ideas are
similar to those of the previous section but the calculations are
slightly more involved. The affine Weyl group $W=W(\tilde{D}_n)$ is a discrete
cocompact group of isometries of $n$-dimensional real Euclidean space
$V^\R=\R^n$. It is generated by $n+1$ reflections, which we number as
in \eqtag{eq-labelled-Dn}. The first $n$ reflections are across the
hyperplanes orthogonal to the simple roots
\eqtag{eq-simple-roots-for-dn}, and the last is the reflection
across the hyperplane whose elements have inner product $-1$ with the
root $(0,\ldots,0,-1,-1)$. We take $C$ to be the open region in $V^\R$
bounded by these hyperplanes, namely the set of points having inner
product $>-1$ with the last root and positive inner product with each
of the other simple roots. We write $\Sigma$ for the set of
reflections in $W$, and for each $s\in\Sigma$ we write $H_s^\R$ for
the mirror of $s$. It turns out that the $W$-translates of $C$
coincide with the components of
$V_0=V^\R\setminus\cup_{s\in\Sigma}H_s^\R$.

Now we take $V$ to be the complexification of $V^\R$ and for
each $s\in\Sigma$ we take $H_s$ to be the complexification of
$H_s^\R$. By this we mean the unique complex hyperplane in $V$
containing $H_s^\R$. (This definition is slightly different from
the corresponding one in section~\tag{sec-dn}, because here $H_s^\R$ and
$H_s$ might not contain the origin.) One can show that $W$ acts
freely on $V_0$, and it turns out that
$\pi_1(V_0/W)\isomorphism\artinDn$. This is the content of the
following theorem, which is essentially due to \nguyen\ \cite{nguyen:regular-orbits-affine-weyl-groups}.

\beginproclaim
Theorem~\Tag{thm-hyperplane-complement-for-Dn}.
Let $p_0\in V_0$ have real part lying in $C$. For each
$s=1,\ldots,n+1$ let $p_s$ be the image of $p_0$ under the $i$th
generating reflection of $W$ (i.e., across the $i$th wall of $C$), and
let $L_s$ be the complex line in $V$ containing $p_0$ and $p_s$. Let
$G_s$ be the path $[0,1]\to L_s$ from $p_0$ to $p_s$ which 
misses the mirrors and
is similar to the path $g_s$ of section~\tag{sec-dn}. That is, 
on
$[0,{1\over3}]$ and $[{2\over3},1]$ it is linear in the real line
through $p_0$ and $p_s$ and on $[{1\over3},{2\over3}]$ it is a positively
oriented semicircle with small radius and center at $(p_0+p_s)/2$. 
Then $\pi_1(V_0/W)$ is
isomorphic to $\artinDn$ and the loops in $V_0/W$ represented by
$G_1,\ldots,G_{n+1}$ may be taken as the standard Artin generators.
\endproclaim

\beginproof{Proof:}
This is mostly proven in
\ecite{nguyen:regular-orbits-affine-weyl-groups}{section~5}, so we
restrict ourselves to a sketch. We may decompose $V_0$ as a union of
disjoint convex sets in the following way. First we observe that the
various intersections of the mirrors of $W$ define a natural
simplicial structure on $V_0$, in which the translates of the the
closure of $C$ are the top-dimensional simplices. Suppose that $E$ is
a simplex of this decomposition, minus its proper faces. Writing $\pi$
for the natural projection $V_0\to V^\R$, the set $\pi^{-1}(E)$ falls
into various components, each of which turns out to be $E$ times a
copy of the Weyl chamber for the stabilizer of $E$, which is a finite
Coxeter group.  Each of these components is one of the sets in our
decomposition of $V_0$. The full decomposition of $V_0$ is obtained
by varying $E$. If $E$ has codimension $k$ in
$V_0$ then each of these subsets in $V_0$ lying over $E$ has
codimension $k$ in $V_0$. The decomposition is preserved by $W$, and
$W$ permutes the `cells' freely. Furthermore, two `cells' are
$W$-equivalent just if the simplices of $V^\R$ over which they lie are
$W$-equivalent. This implies that $V_0/W$ has a sort of cell structure
with one `cell' of codimension 0, represented by $\pi^{-1}(C)$, and one
`cell' of codimension $1$ for each wall of $C$.

For each $s=1,\ldots,n+1$, $G_s$ is a path in $V_0$ beginning at
$p_0$, passing through a codimension one `cell' of $V_0$ lying
over the $s$th wall of $C$, and continuing to the image $p_s$ of
$p_0$ in the neighboring Weyl chamber. These paths represent
loops in $V_0/W$, and the codimension-one `cells' through which
they pass exhaust all the `cells' of that dimension in
$V_0/W$. Therefore they generate $\pi_1(V_0/W)$. One can find
defining relations in a similar manner by considering the
codimension-two stratum of $V^\R$. \nguyen\  does this in detail,
and one can also use Brieskorn's argument \cite{brieskorn:fundamentalgruppe}.

We have sketched the proof here because our $V_0$ is different from
the space \nguyen\ claims to treat, namely
$V\setminus\cup_{s\in\Sigma}(H_s^\R+\sqrt{-1}H_s^\R)$, because not all
of the $H_s^\R$ contain 0. Perhaps this is a notational error. In any
case, his ideas suffice to prove the theorem. The necessary
prerequisites for the manipulations of the `cells' are developed in
his paper. One alternative to developing this machinery is to study
the dual complex, as is done in
\cite{charney:finite-k-pi-1's-for-artin-groups} or
\cite{salvetti:homotopy-type-of-artin-groups}.  We also note that van
der Lek (\cite{van-der-lek:extended-artin-groups},
\cite{van-der-lek:thesis}) has obtained a more general version of this
theorem by using similar techniques  in a slightly
different setting.
\endproof

Another description of $W$, more useful for some purposes, is
the following. We denote by $\Lambda$ the integral span of the
roots of the $D_n$ root system. This is called the $D_n$ root lattice
and consists of all vectors in $\Z^n$ with even coordinate
sum. Then $W$ has structure $\Lambda{\,:\,}W_0$, where $\Lambda$
acts by translations and $W_0$ is the finite Weyl
group $W(D_n)$ studied in the previous section. Since $W_0$ has
structure $2^{n-1}{\,:\,}S_n$ and the $S_n$ normalizes $\Lambda$, we
may indicate the structure of $W$ by
$\Lambda{\,:\,}2^{n-1}{\,:\,}S_n$. The mirrors of $W$ are the
$\Lambda$-translates of those passing through $0$, which in turn
are the mirrors of $W_0$.

This allows us to describe the complex hyperplane arrangement
associated to $W$. A point $x\in\C^n$ lies in the mirror of a
reflection of $W$ if and only if it has integral inner product
with one of the roots of $D_n$. For if $x\cdot r=m\in\Z$ for
a root $r$ then upon choosing another root $r'$ with $r\cdot
r'=1$, we have $(x-mr')\cdot r=0$, so that $x$ lies in a
$\Lambda$-translate of $r^\perp$. On the other hand, if
for some $\lambda\in \Lambda$, $x-\lambda$ lies in $r^\perp$ for
a root $r$, then $r\cdot x= r\cdot \lambda$, which is integral
because all inner products of elements of $\Lambda$ are
integral. Therefore
$$ 
V_0=\{(x_1,\ldots,x_n)\in\C^n
|
\hbox{$x_i\pm x_j\notin\Z$ if $i\neq j$}\}. 
$$
We can now prove theorem~\tag{thm-artinDn}.

\beginproof{Proof of theorem~\tag{thm-artinDn}:}
We recall that $K$ is the orbifold $\C$ with two cone points of
order two. We take the cone points to be at $0$ and $1/2$, and we may
regard the universal orbifold covering map $\C\to K$ to be given by
$x\mapsto\xi(x)=(1-\cos(2\pi x))/4$. On the strip $\re
x\in[0,{1\over2}]$  the map $\xi$ may be visualized as follows. It
fixes each of $0$ and $1/2$, with branching of order 2 there, and each
component of the boundary of the strip is folded upon itself (about
either $0$ or $1/2$) and carried to the real axis.

We have an inclusion of groups from
$W=\Lambda{\,:\,}2^{n-1}{\,:\,} S_n$
into $\Wbar=\Z^n{\,:\,}2^n\mathbin{:}S_n$, where $\Z^n$ acts on $\C^n$ by
translations, $2^{n-1}{\,:\,}S_n$ is $W_0$, and the group $2^n$ acts
by changing coordinates' signs, as in section~\tag{sec-dn}. Because $\Lambda$
has index two in $\Z^n$ and $2^{n-1}$ has index two in $2^n$,
$W$ has index four in $\Wbar$. We will now show that
$(K^n\setminus\D_n)/S_n\isomorphism((V_0/\Z^n)/2^n)/S_n$. The map $V\to
V/\Z^n$ may be described by
$(x_1,\ldots,x_n)\mapsto(y_1,\ldots,y_n)$, where $y_k=\exp(2\pi
ix_k)$. The action of $2^n$ on $V/\Z^n=(\C^\times)^n$ induced by
the action of $\Z^n{\,:\,}2^n$ on $V$ is by replacing some number of
coordinates by their reciprocals. 
The map 
$$
y\mapsto\eta(y)={1\over4}\(1-{1+y^2\over 2y}\)
$$
identifies $y,y'\in\C^\times$ just if they are equal or
reciprocal. Furthermore it carries the branch points (namely $\pm1$)
of the map to $0$ and $1/2$. Therefore we may regard it as an orbifold
covering map (of degree 2) from $\C^\times$ to $K$.
Therefore the map
$V/\Z^n\to(V/\Z^n)/2^n=K^n$ is given by
$(y_1,\ldots,y_n)\mapsto(z_1,\ldots,z_n)$, where
$z_k=\eta(y_k)=\xi(x_k)$. 
It is easy to compute the image of $V_0$ in these models
of $V/\Z^n$ and $V/(\Z^n{\,:\,}2^n)$. We have
$$
V_0/\Z^n 
=\{\,(y_1,\ldots,y_n)\in(\C^\times)^n 
\>|\>
\hbox{$y_i\neq y_j^{\pm1}$ if $i\neq j$}\,\}
\hbox{\quad and}
$$
$$
V_0/(\Z^n{\,:\,}2^n)
=\{\,(z_1,\ldots,z_n)\in K^n 
\>|\>
\hbox{$z_i\neq z_j$ if $i\neq j$}\,\}.
$$
Finally, our map $V\mapsto K^n$ identifies the actions of $S_n$
on its domain and range. Therefore $V_0/\Wbar$ is the $n$-strand
braid space of $K$. This proves that $\artinDn$ has index 4 in
$\braid_n(K)$. 

Next we compute the map on fundamental groups induced by this orbifold
cover.
We take our basepoint in $V_0$ to be
$$ 
p_0={1\over2n-2}(i,1,2,\ldots, n-2,n-1+i). 
$$
It is easy to check that the real part of $p_0$ lies in $C$. We
may take the standard generators for $\artinDn$ to be the paths
$G_s$ of theorem~\tag{thm-hyperplane-complement-for-Dn}; most of the rest of the proof is a
computation of their images in $K^n$. For $s=1,\ldots,n$ the
analysis is similar to that of the proof of theorem~\tag{thm-artindn},
because the reflections are the same. For $s=n+1$, 
the relevant
reflection is across the hyperplane of $\C^n$ whose elements
have inner product $-1$ with $r=(0,\ldots,0,-1,-1)$. 
The image of $p_0$ under this map is
$$ 
p_{n+1}={1\over 2n-2}(i,1,\ldots,(n-1)-i,n), 
$$
which may be verified by observing that $(p_0+p_{n+1})/2$ lies on the
mirror and that $p_0-p_{n+1}$ is proportional to $r$. Here are the
paths $G_1,\ldots,G_{n+1}$, described as the $g_s$ were in
section~\tag{sec-dn}, by simultaneously drawing all the coordinate
projections of the path.
\def\widebox#1{%
\setbox0=\hbox{$G_{n-1}$}%
\llap{\hbox to\wd0{\hfil #1\hfil}%
\qquad}}
$$ 
\beginpicture
\put {\epsfbox[40 445  200 484]{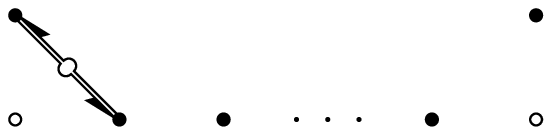}} [bl] at 0 0
\setcoordinatesystem units <1bp,1bp> point at 40 445
\put {\strut$0$} [t] <0pt,-1pt> at 45.0 448.0
\put {$\displaystyle\;{1\over2}$} [l] <0pt,-1.6pt> at 197.0 452.0
\put {\widebox{$G_1$}} [r] <0pt,0pt> at 13.0 450.0
\put {\vrule width 0pt depth 0pt height 0pt} [c] <0pt,0pt> at 240.0 450.0
\endpicture
$$
\medskip
$$ 
\beginpicture
\put {\epsfbox[12 417  200 484]{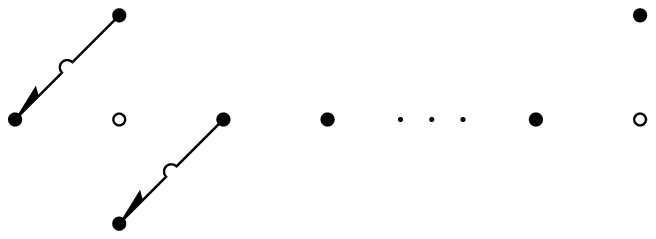}} [bl] at 0 0
\setcoordinatesystem units <1bp,1bp> point at 12 417
\put {\widebox{$G_2$}} [r] <0pt,0pt> at 13.0 465.0
\put {\vrule width 0pt depth 0pt height 0pt} [c] <0pt,0pt> at 240.0 450.0
\endpicture
$$
\medskip
$$ 
\beginpicture
\put {\epsfbox[40 444  200 483]{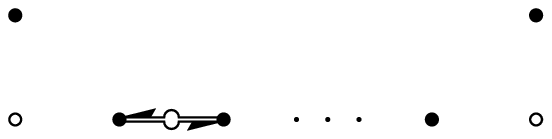}} [bl] at 0 0
\setcoordinatesystem units <1bp,1bp> point at 40 444
\put {\widebox{$G_3$}} [r] <0pt,0pt> at 13.0 465.0
\put {\vrule width 0pt depth 0pt height 0pt} [c] <0pt,0pt> at 240.0 450.0
\endpicture
$$
\medskip
$$ 
\beginpicture
\put {\epsfbox[40 444  200 483]{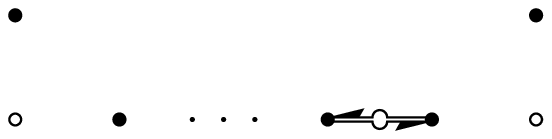}} [bl] at 0 0
\setcoordinatesystem units <1bp,1bp> point at 40 444
\put {\widebox{$G_{n-1}$}} [r] <0pt,0pt> at 13.0 465.0
\put {\vrule width 0pt depth 0pt height 0pt} [c] <0pt,0pt> at 240.0 450.0
\endpicture
$$
\medskip
$$ 
\beginpicture
\put {\epsfbox[40 445  200 484]{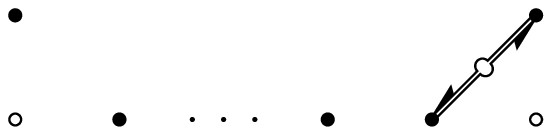}} [bl] at 0 0
\setcoordinatesystem units <1bp,1bp> point at 40 445
\put {\widebox{$G_n$}} [r] <0pt,0pt> at 13.0 465.0
\put {\vrule width 0pt depth 0pt height 0pt} [c] <0pt,0pt> at 240.0 450.0
\endpicture
$$
\medskip
$$ 
\beginpicture
\put {\epsfbox[40 417  229 484]{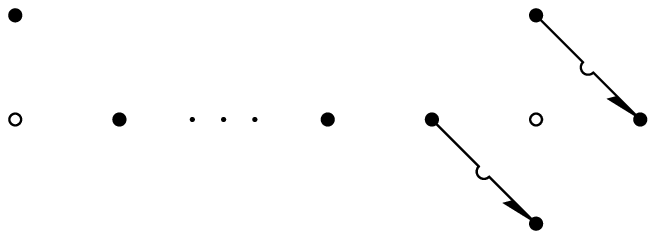}} [bl] at 0 0
\setcoordinatesystem units <1bp,1bp> point at 40 417
\put {\widebox{$G_{n+1}$}} [r] <0pt,0pt> at 13.0 450.0
\put {\vrule width 0pt depth 0pt height 0pt} [c] <0pt,0pt> at 240.0 450.0
\endpicture
$$
Their images in $K^n=V/(\Z^n{\,:\,}2^n)$ are given by applying the map
$\xi$ to the previous pictures. Our curious choice of the
covering map $\xi:\C\to K$ was made to make this step as easy as
possible. 
$$
\beginpicture
\put {\epsfbox[431 2095  618 2113]{paths_Dn_3_A.ps}} [bl] at 0 0
\setcoordinatesystem units <1bp,1bp> point at 431 2095
\put {\strut$0$} [t] <0pt,-1pt> at 450.0 2098.0
\put {\strut$1/2$} [t] <0pt,-1pt> at 600.0 2098.0
\put {\widebox{$G_1$}} [r] <0pt,0pt> at 434.702 2100.0
\endpicture
$$
\smallskip
$$
\beginpicture
\put {\epsfbox[431 2087  618 2105]{paths_Dn_3_B.ps}} [bl] at 0 0
\setcoordinatesystem units <1bp,1bp> point at 431 2087
\put {\widebox{$G_2$}} [r] <0pt,0pt> at 434.702 2100.0
\endpicture
$$
\smallskip
$$
\beginpicture
\put {\epsfbox[432 2093  618 2107]{paths_Dn_3_C.ps}} [bl] at 0 0
\setcoordinatesystem units <1bp,1bp> point at 432 2093
\put {\widebox{$G_3$}} [r] <0pt,0pt> at 434.702 2100.0
\endpicture
$$
\smallskip
$$
\beginpicture
\put {\epsfbox[432 2093  618 2107]{paths_Dn_3_D.ps}} [bl] at 0 0
\setcoordinatesystem units <1bp,1bp> point at 432 2093
\put {\widebox{$G_{n-1}$}} [r] <0pt,0pt> at 434.702 2100.0
\endpicture
$$
\smallskip
$$
\beginpicture
\put {\epsfbox[432 2095  619 2113]{paths_Dn_3_E.ps}} [bl] at 0 0
\setcoordinatesystem units <1bp,1bp> point at 432 2095
\put {\widebox{$G_n$}} [r] <0pt,0pt> at 434.702 2100.0
\endpicture
$$
\smallskip
$$
\beginpicture
\put {\epsfbox[432 2087  619 2105]{paths_Dn_3_F.ps}} [bl] at 0 0
\setcoordinatesystem units <1bp,1bp> point at 432 2087
\put {\widebox{$G_{n+1}$}} [r] <0pt,0pt> at 434.702 2100.0
\endpicture
$$
Drawing the braids associated to these paths and comparing them with
the given braids $H_i$ shows that our orbifold covering
map $V_0/W\to X_n$ carries $G_i$ to $H_i$.

Now we establish the semidirect product structure. The following
two elements of $\braid_n(K)$ generate a group $(\Z/2)^2$.
$$ 
\beginpicture
\put {\epsfbox[29 334  303 379]{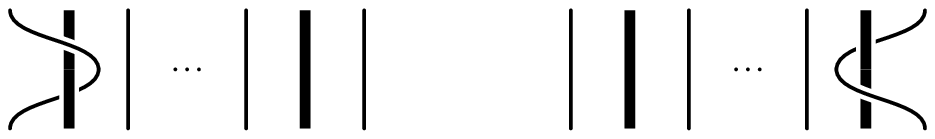}} [bl] at 0 0
\setcoordinatesystem units <1bp,1bp> point at 29 334
\put {$\tau_1\;=$} [r] <0pt,0pt> at 25.5 357.0
\put {$\tau_2\;=$} [r] <0pt,0pt> at 187.0 357.0
\endpicture
$$
The elements of this group are $1$, $\tau_1$, $\tau_2$ and
$\tau_1\tau_2$. To see that of these only the trivial element
lies in $\artinDn$ one can obtain the lifts to $V_0$ of the paths in
$K^n\setminus\D_n$ represented by these braids and determine their
final endpoints. These turn out to be 
${1\over2n-2}(\pm i,1,\ldots,n-2,n-1\pm i)$
Since the real parts all lie in the open Weyl chamber $C$
(indeed they coincide), they are inequivalent under $W$ and
hence only one of the paths represents a loop in $V_0/W$.
Since $\artinDn$ has index 4 in $\braid_n(K)$, we have found a
complete set of coset representatives. It is easy to check that the
conjugation maps of $\tau_1$ and $\tau_2$ have the properties
stated. In particular, this proves that $\langle\tau_1,\tau_2\rangle$
normalizes $\artinDn$, so that the semidirect product decomposition
exists as claimed.
\endproof

\section{\Tag{sec-other-artin-groups}}{The remaining groups}

In this section we sketch the arguments for the rest of the
diagrams listed in table~\tag{tab-results}, namely the
spherical diagrams $A_{n-1}$ and $B_n=C_n$ and the affine diagrams
$\tilde{A}_{n-1}$, $\tilde{B}_n$ and $\tilde{C}_n$. The $A_{n-1}$ diagram appears in
\eqtag{eq-an-1-diagram}. The remaining diagrams are given below.
$$
\beginpicture
\put {\epsfbox[27 697  221 703]{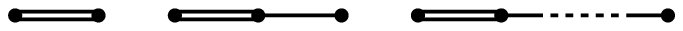}} [bl] at 0 0
\setcoordinatesystem units <1bp,1bp> point at 27 697
\put {\strut$B_2$} [t] <0pt,0pt> at 42.0 698.0
\put {\strut$B_3$} [t] <0pt,0pt> at 100.0 698.0
\put {\strut$B_n$} [t] <0pt,0pt> at 182.0 698.0
\endpicture
$$
$$
\beginpicture
\put {\epsfbox[11 672  237 728]{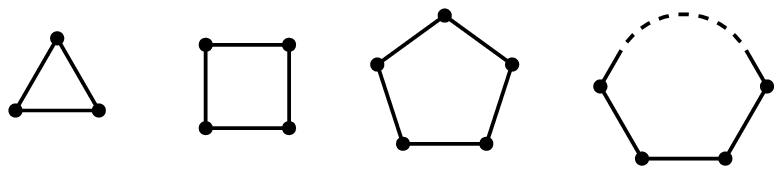}} [bl] at 0 0
\setcoordinatesystem units <1bp,1bp> point at 11 672
\put {\strut$\tilde{A}_2$} [t] <0pt,0pt> at 30.0 677.215
\put {\strut$\tilde{A}_3$} [t] <0pt,0pt> at 84.7844 677.215
\put {\strut$\tilde{A}_4$} [t] <0pt,0pt> at 141.61 677.215
\put {\strut$\tilde{A}_n$} [t] <0pt,0pt> at 210.435 677.215
\endpicture
$$
$$
\beginpicture
\put {\epsfbox[3 676  233 724]{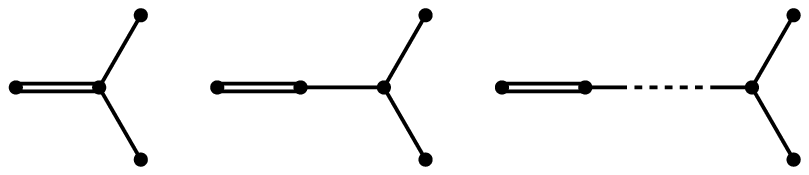}} [bl] at 0 0
\setcoordinatesystem units <1bp,1bp> point at 3 676
\put {\strut$\tilde{B}_3$} [t] <0pt,0pt> at 24.0 677.215
\put {\strut$\tilde{B}_4$} [t] <0pt,0pt> at 94.0 677.215
\put {\strut$\tilde{B}_n$} [t] <0pt,0pt> at 188.0 677.215
\endpicture
$$
$$
\beginpicture
\put {\epsfbox[27 697  293 703]{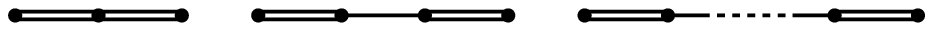}} [bl] at 0 0
\setcoordinatesystem units <1bp,1bp> point at 27 697
\put {\strut$\tilde{C}_2$} [t] <0pt,-2pt> at 54.0 698.0
\put {\strut$\tilde{C}_3$} [t] <0pt,-2pt> at 136.0 698.0
\put {\strut$\tilde{C}_n$} [t] <0pt,-2pt> at 242.0 698.0
\endpicture
$$
\medskip\noindent 
Recall that a diagram $x_n$ has $n$ nodes and a
diagram $X_n$ has $n+1$. The Artin groups are given by taking a
generator for each node and declaring that two such generators $x$ and
$y$ commute (resp. braid, or satisfy $xyxy=yxyx$) if the corresponding
nodes are unjoined (resp. singly joined, or doubly joined). For each of the
diagrams, the pure Artin group is the fundamental group of the
complex hyperplane complement $V_0$ associated to the corresponding
Weyl group. The Artin group itself is the fundamental group of the
quotient of $V_0$ by the Weyl group. The standard Artin generators may
be represented by paths in $V_0$ in the same manner as in
sections~\tag{sec-dn} and \tag{sec-Dn}.
Each of these Artin groups may be realized as a normal subgroup of the
$n$-strand braid group of a suitable 2-orbifold, as indicated in
table~\tag{tab-results}. The results for the cases $\tilde{B}_n$ and $\tilde{C}_n$
are new. In each case the argument has much the same form as that of
section~\tag{sec-dn} or \tag{sec-Dn}. We now treat each case briefly.

\smallskip
{\bf The diagram $A_{n-1}$:} 
This approach to the classical braid group appears in
\cite{fox-neuwirth:braid-groups}. Here,
$$
V_0=\{\,(x_1,\ldots,x_n)\in\C^n
\>|\>
\hbox{$\sum x_i=0$, $x_j\neq x_k$ if $j\neq k$}\,\}
$$
and the Weyl group $W=W(A_{n-1})$ acts by permuting the coordinates. The
space $V_0$ is an $S_n$-equivariant deformation retract of 
$$
\{\,(x_1,\ldots,x_n)\in\C^n
\>|\>
\hbox{$x_j\neq x_k$ if $j\neq k$}\,\}
$$,
showing that $V_0/W$ is homotopy-equivalent to the braid space of $\C$.
It is easy to see that the standard Artin generators correspond to the
standard braid generators. A picture of the fundamental element of
$\artin(A_{n-1})$ appears in \ecite{dbae:wd-proc-gps}{fig.~9.2}.

\smallskip
{\bf The diagram $B_n=C_n$:} The results here are implicit in
\cite{brieskorn:groupes-de-tresses}. Here $W=W(B_n)$ is the group
$2^n{\,:\,}S_n$ used in section~\tag{sec-dn}, so that $W$ contains
$W(D_n)$ as a subgroup of index 2. Each reflection negates a
coordinate, or exchanges two coordinates, or exchanges two coordinates
and negates both. The hyperplane complement is
$
V_0=\{\,(x_1,\ldots,x_n)\in(\C^\times)^n
\>|\>
\hbox{$x_j\neq \pm x_k$ if $j\neq k$}\,\}
$.
The map $V_0\to V_0/2^n=(\C^\times)^n$ given by squaring each
coordinate identifies $V_0/2^n$ with the pure braid space of
$\C^\times$ in an $S_n$-equivariant manner. This shows
$\artin(B_n)=\braid_n(\C^\times)$. In terms of braids, the standard Artin
generators (from left to right) are
$$
\beginpicture
\put {\epsfbox[20 334  107 379]{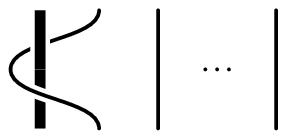}} [bl] at 0 0
\setcoordinatesystem units <1bp,1bp> point at 20 334
\put {\strut$\infty$} [t] <0pt,0pt> at 34.0 340.0
\endpicture
$$
$$
\beginpicture
\put {\epsfbox[49 334  294 379]{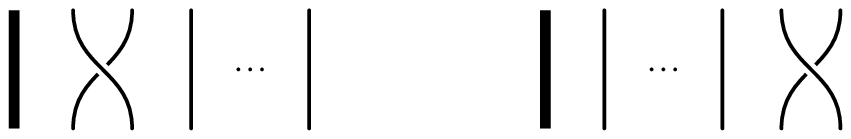}} [bl] at 0 0
\setcoordinatesystem units <1bp,1bp> point at 49 334
\put {\strut$\infty$} [t] <0pt,0pt> at 51.0 340.0
\put {\strut$\infty$} [t] <0pt,0pt> at 204.0 340.0
\endpicture
$$
The `$\infty$' below the heavy line indicates that the line
represents the puncture. It is a pleasing exercise to verify
directly the `four-braid' relation satisfied by the first two
generators. The fundamental element of $\artin(B_n)$ is given by the
$n$th power of the product (from left to right) of the standard
generators. It is easy to check that this element is central.

We remark that tom~Dieck (\cite{tom-dieck:rooted-cylinder-ribbons},
\cite{tom-dieck:symmetrische-brucken}) has studied the problem of
generalizing to this case the machinery of tangle and skein theory
associated to the classical braid group. We wonder whether our results for
the other classical Artin groups will lead to a whole series of skein theories
in various orbifolds.
tom~Dieck's pictures of elements of
$\artin(B_n)$ are braids in $\C$ that are symmetric with respect to
the rotation of order two about the origin. These contain the same
information as our pictures in $\C^\times$ but take more time to draw.
In \ecite{tom-dieck:symmetrische-brucken}{p.~84} he shows that this embeds
$\braid_n(\C^\times)$ into $\braid_{2n}(\C)$. In
\ecite{tom-dieck:rooted-cylinder-ribbons}{p.~37} he also describes
$\artin(B_n)$ as the group of ordinary braids equipped with
framings. This description is not valid: the latter group is the
semidirect product of $\Z^n$ by $\braid_n(\C)$ and has the wrong
cohomological dimension. However, he also considers the natural map
$\braid_n(\C^\times)\to\braid_n(\C)$ given by filling in the puncture, and
observes that the kernel is a certain free group of rank
$n$. The group of framed $n$-strand braids in $\C$ may be
obtained by quotienting by the commutator subgroup of this free group.

\smallskip
{\bf The diagram $\tilde{A}_n$:} Most of the results here are stated in
\cite{tom-dieck:rooted-cylinder-ribbons}.  The Weyl group $W$ is
generated by the reflections in the roots $(1,-1,0,\ldots,0)$,
$\ldots$, $(0,\ldots,0,1,-1)$ in $\C^n$, together with the reflection
across the hyperplane of elements having inner product $-1$ with
$(-1,0,\ldots,0,1)$. Alternatively, $W=\Lambda:S_n$ where $\Lambda$ is
the integral span of the roots, which is the set of elements of $\Z^n$
with vanishing coordinate sum. We have
$$
V_0=\{\,(x_1,\ldots,x_n)\in\C^n
\>|\>
\hbox{$\sum x_j=0$, $x_j-x_k\notin\Z$ if $j\neq k$}
\,\}
\quad.
$$
The easiest way to
compute the quotient of $V_0$ by $\Lambda$ is to take the
quotient of $\C^n$ by $\Z^n$ by applying the map $x\mapsto\exp(2\pi
ix)$ for each coordinate. Then the $n$-strand pure braid space for
$\C^\times$ is
$
\{\,(y_1,\ldots,y_n)\in(\C^\times)^n
\>|\>
\hbox{$x_j\neq x_k$ if $j\neq k$}
\,\}
$,
the map from this space to $\C^\times$ given by
$(y_1,\ldots,y_n)\mapsto y_1\cdots y_n$ is a locally trivial
fibration, and $V_0/\Lambda$ is the fiber over $1$. The long exact
homotopy sequence shows that $\artin(\tilde{A}_{n-1})$ is normal in
$\braid_n(\C^\times)$ with quotient $\Z$. A good choice for a set of
basepoints for $n$-strand braids in $\C^\times$ is the set of $n$th
roots of unity. Then the standard Artin generators correspond to the
standard braid generators; each of these exchanges two adjacent
basepoints in the simplest possible way.

\smallskip
{\bf The diagram $\tilde{B}_n$:} 
The first $n$ generators for $W=W(\tilde{B}_n)$ (all except for the top right
node) may be taken to be the reflections in the roots
$(-1,0,\ldots,0)$, $(1,-1,0,\ldots,0)$, $\ldots$, $(0,\ldots,0,1,-1)$
in $\C^n$, and the last may be taken to be across the hyperplane of
points having inner product $-1$ with $(0,\ldots,0,1,1)$. From this
one can obtain a description of $W$ as $\Lambda{\,:\,}W(B_n)$
where $\Lambda$ is
the $D_n$ root lattice considered in section~\tag{sec-Dn}. 
We observe that $W(\tilde{D}_n)$ has index 2 in $W$. Also,
$$
V_0=\{\,
(x_1,\ldots,x_n)\in\C^n
\>|\>
\hbox{$x_j\pm x_k\notin\Z$ for $j\neq k$ and $x_j\notin\Z$ for all
$j$}
\,\}
\quad.
$$
Since $\Lambda$ has index 2 in $\Z^n$, $W$ has index 2 in
$\Z^n{\,:\,}2^n{\,:\,}S_n$. Closely following the analysis of
section~\tag{sec-Dn}, with $p_0$ replaced by
${1\over2n}(1,2,\ldots,n-1,n+i)$, shows that $\artin(\tilde{B}_n)$ has index
2 in $\braid_n(L)$ where $L$ is $\C^\times$ with an orbifold point of
order two at $1/2$. One can show that the standard Artin generators
correspond in the order given to the $n$-strand braids
$$
\beginpicture
\put {\epsfbox[20 334  154 379]{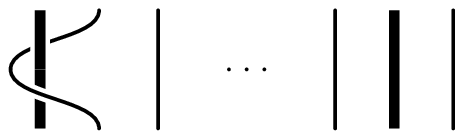}} [bl] at 0 0
\setcoordinatesystem units <1bp,1bp> point at 20 334
\put {\strut$\infty$} [t] <0pt,0pt> at 34.0 340.0
\put {\strut$2$} [t] <0pt,0pt> at 136.0 340.0
\endpicture
$$
$$
\beginpicture
\put {\epsfbox[32 334  341 379]{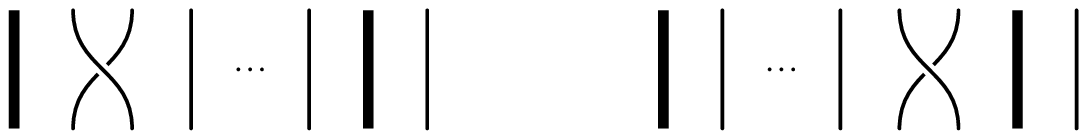}} [bl] at 0 0
\setcoordinatesystem units <1bp,1bp> point at 32 334
\put {\strut$\infty$} [t] <0pt,0pt> at 34.0 340.0
\put {\strut$2$} [t] <0pt,0pt> at 136.0 340.0
\put {\vrule width 0pt height 0pt depth 0pt} [c] <0pt,0pt> at 25.5 357.0
\put {\strut$\infty$} [t] <0pt,0pt> at 221.0 340.0
\put {\strut$2$} [t] <0pt,0pt> at 323.0 340.0
\put {\vrule width 0pt height 0pt depth 0pt} [c] <0pt,0pt> at 348.5 357.0
\endpicture
$$
$$
\beginpicture
\put {\epsfbox[32 334  345 379]{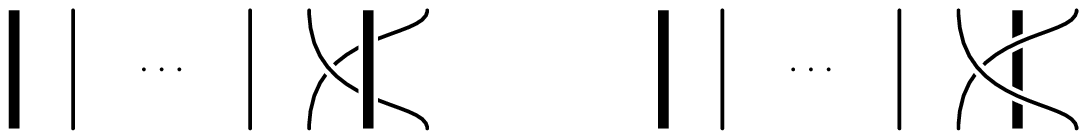}} [bl] at 0 0
\setcoordinatesystem units <1bp,1bp> point at 32 334
\put {\strut$\infty$} [t] <0pt,0pt> at 34.0 340.0
\put {\strut$2$} [t] <0pt,0pt> at 136.0 340.0
\put {\vrule width 0pt height 0pt depth 0pt} [c] <0pt,0pt> at 25.5 357.0
\put {\strut$\infty$} [t] <0pt,0pt> at 221.0 340.0
\put {\strut$2$} [t] <0pt,0pt> at 323.0 340.0
\put {\vrule width 0pt height 0pt depth 0pt} [c] <0pt,0pt> at 348.5 357.0
\endpicture
$$
A subgroup $\Z/2$ complementary to $\artin(\tilde{B}_n)$ may be obtained by
taking the obvious analogue of the element $\tau$ of eq.~\eqtag{eq-def-of-tau}.

\smallskip
{\bf The diagram $\tilde{C}_n$:} 
The first $n$ generators of $W=W(\tilde{C}_n)$ act on $\C^n$ by the
reflections in the roots $(1,0,\ldots,1)$, $(-1,1,0,\ldots,0)$,\dots,
$(0,\ldots,0,-1,1)$, and the last acts by the the reflection
across the hyperplane of vectors having inner product $-1/2$ with
$(0,\ldots,0,-1)$. Alternatively, $W$ is the group $\Z^n{\,:\,}2^n{\,:\,}S_n$,
so that $W(\tilde{D}_n)$ has index 4 and $W(\tilde{B}_n)$ has index 2. We have
$$
V_0=\{\,
(x_1,\ldots,x_n)\in\C^n
\>|\>
\hbox{$x_j\pm x_k\notin\Z$ for $j\neq k$ and $x_j\notin{1\over2}\Z$ for all
$j$}
\,\}
\quad.
$$
The map $V_0\to V_0/(\Z^n{\,:\,}2^n)$ may be described as in
section~\tag{sec-Dn}, and this identifies $V_0/(\Z^n{\,:\,}2^n)$ with the pure
braid space of $\C-\{0,1/2\}$. 
A suitable basepoint for the computation of the braids associated to
the Artin generators is $p_0={1\over2n+1}(1,\ldots,n)$. The
resulting braids are pictured below, in order.
$$
\beginpicture
\put {\epsfbox[20 334  308 379]{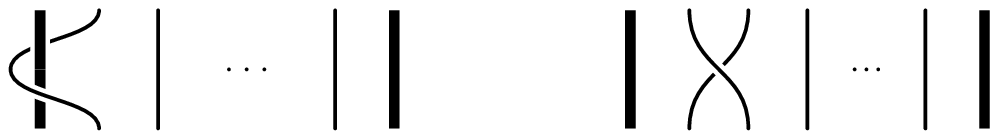}} [bl] at 0 0
\setcoordinatesystem units <1bp,1bp> point at 20 334
\put {\strut$\infty$} [t] <0pt,0pt> at 34.0 340.0
\put {\strut$\infty$} [t] <0pt,0pt> at 136.0 340.0
\put {\vrule width 0pt height 0pt depth 0pt} [c] <0pt,0pt> at 17.0 357.0
\put {\strut$\infty$} [t] <0pt,0pt> at 204.0 340.0
\put {\strut$\infty$} [t] <0pt,0pt> at 306.0 340.0
\put {\vrule width 0pt height 0pt depth 0pt} [c] <0pt,0pt> at 323.0 357.0
\endpicture
$$
$$
\beginpicture
\put {\epsfbox[32 334  320 379]{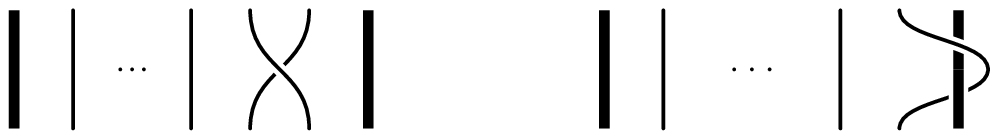}} [bl] at 0 0
\setcoordinatesystem units <1bp,1bp> point at 32 334
\put {\strut$\infty$} [t] <0pt,0pt> at 34.0 340.0
\put {\strut$\infty$} [t] <0pt,0pt> at 136.0 340.0
\put {\vrule width 0pt height 0pt depth 0pt} [c] <0pt,0pt> at 17.0 357.0
\put {\strut$\infty$} [t] <0pt,0pt> at 204.0 340.0
\put {\strut$\infty$} [t] <0pt,0pt> at 306.0 340.0
\put {\vrule width 0pt height 0pt depth 0pt} [c] <0pt,0pt> at 323.0 357.0
\endpicture
$$

\medskip
We close by observing a few curiosities that our orbifold approach
makes visible. We recall that $k$ is the orbifold
$\C/(z\mapsto-z)$. First and most curious is the fact that
$\braid_4(\C)$ has index 2 in $\braid_3(k)$ because of the coincidence
of diagrams $A_3=D_3$. There is a similar coincidence $\tilde{A}_3=\tilde{D}_3$, so
that this affine Artin group has index $2$ in $\braid_3(K)$ and also
index $\infty$ in $\braid_4(\C)$. I do not know of any elementary
explanation for these coincidences.

The other curiosity arises from the natural map
$\braid_n(\C^\times)\to \braid_n(k)$ given by filling in the puncture with
a cone point of order $2$. One can check that
it carries $\artin(\tilde{A}_{n-1})\sset \braid_n(\C^\times)$ onto
$\artin(D_n)\sset \braid_n(k)$. It is not at all clear from the Coxeter
diagrams that there is any surjection $\artin(\tilde{A}_{n-1})\to\artin(D_n)$
at all.
It is natural to wonder whether adjoining
the relation that the leftmost generator of
$\artin(B_n)=\braid_n(\C^\times)$ have trivial square suffices
to reduce $\braid_n(\C^\times)$ to $\braid_n(k)$. 

There is also a natural map $\braid_n(\C^\times)\to\braid_n(\C)$ given
by filling in the puncture with an ordinary point. This map induces a retraction
$\artin(\tilde{A}_{n-1})\to\artin(A_{n-1})$, where we regard $\artin(A_{n-1})$
as embedded in $\artin(\tilde{A}_{n-1})$ by any embedding of the Coxeter
diagrams. This retraction has the pleasant property that it descends
to the natural retraction $W(\tilde{A}_{n-1})\to W(A_{n-1})$ given by
considering the natural action of $W(\tilde{A}_{n-1})$ on the sphere at
infinity of Euclidean space. We remark that the corresponding natural
map $\artin(\tilde{D}_n)\to\artin(D_n)$ obtained by replacing a cone point of
$K$ by an ordinary point does {\it not} have this property. I do not
know if there is any retraction $\artin(\tilde{D}_n)\to\artin(D_n)$ 
which does descend to the natural retraction of Weyl groups.

\reflist

\bibitem{brieskorn:fundamentalgruppe}
E.~Brieskorn.
 Die {F}undamentalgruppe des {R}aumes der regul\"aren {O}rbits einer
  endlichen komplexen {S}piegelungsgruppe.
 {\it Invent. Math.}, 12:57--61, 1971.

\bibitem{brieskorn:groupes-de-tresses}
E.~Brieskorn.
 Sur les groupes de tresses [d'apres {V}. {I}. {A}rnol$'$d].
 In {\it Seminaire Bourbaki, Exp. no. 401}, number 317 in Lec. Notes
  in Math., pages 21--44. Springer, 1973.

\bibitem{brieskorn:artin-gruppen-und-coxeter-gruppen}
E.~Brieskorn and K.~Saito.
 {A}rtin-{G}ruppen und {C}oxeter-{G}ruppen.
 {\it Invent. Math.}, 17:245--271, 1972.

\bibitem{charney:artin-groups-finite-type-biautomatic}
R.~Charney.
 {A}rtin groups of finite type are biautomatic.
 {\it Math. Ann.}, 292:671--683, 1992.

\bibitem{charney:finite-k-pi-1's-for-artin-groups}
R.~Charney and M.~Davis.
 Finite ${K}(\pi,1)$'s for {A}rtin groups.
 In F.~Quinn, editor, {\it Prospects in Topology}, volume 138 of {\it   Annals of Math. Studies}, pages 110--124. Princeton University Press, 1995.

\bibitem{ATLAS}
J.~H. Conway, R.~T. Curtis, S.~P. Norton, R.~A. Parker, and R.~A. Wilson.
 {\it {A}{T}{L}{A}{S} of Finite Groups}.
 Oxford, 1985.

\bibitem{deligne:immuebles-des-groupes-de-tresses}
P.~Deligne.
 Les immuebles des groupes de tresses generalises.
 {\it Invent. Math.}, 17:273--302, 1972.

\bibitem{dbae:wd-proc-gps}
D.~B.~A. Epstein, J.~W. Cannon, S.~V.~F. Levy, M.~S. Paterson, and W.~P.
  Thurston.
 {\it Word Processing in Groups}.
 Jones and Bartlett, Boston, 1992.

\bibitem{fox-neuwirth:braid-groups}
R.~H. Fox and L.~Neuwirth.
 The braid groups.
 {\it Math. Scand.}, 10:110--126, 1962.

\bibitem{nguyen:regular-orbits-affine-weyl-groups}
{\relax Nguy{\^e}{\~n}~V.~D.}
 The fundamental groups of the spaces of regular orbits of the affine
  {W}eyl groups.
 {\it Topology}, 22:425--435, 1983.

\bibitem{salvetti:homotopy-type-of-artin-groups}
M.~Salvetti.
 The homotopy type of {A}rtin groups.
 {\it Math. Res. Lett.}, 1:565--577, 1994.

\bibitem{tom-dieck:symmetrische-brucken}
T.~tom Dieck.
 Symmetrische {B}r{\"u}cken und {K}notentheorie zu den
  {D}ynkin-{D}ia\-grammen von {T}yp ${B}$.
 {\it J. Reine Angew. Math.}, 451:71--88, 1994.

\bibitem{tom-dieck:rooted-cylinder-ribbons}
T.~tom Dieck.
 Categories of rooted cyliner ribbons and their representations.
 {\it J. Reine Angew. Math.}, 494:35--63, 1998.

\bibitem{van-der-lek:extended-artin-groups}
H.~van~der Lek.
 Extended {A}rtin groups.
 In {\it Singularities, part 2}, volume~40 of {\it Proc. Symp. Pure
  Math.}, pages 117--121. A.M.S., 1983.

\bibitem{van-der-lek:thesis}
H.~van~der Lek.
 {\it The Homotopy Type of Complex Hyperplane Arrangements}.
 PhD thesis, Katholieke Universiteit te Nijmegen, 1983.

\endreflist
\parindent=0pt
\references
\bye